\documentclass[11pt,reqno]{amsart}
\usepackage{amsmath,amssymb,geometry,color,tikz-cd}
\usepackage[all]{xy}
\usepackage{hyperref}
\usepackage[alphabetic]{amsrefs}
\usepackage{amssymb}

\usepackage{bbm}
\usepackage{enumitem}
\usepackage{upgreek}

\newcommand{\A}{\mathbb{A}}

\newcommand{\G}{\mathbb{G}}
\newcommand{\Q}{\mathbb{Q}}

\newcommand{\Z}{\mathbb{Z}}
\newcommand{\N}{\mathbb{N}}
\renewcommand{\P}{\mathbb{P}}
\newcommand{\T}{\mathbb{T}}

\newcommand{\cA}{\mathcal{A}}
\newcommand{\cB}{\mathcal{B}}

\newcommand{\cD}{\mathcal{D}}
\newcommand{\cE}{\mathcal{E}}
\newcommand{\cF}{\mathcal{F}}

\newcommand{\cH}{\mathcal{H}}

\newcommand{\cL}{\mathcal{L}}

\newcommand{\cO}{\mathcal{O}}

\newcommand{\cV}{\mathcal{V}}
\newcommand{\cZ}{\mathcal{Z}}
\newcommand{\cX}{\mathcal{X}}
\newcommand{\cY}{\mathcal{Y}}

\newcommand{\fa}{\mathfrak{a}}
\newcommand{\fb}{\mathfrak{b}}
\newcommand{\fc}{\mathfrak{c}}
\newcommand{\fm}{\mathfrak{m}}

\newcommand{\fX}{\mathfrak{X}}

\newcommand{\bP}{\mathbb{P}}
\newcommand{\bC}{\mathbb{C}}
\newcommand{\bR}{\mathbb{R}}
\newcommand{\bA}{\mathbb{A}}
\newcommand{\bQ}{\mathbb{Q}}
\newcommand{\bZ}{\mathbb{Z}}

\newcommand{\bG}{\mathbb{G}}

\newcommand{\bN}{\mathbb{N}}

\newcommand{\hvol}{\widehat{\mathrm{vol}}}
\newcommand{\tc}{\mathrm{tc}}

\newcommand{\Dtc}{\Delta_{\tc}}
\newcommand{\pr}{\mathrm{pr}}
\newcommand{\red}{\mathrm{red}}
\newcommand{\cI}{\mathcal{I}}
\newcommand{\Bl}{\mathrm{Bl}}
\newcommand{\lc}{\mathrm{lc}}
\newcommand{\ac}{\mathrm{ac}}
\newcommand{\s}{\mathrm{s}}
\newcommand{\bfH}{\mathbf{H}}

\DeclareMathOperator{\Aut}{Aut}

\DeclareMathOperator{\Exc}{Exc}
\DeclareMathOperator{\Fut}{Fut}

\DeclareMathOperator{\Hilb}{Hilb}
\DeclareMathOperator{\lct}{lct}
\DeclareMathOperator{\Spec}{Spec}
\DeclareMathOperator{\vol}{vol}
\DeclareMathOperator{\mult}{mult}

\DeclareMathOperator{\ord}{ord}

\DeclareMathOperator{\Proj}{Proj}

\DeclareMathOperator{\Supp}{Supp}
\DeclareMathOperator{\DivVal}{DivVal}
\DeclareMathOperator{\Val}{Val}

\newcommand{\SL}{\mathrm{SL}}

\newcommand{\Chow}{\mathrm{Chow}}
\newcommand{\PGL}{\mathrm{PGL}}
\newcommand{\CM}{\mathrm{CM}}
\newcommand{\Hom}{\mathrm{Hom}}
\newcommand{\bT}{\mathbb{T}}
\newcommand{\tPhi}{\widetilde{\Phi}}
\newcommand{\fI}{\mathfrak{I}}
\newcommand{\bin}{\mathbf{in}}

\newcommand{\mnorm}[1]{\lVert#1\rVert_{\rm m}}

\numberwithin{equation}{section}
\newcommand\numberthis{\addtocounter{equation}{1}\tag{\theequation}}

\newtheorem{prop} {Proposition} [section]
\newtheorem{thm}[prop] {Theorem} 

\newtheorem{lem}[prop] {Lemma}

\newtheorem{prop-def}[prop]{Proposition-Definition}

\newtheorem{conj}[prop]{Conjecture}

\newtheorem{theorem}[prop]{Theorem}
\newtheorem{lemma}[prop]{Lemma}
\newtheorem{corollary}[prop]{Corollary}

\theoremstyle{definition}

\newtheorem{defn}[prop]{Definition}

\newtheorem{remark}[prop]{Remark}

\newtheorem{definition}[prop]{Definition}

\title[Optimal destabilization of K-unstable Fano varieties]{Optimal destabilization of K-unstable Fano varieties via stability thresholds}
\author{Harold Blum}
\address{Department of Mathematics, Stony Brook University, Stony Brook, NY 11794, USA}
\email{harold.blum@stonybrook.edu}

\author{Yuchen Liu}
\address{Department of Mathematics, Northwestern University, Evanston, IL 60208, USA}
\email{yuchenl@northwestern.edu}

\author{Chuyu Zhou}
\address{\'Ecole Polytechnique F\'ed\'erale de Lausanne (EPFL), MA C3 635, Station 8, 1015 Lausanne, Switzerland}
\email{chuyu.zhou@epfl.ch}
\date{\today} 

\begin{document}

\begin{abstract}
We show that for a K-unstable Fano variety, any divisorial valuation 
computing its stability threshold induces a non-trivial special test configuration preserving  the stability threshold. 
When such a divisorial valuation exists, we show that the Fano variety degenerates to a uniquely determined twisted K-polystable Fano variety.
We also show that the stability threshold can be approximated by divisorial valuations induced by special test configurations. 
As an application of the above results and the analytic work of Datar, Sz\'ekelyhidi, and Ross, we deduce that greatest Ricci lower bounds of Fano manifolds of fixed dimension form a finite set of rational numbers. 
As a key step in the proofs, we 
adapt the process of Li and Xu producing special test configurations 
to twisted K-stability in the sense of Dervan. 
\end{abstract}

\maketitle
\setcounter{tocdepth}{1}
\tableofcontents

\section{Introduction}

A classical result of Kempf \cite{Kem78} in Geometric Invariant Theory (GIT) states that any GIT
unstable point admits a unique one-parameter subgroup which destabilizes the point most rapidly.
This one-parameter subgroup is usually called an
optimal destabilization. 
For vector bundles over smooth projective varieties, such an optimal destabilization
is given by the Harder-Narasimhan filtration. Philosophically, K-stability
is closely related to asymptotic GIT stability from their set-ups. Thus, it would be natural to expect
such phenomenon to occur for K-unstable Fano varieties as well.

In this paper, we study optimal destabilization of Fano varieties via the stability threshold.
Recall, the stability threshold (also known as the $\delta$-invariant) of a Fano variety was introduced in \cite{FO18} and characterizes K-semistability. 
Since the invariant may be expressed as an infimum of a function on the space of divisorial valuations, it is natural to view 
valuations computing the infimum as optimal destabilizers. 
Our first main theorem connects such valuations  to destabilization via  test configurations minimizing
$\frac{\Fut(\cdot)}{\mnorm{\cdot}}$ where $\Fut(\cdot)$ is the generalized Futaki invariant and $\mnorm{\cdot}$ is Dervan's  minimum norm \cite{Der16}, which also agrees with the non-Archimedean $I-J$ functional  from \cite{BHJ17}.

\begin{thm}\label{thm:main1}
 Let $(X,\Delta)$ be a log Fano pair that is not uniformly K-stable, or equivalently, $\delta(X,\Delta)\leq 1$. Then we have 
 \begin{equation}\label{eq:main1}
 \inf_{(\cX,\Delta_{\tc};\cL)}\frac{\Fut(\cX,\Delta_{\tc};\cL)}{\mnorm{\cX,\Delta_{\tc};\cL}}=\delta(X,\Delta)-1,
 \end{equation}
 where the infimum runs over all normal non-trivial test configurations $(\cX,\Delta_{\tc};\cL)$ of $(X,\Delta)$. Moreover, we have the following.
 
\begin{enumerate}
\item  If there exists a prime divisor $F$ over $X$ computing the stability threshold $\delta(X,\Delta)$, then
 $F$ is dreamy  and induces a non-trivial special test configuration $(\cX_F,\Delta_{F})$ of $(X,\Delta)$ that achieves the infimum in \eqref{eq:main1}.
\item Conversely, any normal non-trivial test configuration of $(X,\Delta)$ that achieves the infimum in \eqref{eq:main1} is special and induces a prime divisor over $X$ computing $\delta(X,\Delta)$.
\end{enumerate}
Additionally, for any special test configuration $(\cX,\Delta_{\tc})$ achieving the infimum in \eqref{eq:main1}, we have  $\delta(\cX_{0},\Delta_{\tc,0})=\delta(X,\Delta)$ and the values are rational.
\end{thm}

We call a non-trivial test configuration  of $(X,\Delta)$ achieving the infimum in
\eqref{eq:main1}  an \emph{optimal destabilization}. Its central fiber is called an \emph{optimal degeneration} of $(X,\Delta)$. 
It is expected that the assumption in Theorem \ref{thm:main1} (1) on existence of a divisorial valuation computing stability threshold always holds.

\begin{conj}[Optimal Destabilization Conjecture]\label{conj:odc}
 If $(X,\Delta)$ is a log Fano pair that is not uniformly K-stable, then there exists a prime divisor $F$ over $X$ computing the stability
 threshold $\delta(X,\Delta)$. 
\end{conj}

We provide a few remarks on Conjecture \ref{conj:odc}. 
By Theorem \ref{thm:main1}, Conjecture \ref{conj:odc} is equivalent to the existence of an optimal destabilization of $(X,\Delta)$. 
Note that Conjecture \ref{conj:odc} is part of \cite[Conjecture 1.5]{BX18}. Indeed, Theorem \ref{thm:main1} shows that the latter conjecture is implied by the former one.
Although we cannot verify Conjecture \ref{conj:odc}, we show that the stability threshold can always be approximated by divisors induced by special test configurations (see Theorem \ref{thm:specialapprox}). This can be viewed as a global analogue of \cite[Theorem 1.3]{LX16}. We refer to \cite{BLX19} for further results on valuations computing the stability threshold. 

In the above  theory of optimal destabilization for Fano varieties,
optimal destabilizations are not always unique. 
This is quite different from destabilization in GIT. 
Nevertheless, using techniques from \cite{LWX18}, we show that there exists a unique twisted K-polystable optimal degeneration. 

\begin{thm}\label{thm:kpolydegen}
Let $(X,\Delta)$ be a log Fano pair that is not uniformly K-stable. If Conjecture \ref{conj:odc} holds for $(X,\Delta)$, then there exists a unique $\delta(X,\Delta)$-twisted K-polystable optimal degeneration of $(X,\Delta)$.
\end{thm}

Twisted K-stability, which appears in the above theorem, was introduced by Dervan to characterize the existence of twisted cscK metrics \cite{Der16}.  For $\mu\in (0,1]$, $\mu$-twisted K-(semi/poly)stability is defined in terms of the positivity of  $\Fut_{1-\mu}(\cdot)=\Fut( \cdot) + (1- \mu) \mnorm{\cdot }$ (see Section \ref{sec:twistedK}). 
Using the definition of twisted K-polystability and Theorem 1.1, the degeneration in 
Theorem \ref{thm:kpolydegen} can be characterized by the property that 
all its optimal degenerations or equivalently divisorial valuations computing its stability threshold are induced by $\G_m$-actions on the pair.

Conjecture \ref{conj:odc} and Theorem \ref{thm:kpolydegen} are closely related to the analytic work of   
 Datar and Sz\'ekelyhidi \cite{DS16} and Ross and Sz\'ekelyhidi \cite{RS19}
 in which the authors construct degenerations of K-unstable Fano manifolds with twisted K\"ahler-Einstein metrics. 
 It is expected that the unique optimal degeneration in Theorem \ref{thm:kpolydegen} is related to the ones constructed analytically (see Remark \ref{rem:twisted-Kps}). 
 Furthermore, it follows from \cite{DS16, RS19} that 
Conjecture \ref{conj:odc} holds for all complex  Fano manifolds (see Theorem \ref{thm:odc-Fano-manifold}). 
 
 It has been shown in \cite{BBJ18, CRZ19} that the greatest Ricci lower bound of a Fano manifold $X$ is equal to $\min\{1,\delta(X)\}$.
 As an application of Theorem \ref{thm:main1} and the work in \cite{DS16, RS19}, we prove the finiteness of greatest Ricci lower bounds of complex Fano manifolds in any fixed dimension. Note that a different proof is provided by \cite[Theorem 1.1]{BLX19} where a more general statement on finiteness of stability thresholds is proven for log Fano pairs.

  \begin{thm}\label{thm:greatestRic}
 For any fixed positive integer $n$, the set of greatest Ricci lower bounds of 
 $n$-dimensional complex Fano manifolds is a finite subset of $(0,1]\cap \bQ$.
\end{thm}

Many of the arguments in this paper rely  on further developing Dervan's notion  of twisted K-stability \cite{Der16}.
 We prove the following characterization of twisted K-semistability in terms of the stability threshold.
 \footnote{A statement related to Theorem \ref{thm:main-twisted=delta} is proved in \cite{BoJ18} using different methods (see Remark \ref{rem:adjointK}).} The result provides  new justification for calling $\delta(X,\Delta)$ the {stability threshold}.

\begin{theorem}\label{thm:main-twisted=delta}
If $(X,\Delta)$ is a log Fano pair, then 
\[
\sup \{ \mu \in (0,1] \, \vert \, \text{$(X,\Delta)$ is $\mu$-twisted K-semistable} \} 
= \min\{ 1, \delta(X,\Delta) \} . \]
Moreover, the supremum above is a maximum.
\end{theorem}

As a key ingredient in proving Theorem \ref{thm:main1}, 
we adapt the process of Li and Xu \cite{LX14} of modifying test configurations into special ones to the  twisted setting.

\begin{thm}\label{thm:main-twisted-LX}
Let $(X,\Delta)$ be a log Fano pair. Let $\mu\in (0,1]$ be a real number. Then $(X,\Delta)$ is $\mu$-twisted K-semistable (resp. K-stable) if and only if $\Fut_{1-\mu}(\cX,\Delta_{\tc};\cL)\geq 0$ (resp. $> 0$) for any normal non-trivial special test configuration $(\cX,\Delta_{\tc};\cL)$ of $(X,\Delta)$. 
\end{thm}

\begin{remark}
There are different types of optimal destabilization of Fano varieties in the  literature. Note that an optimal destabilization in our sense minimizes the generalized Futaki invariant divided by the minimum norm of Dervan \cite{Der16} (see Theorems \ref{thm:main1},  \ref{thm:main-twisted=delta}, and Proposition \ref{prop:infFut/I-J}).

\begin{enumerate}
    \item In \cite{Don05, Sze08, CSW18, Der19, His19, Xia19} etc., an optimal destabilization of a K-unstable Fano variety refers to a test configuration which minimizes the generalized Futaki invariant divided by the $L^2$-norm. 
    \item In \cite{CSW18, DS16b}, an optimal $\bR$-degeneration of a K-unstable Fano manifold maximizing the $H$-invariant is produced by the K\"ahler-Ricci flow. 
    In \cite{HL20, BLXZ21}, the uniqueness and existence of the $\bR$-degeneration is verified using an algebraic approach.

    \item In \cite{Li18, LX17}, it was conjectured that there exists a unique valuation $v^{\rm m}$ that minimizes normalized volume functional over the affine cone $C(X,-r(K_X+\Delta))$ of a K-unstable Fano variety $(X,\Delta)$, moreover $v^{\rm m}$ is quasimonomial and has finitely generated graded algebra (known as the Stable Degeneration Conjecture).
    All but the finite generation is known by  \cite{Blu18,Xu19,XZ20}.
    Assuming this conjecture, we obtain a unique finitely generated filtration of the anti-pluricanonical ring $R(X,-r(K_X+\Delta))$, which provides an optimal destabilization of $(X,\Delta)$ over an affine toric variety in an appropriate sense.
     
\end{enumerate}
\end{remark}

This paper is organized as follows. In Section \ref{sec:prelim}, we collect preliminary materials on valuations and K-stability. In Section \ref{sec:twistedK}, we study Dervan's notion of twisted K-stability of log Fano pairs when twisting by a scalar multiple of the anti-canonical polarization. We reprove certain results in \cite{Der16} following an algebraic approach. In Section \ref{sec:twistedLX}, we generalize the process of Li and Xu producing special test configurations to the twisted setting and, hence, prove Theorem \ref{thm:main-twisted-LX}. In Section \ref{sec:approx}, we establish results of approximating the stability threshold by divisors with nice geometric properties and use them to deduce a valuative criterion for twisted K-stability. In particular, we prove Theorems \ref{thm:specialapprox} and \ref{thm:main-twisted=delta}. Section \ref{sec:proof1.1} is devoted to proving Theorem \ref{thm:main1}. We study twisted normalized volume in Section \ref{sec:twistednv}. In Section \ref{sec:twistedKps}, we prove Theorem \ref{thm:kpolydegen} by establishing $\Theta$-reductivity type results for twisted K-stability. Finally, we prove Theorem \ref{thm:greatestRic} in Section \ref{sec:greatestRic}.

\medskip

\emph{Postscript remarks.} Since the first version of this article appeared on the arXiv, there has been much progress on the study of optimal destabilizations and valuations computing stability thresholds. We list a few related works below. 
\begin{enumerate}
    \item In \cite{Zhu20}, Zhuang showed that there exists a unique minimal optimal destabilizing center for any K-unstable log Fano pair $(X,\Delta)$. Here an optimal destabilizing center means the center of a valuation computing $\delta(X,\Delta)<1$.
    \item Based on the theory of optimal destabilizations introduced in this paper, in \cite{BHLLX}, the first two authors, Halpern-Leistner, and Xu prove that Conjecture \ref{conj:odc} implies the properness of K-moduli spaces.
    \item In \cite{LXZ21}, the second author, Xu, and Zhuang prove that any valuation computing $\delta(X,\Delta)< \frac{n+1}{n}$ where $n=\dim(X)$ has a finitely generated associated graded ring. This confirms Conjecture \ref{conj:odc} in full generality and hence, combined with \cite{BHLLX}, leads to a proof of the properness of K-moduli spaces.
    \item By \cite{Fuj-Tohoku, LXZ21}, the equality \eqref{eq:main1} also holds when $\delta(X,\Delta)<\frac{n+1}{n}$, and the infimum is a minimum. In fact, by Proposition \ref{prop:Fut-I-J} we have 
    \begin{equation}\label{eq:delta<(n+1)/n}
    \Fut(\cX,\Dtc;\cL)\geq (\delta(X,\Delta)-1)\mnorm{\cX,\Dtc;\cL}
    \end{equation} if a test configuration $(\cX,\Dtc;\cL)$ is special. Thus \cite[Corollary 4.4]{Fuj-Tohoku} implies that \eqref{eq:delta<(n+1)/n} holds for all normal test configurations when $\delta(X,\Delta)\in (1, \frac{n+1}{n})$. On the other hand, \cite{LXZ21} shows that there exists a non-trivial special test configuration  achieving the equality in \eqref{eq:delta<(n+1)/n} which is induced by a divisorial valuation computing $\delta(X,\Delta)$. 
\end{enumerate}

\medskip

\noindent \emph{Acknowledgements:}
We would like to thank Daniel Halpern-Leistner, Chen Jiang, Mattias Jonsson, Chi Li, Song Sun, Chenyang Xu, and Ziquan Zhuang for many helpful discussions. We wish to thank Ved Datar, Ruadha\'i Dervan, Yujiro Kawamata, J\'anos Koll\'ar, Sam Payne, Yanir Rubinstein, G\'abor Sz\'ekelyhidi, and Feng Wang for helpful comments. The third author would like to thank his advisor Prof. Chenyang Xu for his constant support and encouragement.  Much of the work on this paper was completed while the authors enjoyed the hospitality of the MSRI, which is gratefully acknowledged. Finally, we thank the anonymous referees whose comments improved this paper.

HB was partially supported by NSF grant DMS-1803102; 
YL was partially supported by the Della Pietra Endowed Postdoctoral Fellowship of the MSRI (NSF grant DMS-1440140).

\section{Notions and Preliminaries}\label{sec:prelim}

Throughout, we work over an algebraically closed characteristic zero field $k$, and specifically over $\bC$ in Section \ref{sec:greatestRic}. 
We follow the standard terminology in \cite{KM98}.

A \emph{pair} $(X,\Delta)$ is composed of a normal variety $X$ and an effective $\Q$ divisor $\Delta$ 
such that $K_X+\Delta$ is $\Q$-Cartier. For the definitions of \emph{klt}, \emph{lc}, and \emph{plt} pairs, see \cite{KM98}. A klt pair $(X,\Delta)$ is \emph{log Fano} if $X$ is projective and $-K_X-\Delta$ is ample. A variety $X$ is \emph{$\Q$-Fano} if $(X,0)$ is log Fano and a \emph{Fano manifold} if $X$ is in addition smooth.

\subsection{Valuations and associated invariants}

\subsubsection{Divisorial valuations}
Let $X$ be a normal variety.
If $\mu: Y \rightarrow X$ is a proper birational morphism
with $Y$ normal and $F\subset Y$ is a prime divisor (called a \emph{prime divisor over} $X$), then $F$ defines a valuation $\ord_F:K(X)^* \rightarrow \mathbb{Z} $ given by order of vanishing at the generic point of $F$.

Any valuation  $K(X)^\times \to \Q$ of the form $v= c \ord_F$, where $c \in \Q_{>0}$ and $F$ is a prime divisor over $X$, will be a called \emph{divisorial}. The \emph{center} of $v$, denoted $c_X(v)$, is defined to be the generic point of $\mu(F)$. 
We write $\DivVal_X$ for the set of divisorial valuations on $X$.

If $D$ is a $\Q$-Cartier $\Q$-divisor on $X$ and $v=c \ord_F\in \DivVal_X$, 
we set ${v(D):= c \cdot {\rm coeff}_F( \mu^*(D))}$. Note that if $D$ is a Cartier divisor, this equals $v(f)$ where $f$ is a local equation of $D$ at $c_X(v)$. If $\fa\subseteq \cO_X$ is a non-zero ideal, we set $v(\fa): = \min \{ v(f) \, \vert \, f\in \fa\cdot \cO_{X,c_X(v)} \}$.

\subsubsection{Log discrepancy and log canonical thresholds}
Let $(X,\Delta)$ be a pair. 
For $v=c\ord_F \in \DivVal_X$ as above, we call
$$A_{X,\Delta}(v) : =c(1 + {\rm coeff}_F(K_{Y}-\mu^*(K_X+\Delta ))$$
the \emph{log discrepancy} of $v$. When $c=1$, we simply write $A_{X,\Delta}(F)$ for this value. 
More generally, if $\fa\subseteq \cO_X$ is a non-zero ideal 
and $b\in \Q_{\geq0}$, we set  $A_{X,\Delta+ b\cdot \fa}(v):= A_{X,\Delta}(v)-b v(\fa)$.

If $(X,\Delta)$ is a klt pair and $D$ an effective $\Q$-divisor, then 
the \emph{log canonical threshold} of $D$ is given by
$\lct(X,\Delta;D) : = \sup\{ c\in \Q_{\geq 0}\, \vert\, (X,\Delta+cD)
\text{ is lc}\}$ and equals $\min_F \frac{A_{X,\Delta}(F)}{\ord_{F}(D)}$. See \cite[Section 8]{Kol97} for further details.

\subsubsection{Log Fano pairs} 
Let $(X,\Delta)$ be a log Fano pair and $F$ a prime divisor over $X$, where $F$ arises on proper normal model $\mu:Y\to X$. 
We will be interested in  how anti-canonical divisors of $(X,\Delta)$ vanish along $F$. 

Fix a positive integer $r$ so that $r(K_X+\Delta)$ is Cartier.
Following \cite[Def. 1.3.1]{Fuj16}, we say 
$F$ is \emph{dreamy} if the $\N \times \Z$-graded algebra
\[
\bigoplus_{m\in \N} \bigoplus_{p \in \Z}
H^0(Y, -mr\mu^*(K_X+\Delta) - p F)
\]
is finitely generated over $k$; this is independent of the choice of $r$. By \cite[Corollary 1.3.1]{BCHM10}, $F$ is dreamy if $Y$ is Fano type (i.e. there exists a $\Q$-divisor $B$ on $Y$ such that $(Y,B)$ is klt and $-K_{Y}-B$ is big and nef), since a small $\bQ$-factorial modification of $Y$ is a Mori dream space.

We consider the invariants
\[
S_{X,\Delta}(F): = 
\frac{1}{ \vol(-K_X-\Delta) }
\int_0^\infty \vol(- \mu^*(K_X+\Delta)-tF)\, dt
\]
and 
\[
T_{X,\Delta}(F): = 
\sup \{ t \, \vert \,  - \mu^*(K_X+\Delta)-tF \text{ is pseudoeffective}\}
.\]
The invariant $S_{X,\Delta}(F)$ plays  a key role in the work K. Fujita, C. Li, and others on K-stability. 

If $v= c\ord_{F} \in \DivVal_X$, where $c\in \Q_{>0}$, we set $S_{X,\Delta}(v)= c \cdot  S_{X,\Delta}(F)$ and $T_{X,\Delta}(v) = c\cdot T_{X,\Delta}(F)$.
When the pair $(X,\Delta)$ is  clear from context, we simply write $S(v)$ and $T(v)$ for these values.

The values $S$ and $T$ may be written in terms of the vanishing of certain classes of anti-canonical divisors along $F$. 
We write  $$|-K_X-\Delta|_{\Q}:= \{ \textrm{$\bQ$-divisor } D\geq 0 \, \vert \, D\sim_{\Q}-K_X-\Delta \}$$
 for the $\Q$-linear system of $-K_X-\Delta$. 
 Following \cite[Def. 0.1]{FO18}, we say $D\in |-K_X-\Delta|_{\Q}$
is of $m$-basis type if there exists a basis $\{s_1,\ldots, s_{N_m}\}$ for $H^0(X,\cO_X(-m(K_X+\Delta))$
such that 
$
D= 
\frac{1}{mN_m}
\left(
\{s_1=0\} + \cdots + \{s_{N_m}=0\} 
\right).$

With the above notation, for each positive integer $m$ divisible by $r$, we set 
\[
S_m(v) : = \max \{ v(D) \, \vert \, D\in |-K_X-\Delta|_{\Q} \text{ is $m$-basis type} \}
\]
and \[
T_m(v): = \max \{ v(D) \, \vert \, 
\in  \tfrac{1}{m}|-m(K_X+\Delta)|\}
.\]
We have equalities
$S(v)= \lim_{m\to \infty} S_{mr}(v)$ and 
$T(v)  = \sup_m T_{mr}(v) $. See  \cite{FO18,BJ17} for a proof of the first equality, which requires a brief argument.

\subsubsection{Stability threshold}
Let  $(X,\Delta)$ be a log Fano pair and $r$ a positive integer so that $r(K_X+\Delta)$ is a Cartier divisor. 
For each positive integer $m$ divisible by $r$, set
\[
\delta_m(X,\Delta): = \inf \{ \lct(X,\Delta;D) \, \vert \, D
\text{ is $m$-basis type} \}.\footnote{This infimum is indeed a minimum, since $\lct$ is a constructible function in $\bQ$-Gorenstein families (see \cite[Corollary 1.10]{Amb16}), and $m$-basis type divisors are parametrized by a scheme of finite type.}
\]
We call $\delta(X,\Delta):= \limsup_{m \to \infty} \delta_{mr}(X,\Delta)$ the \emph{stability threshold} of $(X,\Delta)$. 

This invariant was introduced by Fujita and Odaka in  \cite{FO18} and coincides with an invariant suggested by Berman \cite[Section 2.9]{BoJ18}. 
As shown in  \cite{BJ17}, 
\begin{equation}\label{eq:delta-infvals}
\delta(X,\Delta) = \inf_{v} \frac{A_{X,\Delta}(v)}{S(v)}
\quad\text{ and } \quad
\delta_{mr}(X,\Delta)= \inf_{v} \frac{A_{X,\Delta}(v)}{S_{mr}(v)},
\end{equation}
where the infimum runs through $v\in \DivVal_X$,  and $\lim_{m\to \infty} \delta_{mr}(X,\Delta)$ exists. 
In light of \eqref{eq:delta-infvals}, we say $v$ \emph{computes}  the stability threshold if $\delta(X,\Delta)= \frac{A_{X,\Delta}(v)}{S(v)}$.

\subsection{K-stability}

In this section, we recall definitions and results on K-stability. See \cite{Tian97,Don02,LX14,Der16,BHJ17}.

\subsubsection{Test configurations}

Let $(X,\Delta)$ be an $n$-dimensional log Fano pair. 

\begin{defn}
A \emph{test configuration} $(\cX,\Delta_{\tc}; \cL)$ of $(X,\Delta)$ consists of the following data:
\begin{enumerate}
    \item a flat proper morphism of schemes $\pi: \cX \to \mathbb{A}^1$,
    \item a $\mathbb{G}_m$-action on $\cX$ extending the standard action on $\mathbb{A}^1$,
    \item a $\pi$-ample $\mathbb{G}_m$-equivariant line bundle $\cL$ on $\cX$,
    \item an isomorphism $(\cX_1, \cL_1) \simeq (X,L)$ where $L$ is a line bundle  on $X$ in the $\Q$-linear equivalence class of $-r(K_X+\Delta)$, where $r\in \Q_{>0}$
    \end{enumerate}
   Via (4), there is a canonical equivariant isomorphism $\cX \vert_{\A^1 \setminus 0} \simeq X \times (\mathbb{A}^1\setminus 0) $, where $\G_m$-acts trivially on $X$ and in the standard way on $\A^1\setminus 0$. The divisor $\Delta_{\tc}$ is defined as the closure of $\Delta\times  ({\A^1 \setminus 0})$ under the above embedding $X \times (\mathbb{A}^1\setminus 0) \hookrightarrow \cX$. 
   If the line bundle $\mathcal{L}$ is only $\pi$-semiample (rather than $\pi$-ample), we call $(\mathcal{X},\Delta_{\tc};\mathcal{L})$ a \emph{semiample test configuration}.

  A test configuration has a compactification $( \overline{\cX},\overline{\Delta}_{\tc};  \overline{\cL})/ \mathbb{P}^1$ defined via gluing $(\cX,\Delta_{\tc}; \cL)$ and $(X \times (\P^1 \setminus 0), \Delta\times (\mathbb{P}^1 \setminus 0); \pr_1^*L )$ along their respective open sets $\cX \setminus \cX_0$ and $X \times (\mathbb{A}^1 \setminus 0)$.

A test configuration is \emph{normal} if $\cX$ is normal. It is \emph{special} if $\cL\sim_{\bQ} -r(K_{\cX}+\Delta_{\tc})$ and $(\cX,\Dtc+\cX_0)$ is plt, which implies that the central fiber $\cX_0$ is integral. \footnote{
If  $(\cX,\Delta_{\tc}; \cL)$ is a special test configuration, we will often omit $\cL$ from the notation. }
It is of \emph{product type} if $(\cX,\Delta_{\tc})$ is isomorphic to $(X ,\Delta)\times \mathbb{A}^1$ over $\A^1$,
and \emph{trivial} if the latter isomorphism is $\mathbb{G}_m$-equivariant. 
\end{defn}

\subsubsection{Invariants of test configurations}
Let $(\cX,\Dtc; \cL)$ be a normal test configuration of $(X,\Delta)$ and set $V: = (-K_X-\Delta)^n$. 
The \emph{generalized Futaki invariant}
of $(\cX,\Delta_{\tc};\cL)$ is given by 
\[
\Fut( \cX, \Dtc;\cL) : =  r^{-n} V^{-1} \left( 
\frac{n}{n+1} r^{-1} \overline{\cL}^{n+1} + 
\overline{\cL}^{n} \cdot (K_{\overline{\cX}/ \P^1} + \overline{ \Delta}_\tc ) \right)
.\]
Note that $\Fut(\cX,\Dtc;\cL)=\Fut(\cX,\Dtc;\cL^{\otimes m})$ for any $m\in \bZ_{>0}$.
When $(\cX,\Dtc; \cL)$ is not normal, one can define $\Fut(\cX,\Dtc;\cL)$ in terms of weights,  rather than intersection numbers. Since the definition will only be used briefly in Section \ref{sec:greatestRic}, we refer the readers to \cite{Don02, Li15, OS15}.

Next, we define Dervan's minimum norm \cite{Der16}, which also agrees with the non-Archimedean $I-J$ functional  by \cite[Remark 7.12]{BHJ17}.
Consider the diagram 
\begin{equation}\label{eq:tc-resolution}
 \begin{tikzcd}[row sep= 1.25 em]
   & \cY   \arrow[rd,"\tau"] \arrow[ld,swap,"\rho"]  &  \\
   X_{\P^1} \arrow[rr,dashed] & & \overline{\cX}
        \end{tikzcd}
        \end{equation}
where $X_{\P^1}:= X\times \P^1$, $L_{\P^1}:=\pr_1^* L$, and  $\cY$ is the normalization of the graph of $X_{\P^1} \dashrightarrow \cX$. 
The \emph{minimum norm}
of $(\cX,\Dtc; \cL)$
is given by
\begin{align*}
\mnorm{\cX, \Dtc;\cL}  :&=  r^{-n-1}V^{-1} \Big( \frac{1}{n+1} \overline{\cL}^{n+1} -  \ \left( \tau ^*\overline{\cL} - \rho^*   L_{\P^1}  \right)    \cdot (\tau^*\overline{\cL})^{n}    \Big),\\
&=
r^{-n-1}V^{-1} \Big(  -\frac{n}{n+1} \overline{\cL}^{n+1} +  \rho^*   L_{\P^1}      \cdot (\tau^*\overline{\cL})^{n}    \Big).
\end{align*}
The above value is positive when the test configuration is non-trivial \cite{BHJ17,Der16} and zero otherwise.

\subsubsection{Definition of K-stability}
\begin{defn}[\cite{Tian97,Don02,LX14}]\label{def:kstable} A log Fano pair $(X,\Delta)$ is
\begin{enumerate}
    \item  \emph{K-semistable} (resp., \emph{K-stable})  if $\Fut(\cX,\Dtc; \cL) \geq 0 $ (resp., $>0$) for all normal non-trivial  test configurations $(\cX,\Dtc; \cL)$ of $(X,\Delta)$;
    \item \emph{K-polystable} if it is \emph{K-semistable} and any normal test configuration of $(X,\Delta)$ satisfying $\Fut(\cX,\Dtc; \cL) = 0 $  is of product type. 
\end{enumerate}
By \cite{LX14}, it is suffices to consider only special test configurations  in the above definition.
\end{defn}

We will also briefly use the notion of uniform K-stability in this paper. 

\begin{defn}[\cite{Der16,BHJ17}]
A log Fano pair $(X,\Delta)$ is uniformly K-stable, if there exists $\varepsilon >0$ so that $
\Fut(\cX,\Delta_{\tc};\cL ) \geq \varepsilon \mnorm{\cX,\Delta_{\tc};\cL }$
for all normal test configurations $(\cX,\Delta_{\tc};\cL)$ of $(X,\Delta)$. 
\end{defn}

\subsubsection{Valuative criterion}
The following valuative characterization of K-stability proved by K. Fujita and C. Li
will play an important role in this paper.

\begin{thm}[\cite{Fuj16,Li17}]
 A log Fano pair $(X,\Delta)$ is
 K-semistable (resp., uniformly K-stable) if and only if
 \[
 \inf_{F} \frac{A_{X,\Delta}(F)}{S(F)}  \geq 1 \, \, \text{ (resp., $>1$)}
 \] where the infimum runs through all prime divisors $F$ over $X$.
\end{thm}
 In Section \ref{sec:twistedkalg}, we will generalize the above valuative criterion to the twisted setting.

\subsection{Test configurations and valuations}
In this subsection, we discuss the relation between test configurations and divisorial valuations \cite[Section 4.1-4.2]{BHJ17}.

\subsubsection{Correspondence}
Let $(\cX,\Delta_{\tc}; \cL)$ be a normal test configuration of a log Fano pair $(X,\Delta)$.
An irreducible component $E$ of $\cX_0$ defines a $\mathbb{G}_m$-equivariant valuation $\ord_E$ of $K(\cX)\simeq K(X\times\bA^1)$. 
Following \cite[\S 4.1-4.2]{BHJ17},
we write $r(\ord_E)$  for the valuation of $K(X)$
defined via restricting $\ord_E$ via the canonical isomorphism
$ K(X \times \mathbb{A}^1 )\simeq K(X)(t)$.

We say $E$ is non-trivial if it is not the birational transform of $X\times 0 $ 
via the map $X_{\P^1} \dashrightarrow \cX$; this is equivalent to the condition that $r(\ord_E)$ is not the trivial valuation.
If $E$ is non-trivial, we set $v_E = b_{E}^{-1} r (\ord_E)$ where $b_E = {\rm coeff}_E (\cX_0)$.
By \cite[Lemma 4.1]{BHJ17}, $v_E\in {\rm DivVal_X}$.

\subsubsection{Test configurations with integral central fiber} 
The following statements were observed and utilized in the work of K. Fujita and C. Li. 

\begin{prop} \label{prop:Fut-I-J}
Let $(\cX,\Delta_{\tc}; \cL)$ be a nontrivial test configuration of  $(X,\Delta)$. 
If $\cX_0$ is integral, then $\cX$ is normal
and $v_{\cX_0}$ is dreamy divisorial valuation.
Furthermore, 
\[
\Fut(\cX, \Dtc; \cL) = A_{X,\Delta}(v_{\cX_0}) - S_{X,\Delta}(v_{\cX_0}), \quad \text{ and } \quad \mnorm{ \cX, \Dtc; \cL } = S_{X, \Delta}(v_{\cX_0}). 
\]
\end{prop} 

\begin{proof} 
Since  the central fiber is integral, $\cX$ is normal \cite[Proposition 2.6.iv]{BHJ17}
and, after replacing $\cL$ with $\cL+d\cX_0$ for some integer $d$, we may assume
$\overline{\cL} \sim_{\Q}  -r (K_{\overline{\cX}/\P^1}+ \overline{\Delta}_{\tc})$  \cite[Lemma 2.6]{Fuj17b}. In this case, it is shown in  \cite[Proof of Theorem 5.1]{Fuj16} that $v_{\cX_0}$ is dreamy,
\begin{equation}\label{eq:FutakiL^n+1}
\Fut(\cX,\Dtc;\cL)= - r^{-n-1} V^{-1} \frac{1}{(n+1)} \overline{\cL}^{n+1}   =
A_{X,\Delta}(v_{\cX_0}) - S_{X,\Delta}(v_{\cX_0}) \end{equation}
and
$$B: = \tau^* \big( -K_{\overline{\cX} / \P^1}  - \overline{\Delta}_{\tc} \big)-  \rho^* \pr_1^* \left(-K_X-\Delta\right)  $$
is supported on $\cY_0$
and has coefficient $- A_{X,\Delta}(v_{\cX_0})$ along the birational transform of $\cX_0$. 
Note that \eqref{eq:FutakiL^n+1} first appeared in  \cite[Lemma 6.20 and Equation (84)]{Li17}.
We compute
\begin{equation}\label{eq:B.L^n}
\left( \tau^*\overline{\cL}- \left( \rho^*   \pr_1^* L  \right)    \right)  \cdot \tau^* \overline{\cL}^{n}  
 = 
 rB \cdot \tau^*\overline{\cL}^{n}   
 =
 -r A_{X,\Delta}(v_{\cX_0}) \cX_0 \cdot \overline{\cL}^{n}  
=
- r^{n+1} V A_{X,\Delta}(v_{\cX_0})
,\end{equation}
where the second equality follows from the projection formula 
and the third from the fact that $\cX_0$ and $\cX_1$ are algebraically equivalent.
Combining \eqref{eq:FutakiL^n+1} and \eqref{eq:B.L^n}, we conclude 
$
\mnorm{\cX,\Dtc;\cL} =S_{X,\Delta}(v_{\cX_0})$.
\end{proof} 

In the reverse direction, dreamy divisors induce test configurations with integral central fiber. 

\begin{lemma}
If $F$ is a dreamy divisor over a log Fano pair $(X,\Delta)$
and $c\in \Z_{>0}$, 
then there is a normal test configuration $(\cX,\Delta_{\tc};\cL)$ 
such that $\cX_0$ is integral and $v_{\cX_0} = c \ord_F$. 
\end{lemma}

\begin{proof}
The test configuration arises from taking the Proj of the finitely generated bigraded algebra associated to $F$ and $c$. 
See \cite[Section 3.2]{Fuj17b} for details. 
\end{proof}

\begin{defn}\label{d:specialdiv}
 A prime divisor $F$ over a log Fano pair $(X,\Delta)$ is called \emph{special} if there exists a non-trivial special test configuration $(\cX,\Delta_{\tc})$ of $(X,\Delta)$ and $c\in \Q_{>0}$ such that $v_{\cX_0}=c\ord_F$.
 Note that special divisors are necessarily dreamy by Proposition \ref{prop:Fut-I-J}.
\end{defn}

\subsubsection{Extension of boundary divisors}\label{ss:extensionboundary}
Let $(\cX,\Delta_{\tc};\cL)$ be a normal test configuration of a log Fano pair $(X,\Delta)$.
Recall the diagram \eqref{eq:tc-resolution}  \[
 \begin{tikzcd}[row sep= 1.2 em]
   & \cY   \arrow[rd,"\tau"] \arrow[ld,swap,"\rho"]  &  \\
   X_{\P^1} \arrow[rr,dashed] & & \overline{\cX}
        \end{tikzcd}
        \]
where $\cY$ is the normalization of the graph of $X_{\P^1}:=X\times \P^1 \dashrightarrow \overline{\cX}$.

Given $D \in |-K_X-\Delta|_{\Q}$,
we write $\overline{\cD}$ for the $\Q$-divisor on $\overline{\cX}$ that is the componentwise
closure of $D_{ \P^1 \setminus 0}:= D\times (\P^1 \setminus0)$
under the open immersion  $X_{\P^1 \setminus 0} \hookrightarrow \overline{\cX}$. Let $D_{\P^1}:=\pr_1^*D$ as a $\bQ$-Cartier $\bQ$-divisor on $X_{\bP^1}$.
The following elementary computations will play a key role later. 

\begin{lemma}\label{l-compat}
With the above notation, 
\begin{enumerate}
    \item $
\tau_* \left( \rho^* D_{\P^1}\right)  = \overline{\cD} + \sum_{E} b_E v_E(D)E
$ and
    \item 
    $\overline{\cD} \cdot \overline{\cL} ^n = r^{-1}(\rho^*L_{\P^1}) \cdot \left(\tau^*\overline{\cL}\right)^n -  \sum_{E} b_E v_E(D) E \cdot \overline{\cL}^n,$
\end{enumerate}
where the sums run through irreducible components $E$ of $\cX_0$.
\end{lemma}
\begin{proof}
Observe that
$
\rho^* D_{\P^1} = \widetilde{ D}_{\P^1} +{ \sum_{E} \ord_E(D_{\P^1}) E}
$,
 where
$ \widetilde{ D}_{\P^1}$ is the strict transformation of $D_{\bP^1}$ on $\cY$ and the sum runs through irreducible components of $\cY_0$. 
Using \cite[Section 3.1]{BHJ17}, we see
\[
\ord_E(D_{\P^1}) = 
r(\ord_E)(D) 
= 
b_E v_E(D) . \]
Since $\tau_*\widetilde{D}_{\P^1} = \overline{\cD}$, (1) holds.

To prove (2), we apply (1) and the projection formula to see 
\begin{align*}
\rho^* D_{\P^1}\cdot (\tau^*\overline{\cL})^n
&= 
\tau_* \left(\rho^* D_{\P^1} \right) \cdot \overline{\cL}^n\\
& =
  \overline{\cD}\cdot \overline{\cL}^n + \sum_{E } b_E v_E(D)\, E
  \cdot \overline{\cL}^n.
  \end{align*}
Since $rD\sim_{\Q}L$, we know 
$r\rho^*D_{\P^1} \sim_{\Q} \rho^* L_{\P^1}$ and the result follows. 
\end{proof}

\subsubsection{Compatible divisors}
The previous computations motivate the next definition.
Below, we continue to use notation from Section \ref{ss:extensionboundary}.

\begin{definition}
A $\Q$-divisor $D \in |-K_X-\Delta|_{\Q}$ is \emph{compatible} with $\cX$
if ${c_X( v_E) \not \subseteq \Supp(D)}$ for all non-trivial irreducible components $E $ of $ \cX_0$. 
Since the condition ${c_X( v_E) \not \subseteq \Supp(D)}$  
is equivalent to  $v_E(D) = 0$,
Lemma \ref{l-compat} (1) implies $D$ is compatible with $\cX$ if and only if $\tau_* \rho^*(D_{\P^1}) = \overline{\cD}$. 
\end{definition}

We proceed to state a few elementary properties of compatible divisors. 

\begin{lemma}\label{l-compatlin}
If $D,D' \in |-K_X-\Delta|_{\Q}$ are compatible with $\cX$, then $\overline{\cD} \sim_{\Q} \overline{\cD'}$. 
\end{lemma}
\begin{proof}
Since $D \sim_{\Q} D'$, $\tau_* \left( \rho^* D_{\P^1} \right) \sim_{\Q}\tau_* \left( \rho^* D'_{\P^1} \right)$. Applying Lemma \ref{l-compat} (1) finishes the proof.
\end{proof}

\begin{lemma}\label{l-compatcontract}
Let $(\cX',\Delta_{\tc}';\cL')$ be a normal test configuration of $(X,\Delta)$ such that the natural map $\cX \dashrightarrow \cX'$ is a birational contraction. 
If $D\in |-K_{X}-\Delta|_{\Q}$ is compatible with $\cX$, then $D$ is compatible with $\cX'$.
\end{lemma}

\begin{proof}
Since $\cX \dashrightarrow\cX'$ is a birational contraction,
any non-trivial component $E'$ of $\cX'_0$ is the birational transform of a non-trivial component $E $ of $ \cX_0$. Since $\ord_E$ and $\ord_{E'}$ induce the same valuation of $K(X\times \P^1)$, $v_E=v_{E'}$. Therefore, compatibility with $\cX$ implies compatibility with $\cX'$.
\end{proof}

\begin{lemma}\label{l-compatgeneral}
Let $m$ be a positive integer such that $-m(K_X+\Delta)$ is a Cartier divisor and  $|-m(K_X+\Delta)|$ is base point free.
\begin{enumerate}
    \item If $D \in m^{-1}|-m(K_X+\Delta)|$ is general, then $ D$ is compatible with $\cX$.
    \item If $D_1,D_2 \in m^{-1}|-m(K_X+\Delta)| $ are general, then $\overline{\cD}_1\sim_{\Q}\overline{\cD}_2$.
\end{enumerate}
\end{lemma}
\begin{proof}
Observe that  a general element $|-m(K_X+\Delta)|$ 
does not intersect $c_{X}(v_E)$ for all non-trivial  components  $E$ of  $\cX_0$. Therefore, a general divisor $D$ is compatible with $\cX$. Statement (2) now follows from  Lemma \ref{l-compatlin}.
\end{proof}

\section{Twisted K-stability via an algebraic approach}\label{sec:twistedkalg}
The notion of twisted K-stability was introduced in the work of Dervan 
to provide an algebraic criterion for the existence of 
twisted constant scalar curvature metrics \cite{Der16}.
In this section, we study twisted K-stability in the setting of log Fano pairs and when
twisting by a scalar multiple of the anti-canonical polarization. 

We will first use an algebraic approach to reprove some properties of twisted K-stability proved in \cite{Der16}.
We will then use this framework to  develop the process due to Li and Xu \cite{LX14}
of modifying test configurations in the twisted setting.

\subsection{Twisted K-stability}\label{sec:twistedK}

Throughout this section, let $(X,\Delta)$ be an $n$-dimensional log Fano pair.


\subsubsection{Definition}

\begin{definition}[Twisted generalized Futaki invariant]\label{def:twistedFut}
For a real number $\mu\in (0,1]$ 
and a test configuration  
$(\cX,\Delta_{\tc};\cL)$ of $(X,\Delta;L)$, the \emph{$\mu$-twisted generalized Futaki invariant} is defined as
\[
\Fut_{1-\mu} (\cX,\Delta_{\tc}; \cL) :=  \sup \big\{  \Fut( \cX,\Delta_{\tc}+(1-\mu) \cD;\cL) \,\, \vert \, \, D \in |-K_X-\Delta|_{\Q} \big\}.
\]
Here, $\cD$ denotes the $\Q$-divisor on $\cX$ that is the componentwise closure of $D_{ \A^1 \setminus 0}$ 
under the open immersion $X_{\A^1 \setminus 0}
\hookrightarrow\cX$.
When writing $\Fut( \cX,\Delta_{\tc}+(1-\mu) \cD;\cL)$, we are viewing 
$(\cX,\Delta_{\tc}+(1-\mu) \cD;\cL)$ as a test configuration 
of the log Fano pair $(X, \Delta + (1-\mu)D)$. 
\end{definition}

\begin{definition}[Twisted K-stability]\label{def-twistedKs}
A log Fano pair $(X,\Delta)$ is
\begin{enumerate}
    \item \emph{$\mu$-twisted K-semistable} (resp., \emph{$\mu$-twisted K-stable}) if $\Fut_{1-\mu}( \cX,\Delta_{\tc};\cL)\geq 0$ (resp., $>0$) 
for all normal non-trivial test configurations 
$(\cX,\Delta_{\tc};\cL)$ of $(X,\Delta)$. 
    \item \emph{$\mu$-twisted K-polystable} if it is $\mu$-twisted K-semistable
    and any normal test configuration of $(X,\Delta)$ satisfying $\Fut_{1-\mu}(\cX, \Dtc ; \cL)=0$ is product type.
    \end{enumerate} 
    
\end{definition}

\begin{remark}\label{rem:logimpliestwisted}

The previous definition immediately implies the following:
\begin{enumerate}
    \item  $1$-twisted K-(semi/poly)stability agrees with K-(semi/poly)stability. 
    \item  If there exists $D\in |-K_X-\Delta|_{\Q}$ such that 
$(X,\Delta+(1-\mu)D)$ is K-semistable (resp., K-stable), then $(X,\Delta)$
is $\mu$-twisted K-semistable (resp.,  $\mu$-twisted K-stable).
\end{enumerate}
\end{remark}

\begin{remark}
The above conventions differ significantly from those in \cite{Der16}.
When $\Delta=0$, the $\mu$-twisted K-stability of $(X,\Delta)$ as defined above  corresponds to the twisted K-stability of the triple $(X,-K_X,-\tfrac{1-\mu}{2}K_X)$ as defined in \cite[Definition 2.7]{Der16}. 
\end{remark}

The following theorem relates twisted K-semistability to the stability threshold. Indeed, it shows that for a log Fano pair $(X,\Delta)$ that is not uniformly K-stable, its stability threshold $\delta(X,\Delta)$ can also be viewed as its ``twisted K-semistability threshold''. 
We will present its proof in Section \ref{sec:optimal}. See \cite[Theorem C]{BL18b} for a related statement.

\begin{theorem}[=Theorem \ref{thm:main-twisted=delta}]\label{thm:twisted=delta}
If $(X,\Delta)$ is a log Fano pair, then 
\begin{equation}\label{eq:sup=delta}
\sup \{ \mu \in (0,1] \, \vert \, \text{$(X,\Delta)$ is $\mu$-twisted K-semistable} \} 
= \min\{ 1, \delta(X,\Delta) \} . 
\end{equation}
Moreover, the supremum above is a maximum.
\end{theorem}

When $X$ is smooth and $\Delta=0$, the result is known to follow from analytic arguments in \cite{BBJ18,CRZ19,DS16}.

\begin{remark}\label{rem:adjointK}
The above notion of twisted K-stability is closely related 
to adjoint K-stability defined in \cite{BoJ18} (in fact, the  notions are expected to coincide \cite[Conjecture 2.5]{BoJ18}). An adjoint version of Theorem \ref{thm:twisted=delta} is proved in \cite[Section 5]{BoJ18} using different methods. 
\end{remark}

\begin{remark}\label{rem:adjointD}
Another important invariant of a test configuration is the Ding invariant,
which was  introduced by Berman in \cite{Ber16} and plays a key role in \cite{BBJ18}.
It is natural to define the twisted Ding invariant and twisted Ding (semi)stability by replacing all mentions of $\Fut$ with ${\rm Ding}$ in Definitions \ref{def:twistedFut} and \ref{def-twistedKs}. 
This definition is closely related to adjoint Ding stability as defined in \cite{BoJ18}. 
We can show that  $\mu$-twisted K-semistability is equivalent to $\mu$-twisted Ding semistability. 
The case when $\mu=1$ holds by \cite{BBJ18, Fuj16}. In general, since ${\rm Ding}(\cdot)\leq \Fut(\cdot)$, $\mu$-twisted Ding semistability implies $\mu$-twisted K-semistability. On the other hand, by \cite[Theorem 7.2 (2)]{BL18b}, for every $\mu\in (0, \min\{1,\delta(X,\Delta)\})\cap\bQ$ there exists $D\in |-K_X-\Delta|_{\bQ}$ such that $(X, \Delta+ (1-\mu) D)$ is K-semistable and hence Ding-semistable by \cite{BBJ18, Fuj16}.  Thus by definition we know that $(X,\Delta)$ is $\mu$-twisted Ding-semistable whenever $\mu\in (0, \min\{1,\delta(X,\Delta)\}]$. Hence Theorem \ref{thm:twisted=delta} shows that $\mu$-twisted K-semistability implies $\mu$-twisted Ding semistability.

\end{remark}

\subsubsection{Properties of the twisted generalized Futaki invariant}

\begin{prop}\label{prop:twistFut}
Let $(\cX,\Delta_{\tc};\cL)$ be a normal test configuration of $(X,\Delta)$. Let $\mu\in (0,1]$ be a real number. The following hold:
\begin{enumerate}
    \item $\Fut_{1-\mu}(\cX, \Dtc;\cL) = \Fut(\cX, \Dtc;\cL) + (1-\mu) \mnorm{ \cX, \Dtc;\cL}$;
    \item If $D\in |-K_X-\Delta|_\mathbb{Q}$ is compatible with $\cX$, then 
    $$\Fut_{1-\mu}(\cX, \Dtc;\cL) = \Fut(\cX, \Dtc + (1-\mu) \cD ;\cL).$$
\end{enumerate}
\end{prop}

Related statements appear in the proofs of Theorem 2.10 (i) and Proposition 3.6 in \cite{Der16}.

\begin{proof}
Fix $D \in |-K_X-\Delta|_{\Q}$. Set $V:= 
(-K_X-\Delta)^n$. Let $r$ be the positive rational number such that $L \sim_{\Q}-r(K_X+\Delta)$. 
Since 
\[
(-K_X-\Delta-(1-\mu)D)^n = \mu^n V\quad \text{ and } \quad L \sim_{\Q} r\mu^{-1}(-K_X-\Delta-(1-\mu)D),\] we see  
\[
\Fut(\cX, \Dtc + (1-\mu) \cD ;\cL)
 = 
r^{-n}V^{-1} \left(
\frac{n}{n+1}r^{-1}\mu \overline{\cL}^{n+1} + \overline{\cL}^n  \cdot ( K_{\overline{\cX}/\mathbb{P}^1}+ \overline{\Delta}_{\tc}+(1-\mu) \overline{\cD} ) \right) 
\]
Using the above expression and then  Proposition \ref{l-compat} (2), 
we see 
\begin{align*}
\Fut(\cX, \Dtc + (1-\mu) \cD ;\cL)-
\Fut(\cX, \Dtc ;\cL) &=
(1-\mu)r^{-n}V^{-1}
\Big( 
 -\frac{n}{n+1} r^{-1}\overline{\cL}^{n+1}
 +  \overline{\cD} \cdot \overline{\cL}^n\Big)\\
 &= (1-\mu) \mnorm{ \cX, \Dtc;\cL} - r^{-n}V^{-1}\sum_{E} b_E v_E(D) 
 E \cdot \overline{\cL}^n.
 \end{align*}
This implies
\[
\Fut(\cX, \Dtc + (1-\mu) \cD ;\cL) 
\leq 
\Fut(\cX, \Dtc ;\cL) 
+
(1-\mu ) \mnorm{ \cX, \Dtc;\cL}\]
and equality holds when $v_E(D)=0$ for all irreducible components $E$ of  $\cX_0$, which holds when $D$  is compatible with $\cX$.
The latter implies (1) and (2) hold.
\end{proof}

Using Proposition \ref{prop:twistFut}, 
we will deduce the following results.

\begin{prop}\label{prop:infFut/I-J}
Let $(X,\Delta)$ be a log Fano pair that is not uniformly K-stable, or equivalently, $\delta(X,\Delta)\leq 1$. 
Set $\mu^{\max}: = \sup \{ \mu \in (0,1] \, \vert \, (X,\Delta) \text{ is $\mu$-twisted K-semistable}\}$, which is a maximum and equals $\delta(X,\Delta)$ according to Theorem \ref{thm:twisted=delta}. 
Then
\begin{equation}\label{eq:infFut/I-J}
\inf_{(\cX,\Delta_{\tc};\cL)} 
\frac{\Fut(\cX,\Delta_{\tc};\cL)}
{\mnorm{\cX,\Delta_{\tc};\cL}}
= \mu^{\max}-1
,\end{equation}
where the infimum runs through all normal non-trivial test configurations of $(X,\Delta)$.
Furthermore, any normal non-trivial test configuration $(\cX,\Delta_{\tc};\cL)$  computes the above infimum if and only if $\Fut_{1-\mu^{\max}}(\cX,\Delta_{\tc};\cL)=0$. 
\end{prop}

\begin{proof}
For $\mu \in (0,1]$, Proposition \ref{prop:twistFut} (1) implies
$(X,\Delta)$ is $\mu$-twisted K-semistable iff
$$
\inf_{(\cX,\Delta_{\tc};\cL)} 
\frac{\Fut(\cX,\Delta_{\tc};\cL)}
{\mnorm{\cX,\Delta_{\tc};\cL}}
\geq \mu-1.$$
Since $(X,\Delta)$ is not uniformly K-stable, the previous infimum is $\leq 0$ and \eqref{eq:infFut/I-J} follows. The second statement is now a consequence of Proposition \ref{prop:twistFut} (1).
\end{proof}


 \begin{prop}\label{prop:twFut=A-muS}
If $(\cX,\Delta_{\tc};\cL)$ is a test configuration of $(X,\Delta)$ with $\cX_0$ integral, then 
 \[\Fut_{1-\mu}(\cX,\Delta_{\tc}; \cL) = A_{X,\Delta}(v_{\cX_0})- \mu S_{X,\Delta}(v_{\cX_0}).\]
 \end{prop}
 
 \begin{proof}
 The statement follows from combining Proposition \ref{prop:twistFut} (1) with  Proposition
\ref{prop:Fut-I-J}.
 \end{proof}

The next statement was originally proven by Dervan  under the slightly stronger assumption that $-m(K_X+\Delta)$ is very ample \cite[Section 2]{Der16}. 

\begin{thm}
Fix a positive integer $m$ such that
$-m(K_X+\Delta)$ is a Cartier divisor and $|-m(K_X+\Delta)|$ is base point free. 
If $(\cX,\Delta_{\tc};\cL)$ is a normal test configuration of $(X,\Delta)$
and $D \in \tfrac{1}{m} |-m(K_X+\Delta)|$ is general, then 
\[
\Fut_{1-\mu} (\cX,\Delta_{\tc}; \cL)  = \Fut(  \cX, \Delta_{\tc}+ (1-\mu)\cD;\cL ).\]
\end{thm}

\begin{proof}
The statement follows from combining Proposition \ref{prop:twistFut} (2) with Lemma \ref{l-compatgeneral}.
\end{proof}

\subsection{Special test configurations and twisted K-stability}
\label{sec:twistedLX}

In this subsection, we will establish the special test configuration theory for twisted K-stability,  as a natural generalization of \cite{LX14}. According to \cite[Corollary 1]{LX14}, to check the K-(semi)stability of a log Fano pair, it suffices to check the sign of the generalized Futaki invariant for all special test configurations. The goal of this section is to prove a twisted version of this statement.

\begin{thm}[=Theorem \ref{thm:main-twisted-LX}]\label{thm:twisted-LiXu}
Let $(X,\Delta)$ be a log Fano pair. Let $\mu\in (0,1]$ be a real number. Then $(X,\Delta)$ is $\mu$-twisted K-semistable (resp. K-stable) if and only if $\Fut_{1-\mu}(\cX,\Delta_{\tc};\cL)\geq 0$ (resp. $> 0$) for any non-trivial special test configuration $(\cX,\Delta_{\tc};\cL)$ of $(X,\Delta)$.
\end{thm}

\begin{proof}

We first recall the work of Li and Xu in \cite{LX14}. Suppose $X$ is a  $\bQ$-Fano variety of dimension $n$, and we want to check whether it is K-semistable. By definition one should check the non-negativity of the generalized Futaki invariant $\Fut(\mathcal{X},\mathcal{L})$ for all  test configurations $(\mathcal{X},\mathcal{L})$ of $(X,L:=-K_X)$. In \cite{LX14}, they create a process to simplify test configurations we need to check, that includes the following steps: semi-stable reduction and log canonical modification, anti-canonical polarization, and $\mathbb{Q}$-Fano extension. Along each step of their process, the generalized Futaki invariant decreases up to a multiple, so it is sufficient to check only special test configurations. 

Our strategy is to use that the twisted generalized Futaki invariant is the log generalized Futaki invariant for a general boundary divisor and use the argument of \cite{LX14} in the log pairs case. Let   $(\mathcal{X},\mathcal{L})\rightarrow C$ be a compactified  test configuration, where $C:=\mathbb{P}^1$ and $C^*:=\mathbb{P}^1\setminus0$, and  $\mathcal{X}_0$ denotes the central fiber of the test configuration.
We modify the given test configuration $(\mathcal{X},\mathcal{L})\rightarrow C$ of $(X,-K_X)$ step by step, such that the twisted generalized Futaki invariant decreases up to a multiple along the way. Since the argument in the $\bQ$-Fano case extends naturally to the log Fano pair case with little change, we will leave the boundary out for convenience.
\smallskip

\textit{Step 1. Finite base change. In this step, we consider $\cX^{(d)}$ as the normalization of $\mathcal{X}\times_C C'$, where $C'\cong\bP^1\rightarrow C$ is a  finite surjective morphism given by $t\mapsto t^d$ for some $d\in\bZ_{>0}$, such that $\mathcal{X}^{(d)}_0$ is reduced. Such a base change always exists by \cite{KKMSD}.
We aim to show that the twisted generalized Futaki invariant, suitably scaled, will strictly decrease after a base change if $\mathcal{X}_0$ is not reduced.}

Let $m$ be a sufficiently divisible positive integer such that $-mK_X$ is base point free. Then we can choose a general element $D\in \frac{1}{m}|-mK_X|$ such that $D$ is compatible with respect to both $\mathcal{X}$ and $\mathcal{X}^{(d)}$. We use $\mathcal{D}$ and $\mathcal{D}^{(d)}$ to denote the extension of $D$ on $\mathcal{X}$ and $\mathcal{X}^{(d)}$ respectively. Then we have 
$$\Fut(\mathcal{X},(1-\mu)\mathcal{D};\mathcal{L})=\frac{n\mu}{n+1} \frac{\mathcal{L}^{n+1}}{(-K_X)^n}+\frac{\mathcal{L}^n(K_{\mathcal{X}/C}+(1-\mu)\mathcal{D})}{(-K_X)^n} $$ 
and
$$ \Fut(\mathcal{X}^{(d)},(1-\mu)\mathcal{D}^{(d)};\mathcal{L}^{(d)})=\frac{n\mu}{n+1} \frac{({\mathcal{L}^{(d)}})^{n+1}}{(-K_X)^n}+\frac{({\mathcal{L}^{(d)})}^n(K_{\mathcal{X}^{(d)}/C'}+(1-\mu)\mathcal{D}^{(d)})}{(-K_X)^n} ,$$
where $\phi:\mathcal{X}^{(d)}\rightarrow \mathcal{X}$, $\mathcal{L}^{(d)}=\phi^*{\mathcal{L}}$, and 
$$K_{\mathcal{X}^{(d)}/C'}+(1-\mu)\mathcal{D}^{(d)} =\phi^*{(K_{\mathcal{X}/C}+(\mathcal{X}_0)_{\red}-\mathcal{X}_0+(1-\mu)\mathcal{D})} .$$
Combining the above relations and applying the projection formula, one has
\begin{align*}
\Fut(\mathcal{X}^{(d)},(1-\mu)\mathcal{D}^{(d)};\mathcal{L}^{(d)})& = {\rm deg}(\phi)\left( \Fut(\mathcal{X},(1-\mu)\mathcal{D};\mathcal{L}) -\frac{\cL^n\cdot (\cX_0-(\cX_0)_{\red})}{(-K_X)^n}\right)\\
& \leq {\rm deg}(\phi) \Fut(\mathcal{X},(1-\mu)\mathcal{D};\mathcal{L}).
\end{align*}
The equality holds if and only if $(\mathcal{X}_0)_{\red}=\mathcal{X}_0$ which means $\mathcal{X}_0$ is reduced.
So $$\Fut_{1-\mu}(\mathcal{X}^{(d)},\mathcal{L}^{(d)})\leq {\rm deg}(\phi) \Fut_{1-\mu}(\mathcal{X},\mathcal{L}) $$ and the equality holds if and only if  $\mathcal{X}_0$ is reduced. Note that the computation here is similar to \cite[Propositions 7.15 and 7.16]{BHJ17}.
\smallskip

\textit{Step 2. Log canonical modification. In this step, we aim to show that the twisted generalized Futaki invariant, suitably scaled, will strictly decrease if $(\mathcal{X},(1-\mu)\cD+\mathcal{X}_0)$ is not log canonical for a general divisor $D\in \frac{1}{m}|-mK_X|$, by running an MMP on its semi-stable reduction.}

Denote by $\cY$ the normalization of the graph of the birational map $X\times\bP^1\dashrightarrow \cX$. Choose $m$ as in Step 1. Denote by $\rho:\cY\to X\times\bP^1$ and $\tau:\cY\to \cX$ the two projection morphism. Consider the linear system $V$ on $\cX$ given as the movable part of  $\tau_*\rho^*(\pr_1^*|-mK_X|)$. Then it is clear that a general divisor in $V$ is exactly the extension $\cH$ of a general divisor $H\in |-mK_X|$ on $\cX$. Denote by $\cI$ the base ideal of $V$ on $\cX$. Since $|-mK_X|$ is base point free, we know that $\Supp(\cI)\subset\cX_0$. Let $\widehat{\cX}$ be the normalization of the blow up $\Bl_{\cI}\cX$. Let $\widehat{V}$ be the strict transform of $V$ on $\widehat{\cX}$. Then from the construction we know that $\widehat{V}$ is a base-point-free linear system on $\widehat{\cX}$. Denote by $\cE$ the reduced exceptional divisor of the morphism $\widehat{\cX}\to \cX$.

Next we will do a semi-stable reduction for $(\widehat{\cX},\cE)\to C$. By \cite[Theorem 7.17]{KM98},
there exist a finite surjective morphism $C'\cong\bP^1\to C$ given by $t\mapsto t^d$ for some $d\in\bZ_{>0}$, and a $\bG_m$-equivariant projective birational morphism $\pi:\cZ\to \cX'$ where $\cX'$ is the normalization of $\cX\times_C C'$, such that the following hold:
\begin{enumerate}[label=(\roman*)]
    \item The pair $(\cZ,\cZ_0+\Exc(\pi))$ has simple normal crossings, and $\cZ_0$ is reduced.
    \item $\pi$ factors through the morphism $\widehat{\cX}'\to \cX'$ where $\widehat{\cX}'$ is the normalization of $\widehat{\cX}\times_C C'$.
\end{enumerate}
Denote by $V_{\cZ}$ the pull-back of the linear system $\widehat{V}$ under the morphism $\cZ\to\widehat{\cX}$. Hence $V_{\cZ}$ is base point free since $\widehat{V}$ is base point free.
From the construction and Bertini's Theorem, we know that a general divisor in $V_{\cZ}$ is the strict transform $\pi_*^{-1}\cH'$ of an extension $\cH'$ on $\cX'$ of a general divisor $H\in |-mK_X|$. Thus again by Bertini's Theorem, we know that $(\cZ,\cZ_0+\Exc(\pi)+\pi_*^{-1}\cH')$ has simple normal crossings for a  general divisor $H\in|-mK_X|$.

Next we run an MMP to get a log canonical modification. For a divisor $D\in \frac{1}{m}|-mK_X|$, denote by $\cD'$ the extension of $D$ on $\cX'$. By \cite{LX14}, we can run a $\bG_m$-equivariant MMP $\mathcal{Z}\dashrightarrow \mathcal{X}^{\lc}$ on $\mathcal{Z}$ over $\mathcal{X}'$ to get a log canonical modification of $(\mathcal{X}',(1-\mu)\mathcal{D}'+\mathcal{X}_0')$ for a general divisor $D\in \frac{1}{m}|-mK_X|$, i.e. $(\mathcal{X}^{\lc},(1-\mu)\mathcal{D}^{\lc}+\mathcal{X}^{\lc}_0)$ is log canonical where $\cD^{\lc}$ is the extension of $D$. Note that the MMP sequence does not depend on the choice of a general $D$ since the divisor $m\cdot\pi_*^{-1}\cD'=\pi_*^{-1}\cH'$ belongs to the same linear system $V_{\cZ}$.
Since $\phi^{\lc} : \mathcal{X}^{\lc}\rightarrow \mathcal{X}'$ is a birational contraction morphism, we have:
 $$\Fut(\mathcal{X}',(1-\mu)\mathcal{D}';\mathcal{L}')= \Fut(\mathcal{X}^{\lc},(1-\mu){\mathcal{D}}^{\lc};(\phi^{\lc})^*\mathcal{L}'),$$
 where $\cL'$ is the pull-back of $\cL$ under the morphism $\cX'\to \cX$.
Now we construct a polarization $\mathcal{L}^{\lc}$ on $\mathcal{X}^{\lc}$ to make the twisted generalized Futaki invariant decrease.
Let
$$E:= \mu\cdot(\phi^{\lc})^*\mathcal{L}'+(1-\mu){\mathcal{D}}^{\lc}+K_{\mathcal{X}^{\lc}/C'}$$
and
$$ \mathcal{L}^{\lc}_t=(\phi^{\lc})^*\mathcal{L}'+tE.$$
Since
$$ \Fut(\mathcal{X}^{\lc},(1-\mu)\mathcal{D}^{\lc};\mathcal{L}^{\lc}_t)=\frac{n\mu}{n+1} \frac{(\mathcal{L}^{\lc}_t)^{n+1}}{(-K_X)^n}+\frac{(\mathcal{L}^{\lc}_t)^n\cdot(K_{\mathcal{X}^{\lc}/C'}+(1-\mu)\mathcal{D}^{\lc})}{(-K_X)^n}  , $$
we have
$$\left.\frac{d}{dt}\right|_{t=0} \Fut(\mathcal{X}^{\lc},(1-\mu)\mathcal{D}^{\lc};\mathcal{L}^{\lc}_t)=\frac{n((\phi^{\lc})^*\mathcal{L}')^{n-1}\cdot E^2}{(-K_X)^n}.$$
Note that $E$ is supported on the central fiber of $\mathcal{X}^{\lc}\rightarrow C'$, so the right hand side of the above equality is non-positive by Zariski's lemma (see \cite[2.4]{LX14}) and it is zero if and only if $E=\mu(\phi^{\lc})^*\mathcal{L}'+(1-\mu){\mathcal{D}}^{\lc}+K_{\mathcal{X}^{\lc}/C'}$ is linearly equivalent to a multiple of $\cX_0^{\lc}$. Recall that $\mathcal{X}^{\lc}$ is the log canonical modification of $(\mathcal{X}',(1-\mu)\mathcal{D}'+\mathcal{X}_0')$,  thus $E \sim_{\mathbb{Q},{\mathcal{X}'}}K_{\mathcal{X}^{\lc}/C'}+(1-\mu)\mathcal{D}^{\lc}$ is relative ample over $\mathcal{X}'$. 
 If $(\mathcal{X}',(1-\mu)\mathcal{D}'+\mathcal{X}_0')$ is not log canonical for a general $D$, then $\mathcal{X}'$ is not isomorphic to $\mathcal{X}^{\lc}$. Thus $E$ cannot be linearly equivalent to a multiple of $\cX_0^{\lc}$, which means that the above derivative is negative. 
Then we can choose a polarization $\mathcal{L}^{\lc}:=\mathcal{L}^{\lc}_t$ for some $0<t\ll1$ such that 
\begin{equation}\label{eq:lcmod}
\Fut(\mathcal{X}^{\lc},(1-\mu)\mathcal{D}^{\lc};\mathcal{L}^{\lc})<\Fut(\mathcal{X}^{\lc},(1-\mu)\mathcal{D}^{\lc};(\phi^{\lc})^*\mathcal{L}')=\Fut(\mathcal{X}',(1-\mu)\mathcal{D}';\mathcal{L}').
\end{equation}
Since $D$ is a general $\bQ$-divisor in $\frac{1}{m}|-mK_X|$, we may assume that it is compatible with both $\cX^{\lc}$ and $\cX'$. Thus \eqref{eq:lcmod} and Step 1 imply
\[
\Fut_{1-\mu}(\mathcal{X}^{\lc},\mathcal{L}^{\lc})<\Fut_{1-\mu}(\mathcal{X}',\mathcal{L}')\leq d\cdot  \Fut_{1-\mu}(\mathcal{X},\mathcal{L}).
\]
Here $d$ is the degree of the finite morphism $C'\to C$. If $(\cX',(1-\mu)\cD'+\cX_0')$ is log canonical for a general $D$ but $(\cX,(1-\mu)\cD+\cX_0)$ is not log canonical, then we know that $\cX\times_C C'$ is not normal. Hence \cite[Proposition 2.6]{BHJ17} implies that $\cX_0$ is not reduced. Therefore, Step 1 implies $$\Fut_{1-\mu}(\cX',\cL')<d\cdot \Fut_{1-\mu}(\mathcal{X},\mathcal{L}).$$
This finishes the proof of Step 2.
\smallskip

\textit{Step 3. Anti-canonical polarization.
By Step 2, we can assume that $(\mathcal{X},(1-\mu)\mathcal{D}+\mathcal{X}_0)$ is  log canonical for general $D\in \frac{1}{m}|-mK_X|$. 
We aim to show that the twisted generalized Futaki invariant will strictly decrease if the polarization is not anti-canonical.}

We run a $\bG_m$-equivariant MMP with scaling (see \cite[4.1]{LX14} or  \cite[section 6]{Fuj16}) for the log Fano pair $(X,(1-\mu)D)$ with general $D$. Then one obtains a finite sequence of birational contractions such that each is an isomorphism over $C\setminus\{0\}$:
$$\mathcal{X}=:\mathcal{X}^0\dashrightarrow \mathcal{X}^1 \dashrightarrow \mathcal{X}^2\dashrightarrow ...\dashrightarrow \mathcal{X}^k \rightarrow \mathcal{X}^{\ac}.$$
Here the last morphism is taking the anti-canonical model of $(\cX^k,(1-\mu)\cD^k)$ for a general $D$. Note that the above sequence of birational contractions does not depend on the choice of general $D$ since $\mathcal{D}$ are mutually $\mathbb{Q}$-linearly equivalent. We use $\mathcal{D}^i$ and $\cD^{\ac}$ to denote the extension of $D$ on $\mathcal{X}^i$ and $\mathcal{X}^{\ac}$ respectively.
Since $D$ is a general $\bQ$-divisor in $\frac{1}{m}|-K_X|$, we may assume that $D$ is compatible with respect to $\cX$. By Lemma \ref{l-compatcontract} we know that $D$ is also compatible with respect to $\cX^{\ac}$. Thus we have
$$\Fut_{1-\mu}(\mathcal{X}^{\ac},\mathcal{L}^{\ac})=\Fut(\mathcal{X}^{\ac},(1-\mu)\mathcal{D}^{\ac};\mathcal{L}^{\ac}) .$$
Here $\mathcal{L}^{\ac}$ is the pushforward of $\mathcal{L}$ to $\cX^{\ac}$, and $\mu \mathcal{L}^{\ac}=-(K_{\mathcal{X}^{\ac}/C}+(1-\mu)\mathcal{D}^{\ac})$. 
Then by \cite{LX14}, one directly has
$$ \Fut_{1-\mu}(\mathcal{X}^{\ac},\mathcal{L}^{\ac})=\Fut(\mathcal{X}^{\ac},(1-\mu)\mathcal{D}^{\ac};\mathcal{L}^{\ac})\leq \Fut(\mathcal{X},(1-\mu)\mathcal{D};\mathcal{L}) 
=\Fut_{1-\mu}(\cX,\cL),$$
and the equality holds if and only if $\mu\mathcal{L}\sim_{\bQ,C}-(K_{\mathcal{X}/C}+(1-\mu)\mathcal{D})$. This finishes the proof of Step 3.
\smallskip

\textit{Step 4. $\bQ$-Fano extension.
In this step, we can assume $(\mathcal{X},(1-\mu)\mathcal{D}+\mathcal{X}_0)$ is log canonical and  $\mu\mathcal{L}\sim_{\bQ,C}-(K_{\mathcal{X}/C}+(1-\mu)\mathcal{D})$ for general $D\in \frac{1}{m}|-mK_X|$.  
We aim to show that the twisted generalized Futaki invariant, suitably scaled, will strictly decrease if $(\mathcal{X},\mathcal{X}_0)$ is not plt.}

We first apply the semi-stable reduction process as in Step 2. Then there exist a finite surjective morphism $C'\cong\bP^1\to C$ given by $t\mapsto t^d$ for some $d\in\bZ_{>0}$, and a $\bG_m$-equivariant projective birational morphism $\pi:\cZ\to \cX'$ where $\cX'=\cX\times_C C'$ is normal, such that $(\cZ,\cZ_0+\Exc(\pi))$ has simple normal crossings, and $\cZ_0$ is reduced. Moreover, there exists a base-point-free linear system $V_{\cZ}$ on $\cZ$ such that a general divisor in $V_{\cZ}$ is exactly $m\cdot\pi_*^{-1}\cD'$ where $\cD'$ is the extension of a general $D\in \frac{1}{m}|-mK_X|$ on $\cX'$. Thus by Bertini's Theorem we have that $(\cZ,\cZ_0+\Exc(\pi)+\pi_*^{-1}\cD')$ has simple normal crossings for a general $D\in \frac{1}{m}|-mK_X|$.  
We then apply a pair version of \cite[Theorem 6]{LX14} for the log Fano pair $(X,(1-\mu)D)$. That is, we run a $\bG_m$-equivariant MMP over $C'$ to obtain a model $\mathcal{Z}\dashrightarrow \mathcal{X}^{\s}$ such that $(\mathcal{X}^{\s}, (1-\mu)\cD^{\s};\mathcal{L}^{\s}:=-(K_{\mathcal{X}^{\s}/C'}+(1-\mu)\mathcal{D}^{\s}))$ is a special test configuration of $(X,(1-\mu)D)$, where $a(\mathcal{X}_0^{\s};\mathcal{X}',(1-\mu)\mathcal{D}')=0$ and $\mathcal{X}_0^{\s}$ is a log canonical place of $(\mathcal{X}',\mathcal{X}_0')$. Here $\mathcal{D}^{\s}$ is the extension of $D$ on $\mathcal{X}^{\s}$. Note that the outcome special test configuration $(\mathcal{X}^{\s},\mathcal{L}^{\s}:=-(K_{\mathcal{X}^{\s}/C'}+(1-\mu)\mathcal{D}^{\s}))$ does not depend on the choice of a general $D\in \frac{1}{m}|-mK_X|$ since the divisor $m\cdot\pi_*^{-1}\cD'$ belongs to the same linear system $V_{\cZ}$, and $\cX_0^{\s}$ is a log canonical place of $(\cX',\cX_0'+(1-\mu)\cD')$ as well.
Then one has
$$\Fut(\mathcal{X}^{\s},(1-\mu)\mathcal{D}^{\s};\mathcal{L}^{\s}) \leq d\cdot  \Fut(\mathcal{X},(1-\mu)\mathcal{D};\mathcal{L}).$$
The equality holds if and only if $(\mathcal{X},(1-\mu)\mathcal{D}+\mathcal{X}_0)$ is plt. Since a general $D\in \frac{1}{m}|-mK_X|$ is compatible with respect to $\cX^{\s}$, the left hand side of above inequality is exactly $\Fut_{1-\mu}(\mathcal{X}^{\s},\mathcal{L}^{\s})$. Hence
$$\Fut_{1-\mu}(\mathcal{X}^{\s},\mathcal{L}^{\s}) \leq d\cdot\Fut(\mathcal{X},(1-\mu)\mathcal{D};\mathcal{L})\leq d\cdot\Fut_{1-\mu}(\mathcal{X},\mathcal{L}) .$$
If $(\mathcal{X},\mathcal{X}_0)$ is  not plt, then neither is $(\mathcal{X},(1-\mu)\mathcal{D}+\mathcal{X}_0)$ and so
$$\Fut_{1-\mu}(\mathcal{X}^{\s},\mathcal{L}^{\s}) < d\cdot\Fut_{1-\mu}(\mathcal{X},\mathcal{L}) .$$
This finishes the proof of Step 4 and, hence also the theorem.
\end{proof}

 We directly have following result:

\begin{corollary}\label{cor:twisted-special}
Let $(X,\Delta)$ be a $\mu$-twisted K-semistable log Fano pair with $0<\mu\leq 1$. Assume $(\mathcal{X},\Delta_{\tc};\mathcal{L})$ is a normal test configuration of $(X,\Delta)$ with $\Fut_{1-\mu}(\mathcal{X},\Delta_{\tc};\mathcal{L})=0$. Then this test configuration must be a special test configuration, and $(\mathcal{X},\mathcal{X}_0+\Delta_{\tc}+(1-\mu)\mathcal{D})$ is plt for general $D\in\frac{1}{m}|-m(K_X+\Delta)|$.
\end{corollary}

\begin{proof}
Assume to the contrary that this test configuration is not special. Then by the proof of Theorem \ref{thm:twisted-LiXu}, we can find a new special test configuration with negative twisted generalized Futaki invariant, which contradicts the fact that $(X,\Delta)$ is $\mu$-twisted K-semistable. The \textit{plt} part also follows from Step 4 of the proof of Theorem \ref{thm:twisted-LiXu}.
\end{proof}

\section{Optimal destabilization of log Fano pairs}\label{sec:optimal}

\subsection{Approximation of stability thresholds and twisted valuative criterion}\label{sec:approx}

In this section, we will construct sequences of divisors whose $A/S$ approximate the stability threshold. We prove Theorem \ref{thm:twisted=delta} using approximation results. We generalize the valuative criterion for K-semistability due to Fujita \cite{Fuj16} and Li \cite{Li17} to the twisted setting. 

Let us start from approximation by dreamy divisors.

\begin{theorem}[Dreamy approximation for $\delta<1$]\label{thm:dreamyapprox}
Let $(X,\Delta)$ be a log Fano pair with $\delta(X,\Delta) < 1$. Then there exists a sequence of dreamy divisors $\{E_i\}_{i\in\bN}$ over $X$ such that $$\lim_{i\to\infty}\frac{A_{X,\Delta}(E_i)}{S_{X,\Delta}(E_i)}=\delta(X,\Delta).$$
\end{theorem}

\begin{proof}
For simplicity, set $\delta:=\delta(X,\Delta)$ and $\delta_m:=\delta_m(X,\Delta)$. 
For each positive natural number $m$, there is an $m$-basis divisor $B_m$ such that $\lct(X,\Delta; B_m)=\delta_m$. Then we choose a prime divisor $E_m$ over $X$ computing $\lct(X,\Delta; B_m)$. Thus we have $$\delta_m=\frac{A_{X,\Delta}(E_m)}{{\rm ord}_{E_m}(B_m)}\geq \frac{A_{X,\Delta}(E_m)}{S_m(E_m)}\geq \delta_m .$$ 
In particular, $S_m(E_m)={\rm ord}_{E_m}(B_m)$.
As $\lim_{m\to\infty}\delta_m=\delta<1$, we have $\delta_m<1$ for $m\gg1$. Consider the log canonical pair $(X,\Delta+\delta_m B_m)$. Since $$A_{X,\Delta+\delta_mB_m}(E_m)=A_{X,\Delta}(E_m)-\delta_m {\rm ord}_{E_m}(B_m)=0,$$  by \cite[Corollary 1.38]{Kollar13} we can extract $E_m$ on a dlt modification $Y_m\rightarrow X$ of $(X,\Delta+\delta_m B_m)$ where $Y_m$ is $\bQ$-factorial. 

It is clear that $-(K_X+\Delta+\delta_m B_m)$ is nef for $m\gg 1$. By Lemma \ref{lem:Fanotype} below, we know that $Y_m$ is of Fano type. Hence by \cite[Corollary 1.3.2]{BCHM10} we know that $Y_m$ is a Mori dream space. This implies that $E_m$ is a dreamy divisor over $X$. Hence it suffices to show that $|\frac{A_{X,\Delta}(E_m)}{S_{X,\Delta}(E_m)}-\delta_m|$ goes to zero. Indeed, \cite[Corollary 2.10]{BJ17} implies that for any given $0<\epsilon<1$, there exists a positive integer $m(\epsilon)$ such that $S_m(v)\leq {(1+\epsilon)S(v)}$ for any $m>m(\epsilon)$ and any valuation $v\in \Val_X$. Thus for any $m>m(\epsilon)$ we have
$$\delta \leq\frac{A_{X,\Delta}(E_m)}{S_{X,\Delta}(E_m)}\leq (1+\epsilon)\frac{A_{X,\Delta}(E_m)}{S_m(E_m)} =(1+\epsilon)\delta_m.$$ 
So $|\frac{A_{X,\Delta}(E_m)}{S_{X,\Delta}(E_m)}-\delta_m|$ goes to zero. The proof is finished.
\end{proof}

\begin{lemma}\label{lem:Fanotype}
Suppose $(X,B)$ is a log Fano pair, and $(X,\Delta)$ is a log canonical pair such that $-(K_X+\Delta)$ is nef. If $Y\rightarrow (X,\Delta)$  is a dlt modification (see e.g. \cite[Theorem 1.34]{Kollar13}), then  $Y$ is of Fano type.
\end{lemma}

\begin{proof}
Denote the dlt modification by $f$, and suppose $F_i$ are exceptional divisors. Then we have 
$$K_Y+f^{-1}_*(B)+\sum{a_iF_i}= f^*(K_X+B)$$ 
and 
$$K_Y+f^{-1}_*(\Delta)+\sum {F_i}=f^*(K_X+\Delta).$$ 
Since $(X,B)$ is klt, for each $i$ we have $a_i<1$.
Then it is easy to see $$(Y,\epsilon(f^{-1}_*(B)+\sum{a_iF_i})+(1-\epsilon)(f^{-1}_*(\Delta)+\sum {F_i}))$$
is a klt weak log Fano pair for $0<\epsilon \ll 1$.
\end{proof}

With the above approximation result, we are able to prove Theorem \ref{thm:twisted=delta}. 

\begin{proof}[Proof of Theorem \ref{thm:twisted=delta}]
We first prove the ``$\geq$'' inequality. By continuity of $\Fut_{1-\mu}$
with respect to $\mu$, it suffices to show that $(X,\Delta)$ is $\mu$-twisted K-semistable
for any rational $0<\mu< \min\{1,\delta(X,\Delta)\}$. By \cite[Theorem 7.2]{BL18b},
there exists $D\in |-K_X-\Delta|_\mathbb{Q}$ depending on $\mu$ such that $(X,\Delta+(1-\mu)D)$ is K-semistable. This implies
the $\mu$-twisted K-semistability of $(X,\Delta)$ by Remark \ref{rem:logimpliestwisted}.

For the ``$\leq$'' inequality, it suffices to show that for any $0<\delta(X,\Delta)<\mu\leq 1$, $(X,\Delta)$ is not $\mu$-twisted K-semistable. Assume to the contrary that $(X,\Delta)$ is $\mu$-twisted K-semistable for some $\mu\in (\delta(X,\Delta), 1]$, then we have $\Fut_{1-\mu}(\mathcal{X}_E,\Delta_E;\mathcal{L}_E)\geq 0$ for any test configuration $(\mathcal{X}_E,\Delta_E;\mathcal{L}_E)$ induced by a dreamy divisor $E$ over $X$. Thus  Proposition \ref{prop:twFut=A-muS} implies that
$A_{X,\Delta}(E)-\mu S_{X,\Delta}(E)\geq 0 $
for any dreamy divisor $E$ over $X$. This contradicts the dreamy approximation for $\delta(X,\Delta)$, as proved in Theorem \ref{thm:dreamyapprox}.

Finally, we show that the supremum in \eqref{eq:sup=delta} is a maximum. Denote this supremum  by $\mu^{\rm sup}$. Let $\{\mu_i\}_{i\in \bZ_{>0}}$ be a sequence of numbers in $(0,1]$ such that $(X,\Delta)$ is $\mu_i$-twisted K-semistable for each $i$, and $\lim_{i\to\infty} \mu_i=\mu^{\rm sup}$. Then $\Fut_{1-\mu_i}(\cX,\Delta_{\tc};\cL)\geq 0$ for any normal test configuration $(\cX,\Delta_{\tc};\cL)$ and any $i$. By Proposition \ref{prop:twistFut} (1), we have \[
\Fut_{1-\mu^{\rm sup}}(\cX,\Delta_{\tc};\cL)=\lim_{i\to \infty} \Fut_{1-\mu_i}(\cX,\Delta_{\tc};\cL)\geq 0.
\]
Thus $(X,\Delta)$ is $\mu^{\rm sup}$-twisted K-semistable, which implies that $\mu^{\rm sup}$ is a maximum.
\end{proof}

Recall from Definition \ref{d:specialdiv} that a dreamy divisor $F$ over a log Fano pair $(X,\Delta)$ is special if $F$ induces a special test configuration of $(X,\Delta)$. 
The following result provides approximation to stability thresholds by special divisors. It is a stronger result than Theorem \ref{thm:dreamyapprox}.

\begin{theorem}[Special approximation for $\delta\leq 1$]\label{thm:specialapprox}
Let $(X,\Delta)$ be a log Fano pair with $\delta(X,\Delta)\leq 1$. Then there exists a sequence of special divisors $\{F_i\}_{i\in\bN}$ over $X$ such that 
$$\lim_{i\to\infty}\frac{A_{X,\Delta}(F_i)}{S_{X,\Delta}(F_i)}=\delta(X,\Delta) .$$
\end{theorem}

\begin{proof}
We first prove the theorem when $\delta(X,\Delta)<1$.
For any $i \gg 1 $ such that $\delta+\frac{1}{i} <1$, we know that $(X,\Delta)$ is not $(\delta+\frac{1}{i})$-twisted K-semistable by Theorem \ref{thm:twisted=delta}. Hence by Theorem \ref{thm:twisted-LiXu}, there must be a special test configuration $(\mathcal{X}_i,\Delta_{i};\cL_i)$  such that $\Fut_{1-\delta-\tfrac{1}{i}}(\mathcal{X}_i,\Delta_i;\cL_i)< 0 $. 
Let $F_i$ be the special divisor over $X$ induced by $(\cX_i,\Delta_i)$, i.e. $v_{\cX_{i,0}}=c_i\cdot \ord_{F_i}$ for some $c_i\in\bZ_{>0}$. Then 
\[
A_{X,\Delta}(F_i)-(\delta+\tfrac{1}{i})S(F_i)=c_i^{-1}\cdot\Fut_{1-\delta-\tfrac{1}{i}}(\mathcal{X}_i,\Delta_i;\cL_i)< 0.
\]
Thus $\delta(X,\Delta)\leq\frac{A_{X,\Delta}(F_i)}{S(F_i)}< \delta(X,\Delta)+\tfrac{1}{i}$. This implies 
\[
\lim_{i\to\infty}\frac{A_{X,\Delta}(F_i)}{S(F_i)}=\delta(X,\Delta) .
\]

Next we treat the case when $\delta(X,\Delta)=1$. This essentially follows from \cite[Theorem 2.9.4]{ZZ19} (see also \cite{Fuj16}). We provide a proof for the readers' convenience. 
Since $(X,\Delta)$ is not uniformly K-stable, by \cite[Corollary 3.4]{Fuj16} and \cite[Proposition 7.8]{BHJ17} there exists a sequence of special test configurations $\{(X_i,\Delta_i;\cL_i)\}_{i\in\bN}$ such that $\Fut(\cX_i,\Delta_i;\cL_i)< \tfrac{1}{i}\mnorm{\cX_i,\Delta_i;\cL_i}$. Let $F_i$ be the special divisor over $X$ induced by $(\cX_i,\Delta_i)$ as before. Then we have
\[
A_{X,\Delta}(F_i)-S(F_i)=c_i^{-1}\Fut(\cX_i,\Delta_i;\cL_i)
<i^{-1} c_i^{-1} \mnorm{\cX_i,\Delta_i;\cL_i}=i^{-1}S(F_i).
\]
Hence we have $1\leq \frac{A_{X,\Delta}(F_i)}{S(F_i)}< 1+\frac{1}{i}$. Thus
\[
\lim_{i\to\infty}\frac{A_{X,\Delta}(F_i)}{S(F_i)}=1=\delta(X,\Delta) .
\]
The proof is finished.
\end{proof}

 

The following theorem generalizes the valuative criterion for K-semistability \cite{Fuj16, Li17} to the twisted setting.


 
\begin{theorem}\label{thm:twisted-valuative-criterion}
 Let $(X,\Delta)$ be a log Fano pair. Let $0<\mu\leq 1$ be a real number. Then the following conditions are equivalent.
 \begin{enumerate}
     \item $(X,\Delta)$ is $\mu$-twisted K-semistable.
     \item $A_{X,\Delta}(E)\geq \mu S_{X,\Delta}(E)$ for any prime divisor $E$ over $X$.
     \item $A_{X,\Delta}(E)\geq \mu S_{X,\Delta}(E)$ for any dreamy prime divisor $E$ over $X$.
     \item $A_{X,\Delta}(E)\geq \mu S_{X,\Delta}(E)$ for any special prime divisor $E$ over $X$.
    \item $\mu\leq \min\{1, \delta(X,\Delta)\}$.
 \end{enumerate}
\end{theorem}

\begin{proof}
It is clear that (2)$\Rightarrow$(3)$\Rightarrow$(4). Theorem \ref{thm:twisted=delta} implies that (1)$\Leftrightarrow$(5). Hence it suffices to show (1)$\Rightarrow$(2) and (4)$\Rightarrow$(5).

For (1)$\Rightarrow$(2), 
suppose $(X,\Delta)$ is $\mu$-twisted K-semistable. Then $\mu\leq \delta(X,\Delta)$ by Theorem \ref{thm:twisted=delta}. So $A_{X,\Delta}(E)\geq \mu S_{X,\Delta}(E)$ for any prime  $E$ over $X$. 

For (4)$\Rightarrow$(5), we may assume that $(X,\Delta)$ is not uniformly K-stable since otherwise it is $\mu$-twisted K-semistable for any $\mu\in (0,1]$ by Theorem \ref{thm:twisted=delta}. Hence Theorem \ref{thm:specialapprox} implies that there exists a sequence of special divisors $F_i$ over $(X,\Delta)$ such that $\delta(X,\Delta)=\lim_{i\to\infty} \frac{A_{X,\Delta}(F_i)}{S(F_i)}$. By assumption we have $A_{X,\Delta}(F_i)\geq \mu S(F_i)$. Hence $\delta (X,\Delta)\geq \mu$.
\end{proof}

\subsection{Proof of Theorem \ref{thm:main1}}\label{sec:proof1.1}

In this section, we will prove Theorem \ref{thm:main1}. Let $(X,\Delta)$ be log Fano pair with $\delta(X,\Delta)\leq 1$.
First of all, by Theorem \ref{thm:twisted=delta} and Proposition \ref{prop:infFut/I-J} we know that the equality \eqref{eq:main1} holds, i.e. $\inf\frac{\Fut(\cdot)}{\mnorm{\cdot}}=\delta(X,\Delta)-1$.
By Proposition \ref{prop:infFut/I-J}, a normal non-trivial test configuration achieves $\inf \frac{\Fut(\cdot)}{\mnorm{\cdot}}$ if and only if it has vanishing $\delta(X,\Delta)$-twisted generalized Futaki invariant. Let $E$ be a prime divisor over $X$ computing $\delta(X,\Delta)$. 
We divide the proof into three parts, namely Theorems \ref{thm:delta-min-special}, \ref{thm:twisted-iff} and \ref{thm:delta-min-preserve}. In Theorem \ref{thm:delta-min-special}, we show that $E$ is dreamy and induces a special degeneration of $(X,\Delta)$ with vanishing twisted generalized Futaki invariant. This confirms Theorem \ref{thm:main1} (1) and the second part of \cite[Conjecture 1.5]{BX18}. Note that a special case of Theorem \ref{thm:delta-min-special} was shown by Blum and Xu in \cite[Theorem 4.1]{BX18} under the assumption of $\delta(X,\Delta)=1$. Later in Theorem \ref{thm:twisted-iff}, we will show that any special degeneration with vanishing twisted generalized Futaki invariant provides a divisorial valuation computing $\delta(X,\Delta)$. This confirms Theorem \ref{thm:main1} (2).  
In Theorem \ref{thm:delta-min-preserve}, we will show that such a special degeneration preserves stability thresholds. Combining Theorems \ref{thm:delta-min-special}, \ref{thm:twisted-iff}, and \ref{thm:delta-min-preserve}, we prove Theorem \ref{thm:main1}.

\begin{thm}\label{thm:delta-min-special}
Let $(X,\Delta)$ be a log Fano pair with $\delta(X,\Delta)\leq 1$. Assume that there exists a divisor $E$ over $X$ computing $\delta(X,\Delta)$. Then $E$ is dreamy and induces a non-trivial special test configuration $(\cX,\Delta_{\tc};\cL)$ of $(X,\Delta)$. In addition, we have $\Fut_{1-\delta(X,\Delta)}(\cX,\Delta_{\tc};\cL)=0$.
\end{thm}

\begin{proof}
We first show that any divisorial valuation $\ord_E$ computing $\delta(X,\Delta)\leq 1$ is dreamy by using the MMP. 
For each $m\in {\mathbb{Z}_{>0}}$ with $-m(K_X+\Delta)$ Cartier and $H^0(X,-m(K_X+\Delta))\neq 0$, we choose an $m$-basis divisor $B_m \in |-K_X-\Delta|_\mathbb{Q}$ such that $S_m(E)={\rm ord}_E(B_m)$.  Then we have $\frac{A_{X,\Delta}(E)}{{\rm ord}_E(B_m)} \rightarrow \delta(X,\Delta)$ as $m\to \infty$ since $S_m(E)\rightarrow S(E)$. Denote by $\delta:=\delta(X,\Delta)$ and $\delta_m:=\delta_m(X,\Delta)$ for simplicity. 
Choose a number $0<\epsilon<\min\{1,(\delta S(E))^{-1}\}$. Since $\lim_{m\to \infty}\delta_m=\delta$, we have $(1-\epsilon)\delta<\delta_m$ when $m\gg 1$. Thus  $(X,\Delta+(1-\epsilon)\delta B_m)$ is a log Fano pair for $m\gg {1}$. Then we have
$$A_{X,\Delta+(1-\epsilon)\delta B_m}(E)=A_{X,\Delta}(E)-{\rm ord}_E((1-\epsilon)\delta B_m)=\delta S(E)-(1-\epsilon)\delta S_m(E).$$
Since $\lim_{m\to \infty} S_m(E)=S(E)$, we know that
\[
\lim_{m\to\infty}A_{X,\Delta+(1-\epsilon)\delta B_m}(E)=\epsilon\delta S(E)<1.
\]
Hence we have $A_{X,\Delta+(1-\epsilon)\delta B_m}(E)<1$ for $m\gg 1$. 
So by \cite[Corollary 1.4.3]{BCHM10} for $m\gg 1$ there exists an extraction  $Y_m\to X$ of $E$. By \cite[Lemma 2.9]{Zhu18} we know that $Y_m$ is of Fano type, which implies that $E$ is dreamy.



Next we show that $E$ induces a special test configuration with vanishing $\delta(X,\Delta)$-twisted generalized Futaki invariant. 
For simplicity, denote by $\delta:=\delta(X,\Delta)$. Since $E$ is dreamy, we know that $E$ induces a test configuration $(\cX,\Delta_{\tc};\cL)$ of $(X,\Delta)$ with integral central fiber $\cX_0$ such that $v_{\cX_0}=\ord_E$. Hence Proposition \ref{prop:twFut=A-muS} implies that 
\[
\Fut_{1-\delta}(\cX,\Delta_{\tc};\cL)=A_{X,\Delta}(\ord_E)-\delta S(E)=0.
\]
Since $(X,\Delta)$ is $\delta$-twisted K-semistable by Theorem \ref{thm:twisted=delta}, the test configuration $(\cX,\Delta_{\tc};\cL)$ must be special by Corollary \ref{cor:twisted-special}. Thus the proof is finished.
\end{proof}

The following result shows that the vanishing of $\delta(X,\Delta)$-twisted generalized Futaki invariant in Theorem \ref{thm:delta-min-special} is not only necessary but also sufficient for the existence of a divisorial valuation computing $\delta(X,\Delta)$.

\begin{theorem}\label{thm:twisted-iff}
Let $(X,\Delta)$ be a log Fano pair with $\delta(X,\Delta)\leq 1$. If there is a normal non-trivial test configuration $(\mathcal{X},\Delta_{\tc};\mathcal{L})$ of $(X,\Delta)$ such that $\Fut_{1-\delta(X,\Delta)}(\mathcal{X},\Delta_{\tc};\mathcal{L})=0$, then this test configuration is special and $v_{\cX_0}$ is a divisorial valuation computing $\delta(X,\Delta)$. 
\end{theorem}

\begin{proof}
For simplicity, denote by $\delta:=\delta(X,\Delta)$.
Since $\Fut_{1-\delta}(\mathcal{X},\Delta_{\tc};\mathcal{L})=0$, the test configuration $(\mathcal{X},\Delta_{\tc};\mathcal{L})$ is special by Corollary \ref{cor:twisted-special}. By Proposition \ref{prop:twFut=A-muS}, we know that 
\[
A_{X,\Delta}(v_{\cX_0})=\delta S(v_{\cX_0})+\Fut_{1-\delta}(\mathcal{X},\Delta_{\tc};\mathcal{L})=\delta S(v_{\cX_0}).
\]
Hence, $v_{\cX_0}$ computes $\delta(X,\Delta)$. This finishes the proof.
\end{proof}

Next we turn to the second part of the proof of Theorem \ref{thm:main1}, that is, stability thresholds are preserved under the special degeneration induced by a divisorial valuation computing $\delta(X,\Delta)$. 

\begin{thm}\label{thm:delta-min-preserve}
Let $(X,\Delta)$ be a log Fano pair with $\delta(X,\Delta)\leq 1$. Assume that there exists a non-trivial special test configuration $(\cX,\Delta_{\tc};\cL)$ of $(X,\Delta)$ such that $\Fut_{1-\delta(X,\Delta)}(\cX,\Delta_{\tc};\cL)=0$. Then $\delta(\cX_0,\Delta_{\tc,0})=\delta(X,\Delta)$ are rational numbers.
\end{thm}

The key lemma to prove Theorem \ref{thm:delta-min-preserve} is the following generalization of \cite[Lemma 3.1]{LWX18} to twisted K-semistability.

\begin{lem}\label{lem:lwx3.1}
Let $(X,\Delta)$ be a $\mu$-twisted K-semistable log Fano pair for some $0<\mu\leq 1$. Let $(\cX,\Delta_{\tc};\cL)$ be a special test configuration of $(X,\Delta)$ satisfying $\Fut_{1-\mu}(\cX,\Delta_{\tc};\cL)=0$. Then $(\cX_0,\Delta_{\tc,0})$ is also $\mu$-twisted K-semistable. 
\end{lem}

\begin{proof}
Assume to the contrary that $(\mathcal{X}_0,\Delta_{\tc,0})$ is not $\mu$-twisted K-semistable. Since $(\mathcal{X}_0,\Delta_{\tc,0})$ admits a natural $\bG_m$-action from the test configuration $(\cX,\Delta_{\tc};\cL)$,  by Theorem \ref{thm:Ttwisted=twisted} we know that $(\mathcal{X}_0,\Delta_{\tc,0})$ is not $\bG_m$-equivariantly $\mu$-twisted K-semistable. Thus there is a $\bG_m$-equivariant special test configuration $(\mathcal{X}', \Delta_{\tc}';\cL')$ of $(X',\Delta'):=(\mathcal{X}_0, \Delta_{\tc,0})$ such that $\Fut_{1-\mu}(\mathcal{X}', \Delta_{\tc}';\cL')<0$. 
Choose a sufficiently divisible positive integer $m$ such that $-m(K_{X}+\Delta)$ is Cartier and base point free. 
Then one can choose a general $D\in\frac{1}{m}|-m(K_X+\Delta)|$ such that $D$ is compatible with $\cX$, i.e.
$$\Fut_{1-\mu}(\mathcal{X},\Delta_{\tc};\cL)=\Fut(\mathcal{X},\Delta_{\tc}+(1-\mu)\mathcal{D};\cL)=0,$$
where $\mathcal{D}$ is the extension of $D$ on $\mathcal{X}$. By Corollary \ref{cor:twisted-special}, we have that $(\mathcal{X}_0,\Delta_{\tc,0}+(1-\mu)\mathcal{D}_0)$ is a special degeneration of $(X,\Delta+(1-\mu)D)$ for general $D\in\frac{1}{m}|-m(K_X+\Delta)|$.
Then we have
$$\Fut(\mathcal{X}',\Delta_{\tc}'+(1-\mu)\mathcal{D}';\cL')\leq\Fut_{1-\mu}(\mathcal{X}',\Delta_{\tc}';\cL')<0,$$
where $\mathcal{D}'$ is the extension of $D':=\mathcal{D}_0$ on $\mathcal{X}'$.

By \cite[Proof of Lemma 3.1]{LWX18}, for $k\gg 1$ we can construct a new special test configuration $(\mathcal{X}'',\Delta_{\tc}'';\cL'')$ of $(X,\Delta)$ independent of the choice of a general $D\in \frac{1}{m}|-(mK_X+\Delta)|$, such that 
$$\Fut(\mathcal{X}'',\Delta_{\tc}''+(1-\mu)\mathcal{D}'';\cL'')= k\cdot \Fut(\mathcal{X},\Delta_{\tc}+(1-\mu)\mathcal{D};\cL)+\Fut(\mathcal{X}',\Delta_{\tc}'+(1-\mu)\mathcal{D}';\cL')<0 ,$$
where $\mathcal{D}''$ is the extension of $D$ on $\mathcal{X}''$. Since $D$ is a general divisor in $\frac{1}{m}|-m(K_X+\Delta)|$, we may assume that it is also compatible with $\cX''$. Thus
$$ \Fut_{1-\mu}(\mathcal{X}'',\Delta_{\tc}'';\cL'')=\Fut(\mathcal{X}'',\Delta_{\tc}''+(1-\mu)\mathcal{D}'';\cL'')<0,$$
which is a contradiction to $\mu$-twisted K-semistability of $(X,\Delta)$. The proof is finished.
\end{proof}

\begin{proof}[Proof of Theorem \ref{thm:delta-min-preserve}]
For simplicity, denote by $\delta:=\delta(X,\Delta)$
By Theorem \ref{thm:twisted=delta}, we know $(X,\Delta)$ is $\delta$-twisted K-semistable. Thus by Lemma \ref{lem:lwx3.1} we know that $(\cX_0,\Delta_{\tc,0})$ is also $\delta$-twisted K-semistable, which implies $\delta(\cX_0,\Delta_{\tc,0})\geq \delta$ by Theorem \ref{thm:twisted=delta}. Since the stability thresholds are lower semicontinuous in families by \cite[Theorem B]{BL18b}, we know $\delta(\cX_0,\Delta_{\tc,0})\leq \delta$. Therefore, we have $\delta(\cX_0,\Delta_{\tc,0})= \delta$. By Proposition \ref{prop:Fut-I-J}, we know 
\[
\Fut(\cX,\Delta_{\tc};\cL)=A_{X,\Delta}(v_{\cX_0})-S(v_{\cX_0}).
\]
It is clear that both $\Fut(\cX,\Delta_{\tc};\cL)$ and $A_{X,\Delta}(v_{\cX_0})$ are rational. Thus $S(v_{\cX_0})$ is rational which implies $\delta=A_{X,\Delta}(v_{\cX_0})/S(v_{\cX_0})$ is also rational. The proof is finished. 
\end{proof}

\begin{proof}[Proof of Theorem \ref{thm:main1}]
The equality \eqref{eq:main1} follows from 
Theorem \ref{thm:twisted=delta} and Proposition \ref{prop:infFut/I-J}.
By Proposition \ref{prop:infFut/I-J}, a normal non-trivial test configuration $(\cX,\Delta_{\tc};\cL)$ of $(X,\Delta)$ achieves the infimum in \eqref{eq:main1} if and only if $\Fut_{1-\delta(X,\Delta)}(\cX,\Delta_{\tc};\cL)=0$.
Hence part (1) follows directly from Theorem \ref{thm:delta-min-special}. Part (2) follows directly from Theorem \ref{thm:twisted-iff}. The last statement follows directly from Theorem \ref{thm:delta-min-preserve}.
\end{proof}

Note that the proof of Lemma \ref{lem:lwx3.1} uses a result, namely Theorem \ref{thm:Ttwisted=twisted}, on the equivalence between twisted K-semistability and its torus equivariant analogue. This is treated in the below and Section \ref{sec:twistednv}. 

\begin{definition}
Let $(X,\Delta)$ be a log Fano pair with a $\bT$-action, where $\bT=\bG_m^r$ with $r\in\bZ_{\geq 0}$ is an algebraic torus.
\begin{enumerate}
    \item A test configuration $(\cX,\Delta_{\tc};\cL)$ of $(X,\Delta)$ is said to be  a \emph{$\bT$-equivariant test configuration} if there exists a $\bT$-action on $(\cX,\Delta_{\tc};\cL)$ commuting with its $\bG_m$-action, such that 
    the canonical isomorphism $(\cX\setminus\cX_0,\Delta_{\tc}|_{\cX\setminus\cX_0};\cL|_{\cX\setminus\cX_0})\xrightarrow{\cong} (X,\Delta;L)\times(\bA^1\setminus\{0\})$ is $(\bG_m\times \bT)$-equivariant with the trivial $\bT$-action on $\bA^1\setminus\{0\}$.

    If in addition that $(\cX,\Delta_{\tc};\cL)$ is a special test configuration, then we say that it is a \emph{$\bT$-equivariant special test configuration}.

\item For a real number $0<\mu\leq 1$, we say that  $(X,\Delta)$ is \emph{$\bT$-equivariantly $\mu$-twisted K-semistable} if for any $\bT$-equivariant test configuration  $(\cX,\Delta_{\tc};\cL)$, we have
$$\Fut_{1-\mu}(\cX,\Delta_{\tc};\cL)\geq 0.$$
\end{enumerate}
\end{definition}

The following result generalizes \cite[Theorem E]{LX16} to twisted K-stability. We will present its proof in Section \ref{sec:twistednv} using twisted normalized volume on the affine cone.

\begin{theorem}\label{thm:Ttwisted=twisted}
Let $(X,\Delta)$ be a log Fano pair with a $\bT$-action. Let $0<\mu\leq 1$ be a real number.
The following are equivalent.
\begin{enumerate}
    \item $(X,\Delta)$ is $\mu$-twisted K-semistable.
    \item $(X,\Delta)$ is $\bT$-equivariantly $\mu$-twisted K-semistable.
    \item For any $\bT$-equivariant special test configuration $(\cX,\Delta_{\tc};\cL)$ of $(X,\Delta)$, we have
$$\Fut_{1-\mu}(\cX,\Delta_{\tc};\cL)\geq 0.$$
\end{enumerate}
\end{theorem}

\subsection{Twisted normalized volume}\label{sec:twistednv}

In this section, we will study normalized volumes of klt singularities twisted by ideal boundaries. We establish a result on approximation to twisted normalized volume using $\bT$-equivariant Koll\'ar componenents (see Lemma \ref{lem:twistednv-approx}). We show that twisted K-semistability can be characterized by twisted normalized volume on the affine cone (see Theorem \ref{t:deltahvol}). These results enable us to prove Theorem \ref{thm:Ttwisted=twisted}, i.e. twisted K-semistability of a log Fano pair with a $\bT$-action can be tested by $\bT$-equivariant special test configurations.

\begin{defn}
 Let $x\in (X,\Delta)$ be an $n$-dimensional klt singularity. Let $\fa\subset\fm:=\fm_x$ be an ideal sheaf. Let $c\geq 0$ be a real number. We define the \emph{log discrepancy} of a valuation  $v\in\Val_X$ with respect to the pair $(X,\Delta+c\cdot\fa)$ as
 \[
 A_{(X,\Delta+c\cdot\fa)}(v):=A_{(X,\Delta)}(v)-c\cdot v(\fa).
 \]
 Assume in addition that $(X,\Delta+c\cdot\fa)$ is klt, i.e. $A_{(X,\Delta+c\cdot\fa)}(v)>0$ for any $v\in\Val_X$. 
 Then we define the \emph{normalized volume function} 
 $\hvol_{(X,\Delta+c\cdot\fa),x}:\Val_{X,x}\to (0,+\infty]$ as
 \[
  \hvol_{(X,\Delta+c\cdot\fa),x}(v):=\begin{cases}
  A_{(X,\Delta+c\cdot\fa)}(v)^n\cdot\vol_{X,x}(v), &\textrm{ if }A_{(X,\Delta)}(v)<+\infty;\\
  +\infty, &\textrm{ if }A_{(X,\Delta)}(v)=+\infty.
  \end{cases}
 \]
 We define the \emph{local volume} of a klt singularity $x\in (X,\Delta+c\cdot\fa)$ as
 \[
 \hvol(x,X,\Delta+c\cdot\fa):=\inf_{v\in\Val_{X,x}}\hvol_{(X,\Delta+c\cdot\fa),x}(v).
 \]
\end{defn}

The above definition can be viewed as an analogue of the usual normalized volume function defined by Chi Li \cite{Li18} \emph{twisted} by an ideal boundary. Hence we also call it \emph{twisted normalized volume}. 

We first give a result on approximation to twisted normalized volume. This is a generalization of \cite[Proposition 4.4]{LX16}
and \cite[Theorem 27]{Liu18}.

\begin{lemma}\label{lem:twistednv-approx}
Let $x\in (X,\Delta)$ be a klt singularity with a $\bT$-action where $\bT$ is an algebraic torus. Let $\fa\subset \fm$ be a $\bT$-invariant ideal sheaf. Let $c\geq 0$ be a real number such that $(X,\Delta+c\cdot\fa)$ is klt. Then we have
\begin{equation}\label{eq:lx16-approx}
\hvol(x, X,\Delta+c\cdot\fa)=\inf_{\fb}\lct(X,\Delta+c\cdot\fa;\fb)^n\cdot \mult(\fb)=\inf_S\hvol_{(X,\Delta+c\cdot\fa),x}(\ord_S),
\end{equation}
where $\fb$ runs over all $\bT$-invariant $\fm$-primary ideals, and $S$ runs over all $\bT$-equivariant Koll\'ar components over $x\in (X,\Delta)$.
\end{lemma}

\begin{proof}
Firstly, we have a direct generalization of \cite[Theorem 27]{Liu18} that
\begin{equation}\label{eq:twisted-liu}
\hvol(x, X,\Delta+c\cdot\fa)=\inf_{\fc}\lct(X,\Delta+c\cdot\fa;\fc)^n\cdot \mult(\fc)=\inf_{\fc_\bullet}\lct(X,\Delta+c\cdot\fa;\fc_\bullet)^n\cdot \mult(\fc_\bullet),
\end{equation}
where $\fc$ runs over all $\fm$-primary ideals, and $\fc_\bullet$ runs over all graded sequence of $\fm$-primary ideals. Indeed, the proof of \eqref{eq:twisted-liu} is  the same as \cite[Proof of Theorem 27]{Liu18} after replacing $X$ by $(X,\Delta+c\cdot\fa)$. 
Let $\fc$ be an arbitrary $\fm$-primary ideal. Denote by $\fc_i:=\fc^i$. Assume $\bT\cong(\bG_m)^r$. Fixing a lexicographic order on $\bZ^r$, we can degenerate $\fc_i$ to its initial ideal $\fb_i:=\bin(\fc_i)$. By \cite[Lemmas 4.1, 4.2, and 4.3]{LX16}, we know that $\fb_\bullet$ is a graded sequence of $\bT$-invariant $\fm$-primary ideals, and 
\begin{align*}
\lct(X,\Delta+c\cdot\fa;\fc_\bullet)^n\cdot \mult(\fc_\bullet)&\geq \lct(X,\Delta+c\cdot\fa;\fb_\bullet)^n\cdot \mult(\fb_\bullet)\\
&=\lim_{i\to\infty}\lct(X,\Delta+c\cdot\fa;\fb_i)^n\cdot \mult(\fb_i).
\end{align*}
This combined with \eqref{eq:twisted-liu} implies the first equality of \eqref{eq:lx16-approx}. 

For the second equality of \eqref{eq:lx16-approx}, it suffices to show that for any $\bT$-invariant $\fm$-primary ideal $\fb$, there exists a $\bT$-equivariant Koll\'ar component $S$ over $(X,\Delta)$ such that 
\[
\lct(X,\Delta+c\cdot\fa;\fb)^n\cdot \mult(\fb)\geq \hvol_{(X,\Delta+c\cdot\fa),x}(\ord_S).
\]
This follows from the same argument as \cite[Proof of Proposition 4.4]{LX16} after replacing $(X,D,\fa)$ by $(X,\Delta+c\cdot\fa,\fb)$. Indeed, we can take a $\bT$-equivariant log resolution of $(X,\Delta+c\cdot\fa,\fb)$, and then run $\bT$-equivariant MMP to get a dlt modification $Y\to X$ and then a Koll\'ar component $S$ over $x\in (X,\Delta)$, such that 
\[
\lct(X,\Delta+c\cdot\fa;\fb)^n\cdot \mult(\fb)\geq  \vol_{(X,\Delta+c\cdot\fa),x}(Y)\geq \hvol_{(X,\Delta+c\cdot\fa),x}(\ord_S).
\]
Hence the proof is finished.
\end{proof}

In the rest of this subsection, we denote by $(V,D)$ an $(n-1)$-dimensional log Fano pair. Choose $r\in\bN$ sufficiently divisible such that $-r(K_V+D)$
is Cartier, and the section ring $R(V,-r(K_V+D))$ is generated in degree $1$. Let $(X,\Delta)$ be the affine cone $C(V,D;-r(K_V+D))$ with cone vertex $x\in X$. Denote by $\ord_{V}$ the divisorial valuation induced by blowing up the cone vertex $x$. 

The following result characterizes twisted K-semistability of a log Fano pair using twisted normalized volume on its affine cone. It generalizes the normalized volume characterization of K-semistability in \cite[Theorem 4.5]{LX16} (see also \cites{Li17, LL19}).

\begin{thm}\label{t:deltahvol}
The log Fano pair $(V,D)$ is $\mu$-twisted K-semistable for $0<\mu\leq 1$ if and only if $\ord_V$ minimizes $\widehat{\vol}_{(X,\Delta+\frac{1-\mu}{r}\mathfrak{m}),x}$.
\end{thm}

\begin{proof}
By Theorem \ref{thm:twisted=delta}, we know that $(V,D)$ is $\mu$-twisted K-semistable if and only if $0<\mu\leq \min\{1,\delta(V,D)\}$. Hence by taking $c:=\frac{1-\mu}{r}$ it suffices to show that $\max\{0,1-\delta(V,D)\}\leq cr<1$ if and only if  $\ord_V$ minimizes $\hvol_{(X,\Delta+c\cdot\fm),x}$. 

For the ``only if'' part, we may assume in addition that $\max\{0,1-\delta(V,D)\}< cr<1$ and $c$ is rational since $\hvol_{(X,\Delta+c\cdot\fm),x}$
is continuous in $c$. Denote by $\gamma:=1-cr\in (0,\min\{1,\delta(V,D)\})$.
Then by \cite[Theorem 7.2]{BL18b}, there exists an
effective $\bQ$-divisor $B\sim_{\bQ}-K_V-D$ on $V$
such that $(V,D+(1-\gamma)B)$ is log K-semistable.
Hence \cite[Theorem 4.5]{LX16} implies that 
$\ord_V$ minimizes $\hvol_{(X,\Delta+(1-\gamma)C(B)),x}$. It is clear that
$\ord_V(C(B))=r^{-1}$, hence $v(C(B))\geq r^{-1}v(\fm)$ for any $v\in\Val_{X,x}$.
Thus for any $v\in\Val_{X,x}$, we have
\begin{align*}
\hvol_{(X,\Delta+c\cdot\fm),x}(v)
&\geq \hvol_{(X,\Delta+(1-\gamma)C(B)),x}(v)\\
&\geq \hvol_{(X,\Delta+(1-\gamma)C(B)),x}(\ord_V)\\
&=\hvol_{(X,\Delta+c\cdot\fm),x}(\ord_V).
\end{align*}
This finishes the proof of the ``only if'' part.

For the ``if'' part, we follow the strategy of
\cite{Li17}. Let $\ord_E$ be any divisorial valuation of $K(V)$. Then $\ord_E$ extends trivially to a valuation $\ord_{C(E)}$ of $K(X)\cong K(V)(s)$. For any $t\in[0,\infty)$, denote by $v_t$ the quasi-monomial valuation of weights $(1,t)$ along $(V,C(E))$. Then the computations in \cite{Li17}
shows that 
\[
\left.\frac{d}{dt}\right|_{t=0}\hvol_{(X,\Delta),x}(v_t)
=n\beta(\ord_E).
\]
Here $\beta(\ord_E):=(-K_V-D)^{n-1}(A_{(V,D)}(E)-S_{(V,D)}(E))$ is Fujita's $\beta$-invariant (see \cite[Definition 1.3]{Fuj16}).
It is clear that $v_t(\fm)=1$ for any $t$. Hence we have 
\[
A_{(X,\Delta)}(v_t)=r^{-1}+tA_{(V,D)}(\ord_E), \quad
A_{(X,\Delta+c\cdot\fm)}(v_t)=r^{-1}-c+tA_{(V,D)}(\ord_E).
\]
Hence we have
\begin{align*}
    \left.\frac{d}{dt}\right|_{t=0}\hvol_{(X,\Delta+c\cdot\fm),x}(v_t)& =\left.\frac{d}{dt}\right|_{t=0}
    \left(\left(\frac{A_{(X,\Delta+c\cdot\fm)}(v_t)}{A_{(X,\Delta)}(v_t)}\right)^n\cdot\hvol_{(X,\Delta),x}(v_t)\right)\\
    &=n(1-cr)^n\left(\beta(\ord_E)+\frac{c}{r^{-1}-c}A_{(V,D)}(\ord_E)\vol_V(-K_V-D)\right)\\
    &=n(1-cr)^{n-1}\vol_V(-K_V-D)\left(A_{(V,D)}(\ord_E)-(1-cr)S(\ord_E)\right). \numberthis \label{eq:Li-derivative}
\end{align*}
Since $\ord_V$ minimizes $\hvol_{(X,\Delta+c\cdot\fm),x}$, the above derivative
is non-negative. This implies that $A_{(V,D)}(E)\geq (1-cr)S(\ord_E)$ for any prime divisor $E$ over $V$, i.e. $\delta(V,D)\geq 1-cr$. The proof is finished.
\end{proof}

Based on above results, we are ready to prove Theorem \ref{thm:Ttwisted=twisted}. 

\begin{proof}[Proof of Theorem \ref{thm:Ttwisted=twisted}]
The implications (1)$\Rightarrow$(2) and (2)$\Rightarrow$(3) are obvious. Thus, it suffices to show that (3) implies (1). For consistency of notation, we denote the log Fano pair with a $\bT$-action by $(V,D)$. Then we know $\Fut_{1-\mu}(\cV,\cD;\cL)\geq 0$ for any $\bT$-equivariant special test configuration $(\cV,\cD;\cL)$ of $(V,D)$. By Theorem \ref{t:deltahvol} it suffices to show that $\ord_V$ minimizes $\hvol_{(X,\Delta+c\cdot\fm),x}$ with $c:=\frac{(1-\mu)}{r}$. By Lemma \ref{lem:twistednv-approx}, this is equivalent to showing that 
\begin{equation}\label{eq:T-Kc}
\hvol_{(X,\Delta+c\cdot\fm),x}(\ord_S)\geq \hvol_{(X,\Delta+c\cdot\fm),x}(\ord_V)    
\end{equation}
for any $(\bT\times\bG_m)$-equivariant Koll\'ar component $S$ over $x\in (X,\Delta)$. 

We follow the argument of \cite[Proof of Theorem 4.5]{LX16}. Denote by $L:=-r(K_V+D)$. Let $R:=\oplus_{j=0}^\infty R_m=\oplus_{j=0}^\infty H^0(X,mL)$
be the section ring of $(V, L)$, then $X=\Spec R$.
For a  $(\bT\times\bG_m)$-equivariant Koll\'ar component $S$ over $x\in (X,\Delta)$, consider the filtration of $R$ as 
\[
\cF^p R_m:= \{f\in R_m\mid \ord_S(f)\geq p\}.
\]
Let $Y\to X$ be the blow up of the cone vertex $x$ whose exceptional divisor is still denoted by $V$. Denote by $\cI_V$ the ideal sheaf of $V$ in $Y$. We define
\[
c_1:=\min\{\ord_S(f)\mid f\in \cI_V(U),~ U\subset Y\textrm{ is open, and } U\cap c_{Y}(\ord_S)\neq\emptyset\}.
\]
For $t\in\bR$, $\lambda\in \bR_{>0}$, and $s\in[0,1]$, define
\[
\vol(\cF R^{(t)}):=\limsup_{m\to\infty}\frac{\dim(\cF^{mt}R_m)}{m^n/n!}, \qquad
\Phi(\lambda, s):=\int_{c_1}^{+\infty}\frac{-d\vol(\cF R^{(t)})}{((1-s)+\lambda st)^n}.
\]
By \cite{Li17} we know that $\Phi(\lambda,s)$ satisfies the following properties.
\begin{enumerate}
    \item For any $\lambda\in\bR_{>0}$ we have $\Phi(\lambda,1)=\lambda^{-n}\vol(\ord_S)$ and $\Phi(\lambda,0)=\vol(\ord_V)$;
    \item For any $\lambda\in \bR_{>0}$ the function $s\mapsto\Phi(\lambda,s)$ is continuous and convex for $s\in[0,1]$.
    \item The directional derivative of $\Phi(\lambda,s)$ at $s=0$ is equal to 
    \[
    \Phi_s(\lambda,0)=n \lambda (L^{n-1})\left(\lambda^{-1}-c_1-\frac{1}{(L^{n-1})}\int_{c_1}^{+\infty}\vol(\cF R^{(t)})dt \right).
    \]

\end{enumerate}
Let $\lambda_*:=\frac{\mu}{r(A_{(X,\Delta+c\cdot\fm)}(\ord_S))}$. Since $A_{(X,\Delta+c\cdot\fm)}(\ord_V)=\frac{1}{r}-c=\frac{\mu}{r}$, we know that 
\[
    \frac{\hvol_{(X,\Delta+c\cdot\fm),x}(\ord_S)}{\hvol_{(X,\Delta+c\cdot\fm),x}(\ord_V)}=\frac{\Phi(\lambda_*,1)}{\Phi(\lambda_*,0)}.
\]
Therefore, by convexity of $\Phi(\lambda,\cdot)$, to show \eqref{eq:T-Kc} it suffices to show $\Phi_s(\lambda_*,0)\geq 0$.

Consider the restriction valuation $\ord_S|_{K(V)}$ under the field embedding $K(V)\hookrightarrow K(X)$. From \cite[Proof of Theorem 4.5]{LX16} and references therein, there exists a prime divisor $E$ over $V$ and $b\in\bZ_{\geq 0}$ such that $\ord_S|_{K(V)}=b\cdot\ord_E$. Moreover, if $\tilde{V}\to V$ is a birational model containing $E$ as a divisor, then $\ord_S$ is a quasi-monomial combination of $\tilde{V}$ and $E\times_{\tilde{V}}\tilde{Y}$ in $\tilde{Y}:=Y\times_V \tilde{V}$ of weights $(c_1, b)$. In particular, we have $\ord_S(\fm)=c_1$, $A_{(X,\Delta)}(\ord_S)=\frac{c_1}{r}+b A_{(V,D)}(\ord_E)$, and 
\[
\cF^p R_m=\{f\in R_m\mid b\cdot \ord_E(f)\geq p -c_1m\}.
\]
Thus straightforward computation shows that $\vol(\cF R^{(t)})= \vol(L-\frac{t-c_1}{b}E)$. Hence we have 
\begin{align*}
   \frac{ \Phi_s(\lambda_*,0) }{n\lambda_*(L^{n-1})}& = \frac{r}{\mu}\left(A_{(X,\Delta)}(\ord_S)-\frac{1-\mu}{r}\ord_S(\fm)\right)-c_1-\frac{1}{(L^{n-1})}\int_{c_1}^{+\infty}\vol(L-\tfrac{t-c_1}{b}E)dt\\
   & = \frac{rb}{\mu} A_{(V,D)}(\ord_E)+\frac{c_1}{\mu}-\frac{1-\mu}{\mu}c_1-c_1-\frac{b}{(L^{n-1})}\int_0^{+\infty}\vol(L-tE)dt\\
   & = \frac{rb}{\mu}(A_{(V,D)}(\ord_E)-\mu S(\ord_E)).
\end{align*}
Thus it suffices to show $A_{(V,D)}(\ord_E)-\mu S(\ord_E)\geq 0$. Denote by $c_2:= c_1+ rbA_{(V,D)}(\ord_E)$. We define another filtration $\cF'$ of $R$ as  
\[
\cF'^p R_m:= \cF^{p+c_2 m} R_m= \{f\in R_m \mid b\cdot\ord_E(f)\geq p+ rm A_{(V,D)}(b\cdot\ord_E)\}.
\]
Since $S$ is a $\bG_m$-equivariant invariant Koll\'ar component, we know $\cF$ (and hence $\cF'$) is finitely generated.
By \cite[Proposition 2.15 and Lemma 5.17]{BHJ17}, we know that the filtration $\cF'$ induces a test configuration $(\cV,\cD;\cL)$ of $(V,D)$, such that $\cL\sim_{\bQ} -r(K_{\cV}+\cD)$,  $\cV_0$ is integral,  and $v_{\cV_0}=b\cdot \ord_E$. Since $\ord_S$ is $(\bT\times\bG_m)$-invariant, we know that $\ord_E$ is $\bT$-invariant which implies that $(\cV,\cD;\cL)$ is $\bT$-equivariant. 
Moreover, we know 
\[
\cV_0\cong \Proj\bigoplus_{m=0}^\infty\bigoplus_{p=-\infty}^{\infty}\cF'^p R_m/\cF'^{p+1} R_{m}=\Proj \bigoplus_{m=0}^\infty\bigoplus_{p=-\infty}^{\infty}\cF^{p+c_2 m} R_m/\cF^{p+c_2 m-1} R_m.
\]
Let $(X_0,\Delta_0)$ be the klt special degeneration of $(X,\Delta)$ induced by $\ord_S$. Then we know that 
\[
X_0=\Spec\bigoplus_{m=0}^\infty\bigoplus_{p=-\infty}^{\infty}\cF^p R_m/\cF^{p+1} R_{m}.
\]
Thus, $(\cV_0,\cD_0)$ is a quotient of the klt pair $(X_0,\Delta_0)$ by a $\bG_m$-action corresponding to the grading of $p+c_2m$. So $(\cV_0,\cD_0)$ is a log Fano pair which implies that $(\cV,\cD;\cL)$ is special. Hence, Proposition \ref{prop:twFut=A-muS} implies that 
\[
b(A_{(V,D)}(\ord_E)-\mu S(\ord_E))=\Fut_{1-\mu}(\cV,\cD;\cL)\geq 0. 
\]
This finishes the proof.
\end{proof}

\section{Twisted K-polystable degenerations}\label{sec:twistedKps}

In this section, we prove Theorem \ref{thm:kpolydegen} on the existence and uniqueness of twisted K-polystable optimal degenerations. Before stating it, we recall the definition of optimal degenerations. 

\begin{defn}
Let $(X,\Delta)$ be a log Fano pair with $\delta(X,\Delta)\leq 1$. A non-trivial special test configuration $(\cX,\Delta_{\tc};\cL)$ of $(X,\Delta)$ is called an \emph{optimal destabilization} if
$\Fut_{1-\delta(X,\Delta)}$ vanishes.
We call its central fiber $(\cX_0,\Delta_{\tc,0})$ an \emph{optimal degeneration} of $(X,\Delta)$.
\end{defn}

\begin{thm}[=Theorem \ref{thm:kpolydegen}]\label{thm:kpolydegen2}
 Let $(X,\Delta)$ be a log Fano pair that is not uniformly K-stable.
 If Conjecture \ref{conj:odc} holds for $(X,\Delta)$, then there exists a unique $\delta(X,\Delta)$-twisted K-polystable optimal degeneration of $(X,\Delta)$.
\end{thm}

The theorem may be viewed as a generalization of a result from \cite{LWX18} to the case when $(X,\Delta)$ is not necessarily K-semistable. 
Indeed, if  $(X,\Delta)$ is K-semistable, but not K-stable, the theorem implies
that there exists a test configuration $(\cX,\Delta_{\tc})$ of $(X,\Delta)$ with $\Fut=0$ and K-polystable special fiber (cf. \cite[Theorem 1.3]{LWX18}).

The key to proving the uniqueness of the twisted K-polystable optimal degeneration is to show that if $(\cX^{1},\Dtc^{1})$ and $(\cX^{2},\Dtc^{2})$ are optimal destabilizations of $(X,\Delta)$, then there exists a $\G_m^2$-equivariant family of log Fano pairs $(\fX,\Delta_{\fX}) \to\A^2$ such that
\[
(\fX,\Delta_{\fX})_{\A^1\times 1} \simeq (\cX^1,\Dtc^1) \quad \text{ and } \quad
(\fX,\Delta_{\fX})_{1\times \A^1} \simeq (\cX^2,\Dtc^2)
\]
(see Theorem \ref{c:familyoverA^2}).
Hence, the restrictions of $(\fX,\Delta_{\fX})_{0\times \A^1}$ and $(\fX,\Delta_{\fX})_{\A^1\times 0}$  are special degenerations  of  $(\cX^1_0,\Delta_{\tc, 0}^{1})$ and $(\cX^2_0,\Delta_{\tc, 0}^2)$ to a common log Fano pair $(\fX,\Delta_{\fX})_{0,0}$. 
In order to construct the above family over $\A^2$, 
we first prove a more general in Theorem \ref{t-thetared} concerning how test configurations degenerate along the germ of a curves.

The construction of the family $(\fX,\Delta_{\fX}) \to \A^2$ relies on utilizing finite generation results from the MMP. 
This construction of the family is carried out in \cite{LWX18}, when $(X,D)$  is K-semistable
using local arguments, in particular,  the cone construction and Chi Li's normalized volume function. 
The approach in this paper is global and uses the $\delta$-invariant. 

\subsection{Extension of test configurations}\label{s:thetareductivity}
In this section, we consider the following  problem: given a family of log Fano pairs over a pointed curve and a test configuration of the general fiber, does the test configuration extend to the entire family over the curve?
When the test configuration is an optimal degeneration, we show the existence of an extension (see Theorem \ref{t-thetared}).
This result is the key technical step in proving the uniqueness of twisted K-polystable degenerations in Theorem \ref{thm:kpolydegen2}.

\medskip

Let $R$ be the local ring of a smooth pointed curve, $K$ the fraction field of $R$, and $\kappa$ the residue field.
 Denote by $x$ the parameter for $\A^1$. Write $0\in \A^1$, $0_K \in \mathbb{A}^1_K$, and $0_\kappa\in \A^1_\kappa$ for the closed points defined by the vanishing of $x$.

Consider a $\Q$-Gorenstein family of log Fano pairs (see \cite[Definition 2.13]{ABHLX19}) $(X,\Delta) \to \Spec(R)$  and  a special test configuration
$(\cX_K, \Delta_{\cX_K})$  of $(X_K,\Delta_K)$.
From this data, we can construct a $\mathbb{G}_m$-equivariant $\Q$-Gorenstein family of log Fano pairs 
\begin{equation}\label{e-famXD}
f:(\cX^\circ ,\Delta_{\cX^\circ}) \to \mathbb{A}^1_R \setminus 0_\kappa
\end{equation}
defined by gluing the $\mathbb{G}_m$-equivariant families $(X,\Delta)\times(\mathbb{A}^1 \setminus 0)\to (\mathbb{A}^1 \setminus 0)_R$
and $(\cX_K,\Delta_{\cX_K})\to \mathbb{A}^1_K$ along their restrictions 
to $\mathbb{A}^1_K \setminus 0_K$. \footnote{The $\G_m$-action on $X\times(\mathbb{A}^1\setminus 0)$ is the product of the  trivial  action on $X$ and the standard action on $\mathbb{A}^1\setminus 0$.}

\begin{theorem}\label{t-thetared}
If $(X_\kappa,\Delta_\kappa)$ is $\mu$-twisted K-semistable and $\Fut_{1-\mu}(\cX_K,\Delta_{\cX_K}) = 0$, 
then $f^\circ$ extends uniquely to a $\mathbb{G}_m$-equivariant 
$\Q$-Gorenstein family of log Fano pairs
\[
{f}:({\cX},\Delta_{{\cX}}) \to \mathbb{A}^1_{R}
.\]
Furthermore, the fiber over $0_\kappa$
is $\mu$-twisted K-semistable and
$\Fut_{1-\mu}( { \cX}_\kappa , \Delta_{{\cX}_\kappa })=0$.
\end{theorem}

The following diagram shows the relationship between 
the different log Fano pairs that appear in the extension $f$. 
The symbol $\rightsquigarrow$ correspond to degeneration and the vertical ones are degenerations via the test configurations written on the side.
\begin{equation*}
\xymatrix@1 @R=1.1pc @C=0.9pc
{
(X_K, \Delta_K)
 \ar@{~>}^{}[rrr]
\ar@{~>}_{(\cX_K, \Delta_{\cX_K})/\A^1_K}[ddd]
 && & (X_\kappa,\Delta_\kappa) \ar@{~>}^{(\cX_\kappa, \Delta_{\cX_\kappa})/\A^1_\kappa}[ddd]
\\ \\ \\
(\cX_{0_K}, \Delta_{\cX_{0_K}}) \ar@{~>}_{}[rrr]
&&& (\cX_{0_\kappa},\Delta_{\cX_{0_\kappa}})
}
\end{equation*}

In the case when $\mu=1$, the result is precisely  \cite[Theorem 5.2]{ABHLX19}, which is an extension of results in \cite{LWX18}, and plays a key role in the construction of the moduli space of K-polystable Fano varieties as it verifies the $\Theta$-reductivity of  K-moduli stacks. 
In  \cite{BHLLX},  Theorem \ref{t-thetared} is used to show  that the moduli stack of log Fano pairs  admits a $\Theta$-stratification in the sense of \cite{HL18} assuming Conjecture \ref{conj:odc} holds.

To prove the theorem, we use the equivalence between test configurations and finitely generated filtrations of the section ring. 
Hence, proving the result will require us to show that the filtration of the section ring of $(X_K,\Delta_K)$ induced by  $(\cX_K,\Delta_{\cX_K})$
extends to a finitely generated filtration of the relative section ring of $(X,\Delta)$; see \cite[Section 5.2]{ABHLX19}.

\subsubsection{Extending filtrations induced by a divisor} 
Fix a positive integer $r$ so that $L:=-r(K_X+\Delta)$ is a very ample Cartier divisor. 
Let 
\[
V: = \bigoplus_{m \in \N}  V_m =
\bigoplus_{m \in \N}  H^0(X, \cO_{X}(mL)) 
\]
denote the section ring of $X$ with respect to $L$.
Write $V_\kappa= \bigoplus_m V_{\kappa,m}$ and $V_K= \bigoplus_m V_{K,m}$ for the restrictions of $V$ to $\Spec(\kappa)$ and $\Spec(K)$. Since $H^1(X, \cO_X(L)) =0$ by  \cite[Theorem 10.37]{Kollar13}, each $V_m$ is a  free $R$-module, $V_{K} = H^0(X_K,\cO_{X_K}(L_K))$, and $H^0(X_\kappa, \cO_{X_\kappa}(L_\kappa))$. 

Let $F_K$ be a divisor over $X_K$ 
that remains prime after base change to $\overline{K}$ 
and the consider the filtration
\[
\cF_K^p V_{K,m} : = \{ s \in V_{K,m}\, \vert \, \ord_{F_K}(s)\geq p \} .\]
The divisor $F_K$ extends to a divisor $F$ over $X$.  The filtration
$\cF_K$ extends to a filtration of $V_m$ by $R$-submodules
defined by 
\[
\cF^p V_m := \cF_K^p V_{K,m}\cap V_m
.\]
Note that $\cF^p V_m = \{s\in V_m \, \vert \, \ord_F(s)\geq p \}$,
since the vanishing of a function $f\in \cO_X$ along $F$ 
is equivalent to the vanishing of its restriction to $\cO_{X_K}$ along $F_K$.

When the $R[x]$-algebra 
\begin{equation}\label{eq:relativerees}
\bigoplus_{m \in \N}  \bigoplus_{p\in \Z}    \cF^p V_m  x^{-p}
\end{equation}
is of  finite type, 
we set $\cX:= \Proj_{\mathbb{A}^1_R} \left( \bigoplus_m \bigoplus_p   \cF^p V_m  x^{-p} \right)$, 
which has a $\mathbb{G}_m$-action given by the grading by $p$. 
The natural isomorphism
\[
\Big( \bigoplus_{m\in \N} \bigoplus_{p\in \Z} \cF^p V_m x^{-p} \Big) \otimes R[x^{\pm 1}]
\simeq
V \otimes_{R} R[x^{\pm1}]
\]
induces an isomorphism 
${\cX}\vert_{ \mathbb{A}^1\setminus 0}  \simeq X\times (\mathbb{A}^1\setminus 0)$. 
We write $\Delta_{\cX}$ for the  closure of $\Delta\times (\mathbb{A}^1\setminus 0)$ in  $\cX$.
Note that $\cX\to \A^1_R$ is flat, since $\bigoplus_p \cF^p V_m x^{-p}$ is 
a flat $R[x]$-module for each $m$ (see \cite[Section 3.1]{ABHLX19}).

\begin{prop}\label{p-extendval}
Fix $\mu \in (0,1]$.
If  $(X_\kappa,\Delta_\kappa)$ is $\mu$-twisted K-semistable and
$A_{X,\Delta}(F_K)- \mu S(F_K) =0$, 
then
\begin{enumerate}
    \item  \eqref{eq:relativerees} is a finite type $R[x]$-algebra;
    \item $( {\cX},\Delta_{{\cX}}) \to \mathbb{A}^1_R$ is a $\mathbb{G}_m$-equivariant $\Q$-Gorenstein family of log Fano pairs, the fiber over $0_\kappa$ is $\mu$-twisted K-semistable, and $\Fut_{1-\mu}({\cX}_\kappa ,\Delta_{{\cX}_\kappa} ) =0$.
\end{enumerate}
\end{prop}

\begin{proof}
Since $(X_{\kappa}, \Delta_\kappa)$ is $\mu$-twisted K-semistable, $\mu \leq \delta(X_{\kappa}, \Delta_{\kappa})$ by  Theorem \ref{thm:twisted=delta}.
Using  the lower semicontinuity of the stability threshold \cite{BL18b}
and our assumption on $F_K$, we see 
\[
\delta(X_\kappa, \Delta_\kappa) \leq \delta(X_{\overline{K}}, \Delta_{\overline{K}})
\leq 
\frac{A_{X_K,\Delta_K}(F_{K}) }{ S(F_K)}= \mu
.\]
Therefore, $ \delta(X_\kappa, \Delta_\kappa)=\mu\leq 1$.
\\

\noindent \emph{Claim 1}: For any $\epsilon>0$,
there exists $D_\epsilon  \in |-K_X-\Delta|_{\Q}$ so that
$(X,\Delta+ \mu D_\epsilon +X_\kappa)$ is lc and  
$A_{X,\Delta+\mu D_\epsilon +X_{\kappa}}(F)<\epsilon $.\\

Fix a positive integer $m$ and choose a basis $\{s_1,\ldots, s_{N_m} \}$ for $V_m$ compatible with the filtration $\cF  $. 
By the latter, we mean
$
\cF^p V_m = 
{\rm span} \{ s_1,\ldots, s_{d_p} \}$,
where $d_p:= {\rm rank} ( \cF^p V_m) $.
Consider the divisor
\[
B_m: 
=  
\frac{1}{mrN_m}
\big( \{s_1=0 \} +\cdots + \{s_{N_m} =0\} \big) 
 \in 
 |-K_X-\Delta|_{\Q}.\]
Note that $B_m$ restricts to an $mr$-basis type divisor on $(X_\kappa, \Delta_\kappa)$.
In addition, 
\[
\ord_{F}(B_m)
=
\ord_{F_K}(B_{m,K}) 
=
\frac{1}{mrN_m}
\sum_{p \geq 1} p \left( \dim \cF^p_K V_{m,K} - \dim \cF^{p+1}_K V_{m,K}\right)
=
S_{mr}(F_K)
, \]
where the second inequality follows from the fact that 
 $\{s_{1}, \ldots, s_{N_m}\}$   restricts to  a basis for $V_{m,K}$ that is compatible with the filtration $\cF_K$ 
 and the third from the equivalent definition of $S_{mr}$ in \cite[Section 3]{BJ17}.

Set $c_m : = \min \{ \delta_{mr}( X_{\kappa}, \Delta_\kappa),\delta(X_\kappa,\Delta_\kappa)\} $.  Since $c_m \leq \delta_{mr}(X_\kappa, \Delta_\kappa)$, the pair ${(X_\kappa, \Delta_\kappa+ c_m B_{m,\kappa})}$ is lc.
Therefore,
$(X, \Delta+ c_m B_m + X_\kappa)$ is lc  by inversion of adjunction. 
Now,  observe
\begin{align*}
A_{X,\Delta+c_m B_m +X_\kappa}(F) 
&= A_{X,\Delta+X_\kappa }(F) - c_m \ord_{F}(B_m) \\
&= A_{X_K,\Delta_K}(F_K)- c_m S_{mr}(F_K) .
\end{align*} 
Since $\lim_{m \to \infty} c_m  = \delta(X_\kappa,\Delta_\kappa)= \mu$ and
$\lim_{m \to \infty} S_{mr}(F_K)= S(F_K)$,  it follows that
\[
\lim_{m \to \infty} 
A_{X,\Delta+c_m B_m +X_\kappa}(F)
= 0.
\]
Hence, we may fix $m\gg0$ so that
\[A_{X,\Delta+c_{m} B_{m}  +X_\kappa}(F) < \epsilon .
\]
Using that 
 ${-K_{X}-\Delta}$ is ample, we may choose 
 $G\in |-K_X-\Delta|_{\Q}$ 
 so that 
\[
\left(X, \Delta+ c_{m} B_{m}  + (\mu-c_{m}) G + X_\kappa\right)
\]
remains lc. Therefore, the claim  holds with $D_\epsilon =\tfrac{1}{\mu} \left( c_{m}B_{m}+(\mu -c_{m})G\right)$.\\

Now, fix $0<\epsilon \ll1$. Set $D:= \mu D_\epsilon +(1-\mu)H \in |-K_X-\Delta|_{\Q}$, 
where $D_\epsilon$ satisfies the conclusion of Claim 1 and $H\in \frac{1}{r} |-r(K_X+\Delta)|$ is general.
Note that
\[
(X,\Delta+D+X_\kappa) \text{ is lc } \quad \text{ and  } \quad A_{X,\Delta+D+ X_\kappa}(F) <\epsilon.
\] 
By  \cite[Corollary 1.4.3]{BCHM10} (applied to the klt pair 
$(X,\Delta+(1-\epsilon') D_\epsilon +X_{\kappa})$ with $0<\epsilon'\ll 1$),
there  exists a proper birational morphism of normal varieties $Y\to X$ so that, $Y$ is $\Q$-factorial,  $F$ appears as a divisor on $Y$, and 
 ${\rm ExcDiv}(\mu) = F$.  
Note that $Y$ is Fano type over $\Spec(R)$ by  \cite[Lemma 2.9]{Zhu18}.
Therefore, $Y$ is a Mori dream space over $\Spec(R)$ \cite[Corollary 1.3.1]{BCHM10}. 
This implies       \eqref{eq:relativerees}  is a finite type $R[x]$-algebra and (1) is complete.
\\

\noindent \emph{Claim 2}: The pair $(\cX, \Delta_{\cX}+(1-\mu)\cH+ \cX_\kappa+ \cX_0)$ is lc, where $\cH$ is the closure of $H\times (\A^1\setminus 0)$ on $\cX$. 

Note that  $\cX$ is normal, $\cX_0$ is integral, and $K_{\cX}+\Delta_{\cX}$ is $\Q$-Cartier
by the  proofs of \cite[Lemma 3.8]{Fuj17b} and \cite[Lemma 2.6]{Fuj17b}, which extend   this relative setting.
Now, consider the proper birational morphism $Y_{\A^1}\to X_{\A^1}$ constructed by taking the product of $Y\to X$ with $\A^1$.
Note that the pair
\begin{equation}\label{e:CYpair}
(X_{\A^1},\Delta_{\A^1}+ D_{\A^1}+X_{\kappa, \A^1}+X\times \{0 \})
\end{equation}
 is log Calabi--Yau, meaning the pair is lc and $K_{X_{\A^1}}+\Delta_{\A^1}+D_{\A^1}+X_{\kappa,\A^1} +X\times \{0 \}\sim_{\Q} 0$.
Let $E$ be the divisor over $X_{\A^1}$ so that $\ord_{E}$ is the quasimonomial valuation with weights $(1,1)$ along $F\times \A^1$ and $Y\times \{0\}$. 
Since
\[
a_\epsilon : = A_{X_{\A^1},\Delta_{\A^1}+ D_{\A^1}+X_{\kappa,\A^1}+X\times \{0 \}} (E) = A_{X,\Delta+ D+X_\kappa}(F)< \epsilon <1
,\]
 we may find  a proper birational morphism of normal varieties $\rho: \cY\to X_{\A^1}$ 
so that ${\rm ExcDiv}(\rho) =E$  \cite[Corollary 1.4.3]{BCHM10}. Therefore,
\begin{equation}\label{e:CYpair2}
(\cY, \rho_*^{-1} (\Delta_{\A^1}+ D_{\A^1}+X_{\kappa, \A^1}+X\times \{0 \}) + (1-a_\epsilon) E) 
\end{equation}
is log Calabi--Yau, since it is the log pullback \eqref{e:CYpair}.

Now, consider the birational map $\tau:\cY \dashrightarrow \cX$.
Since $\tau$ is an isomorphism over $\A^1\setminus 0$, $\tau_*E = \cX_0$ and $\tau_*(\rho_*^{-1}(X\times \{0 \}) ) = 0$,
$\mathcal{Y}\dashrightarrow \mathcal{X}$ is a birational contraction. 
This implies
\[
( \cX, \Delta_{\cX}+  \cD + \mathcal{X}_\kappa + (1- a_\epsilon ) \mathcal{X}_0) 
\]
is log Calabi--Yau, where $\cD$ is the closure of $D\times (\A^1\setminus 0)$ on $X$,
since it is the  push forward of the log Calabi--Yau pair \eqref{e:CYpair2} under a birational contraction. 
Using that $(1-\mu) \cH \leq \cD$, we see
 \[
 (\cX, \Delta_{\mathcal{X}}+(1-\mu)\cH+ \mathcal{X}_\kappa+(1-a_\epsilon )\mathcal{X}_0)\]
is lc.
Since $\lim_{\epsilon \to 0} a_{\epsilon }  =0$, we conclude 
$( \cX, \Delta_{\cX}+(1-\mu) \cH+ \mathcal{X}_\kappa +\mathcal{X}_0)$ is lc, which completes the claim.\\

Now, consider the $\G_m$-linearized line bundle  $\cL:= \cO_{\cX}(d)$, where $d$ is a sufficiently divisible positive integer, and   the test configuration $(\cX_K, \Delta_{\cX_K}; \cL_K)$.
Since $\cX_K$
is equivariantly isomorphic to $\Proj  \big( \bigoplus_{m,p} \cF^p_K V_{m,K}x^{-p}\big)$,
 $\cX_{0_K}$ is integral and $v_{\cX_{K,0_K}} = \ord_{F_K}$ \cite[Section 3.2]{Fuj17b}. Therefore,  
Proposition \ref{prop:twFut=A-muS} implies 
\[\Fut_{1-\mu}( \cX_K, \Delta_{\cX_K};\cL_K)=A_{X_K,\Delta_K}(F_K)- \mu S(F_K)=0.\]

We now analyze the restriction of $(\cX, \Delta_{\cX}; \cL)$
to  $\Spec(\kappa)$. 
Since  ${(\cX,\Delta_{\cX}+(1-\mu) \cH+\cX_\kappa+\cX_0)}$ is lc, $(\cX_\kappa,\Delta_{\cX_\kappa} +(1-\mu) \cH_\kappa +\cX_{0_\kappa})$ is slc by adjunction. 
The latter implies   $\Supp(\Delta_{\cX_{\kappa}} +(1-\mu)\cH_{\kappa})$ does not contain any irreducible component of $\cX_{0_\kappa}$
and, using \cite[Proposition 2.6]{BHJ17}, $\cX_\kappa$ is normal.  
Since  $H\in \frac{1}{r}|-r(K_X+\Delta)|$ is general, we may assume  $H_K$ and $H_\kappa$ are  compatible with $\cX_K$ and $\cX_{\kappa}$ by Lemma \ref{l-compatgeneral}.
Therefore, 
\begin{align*}
0=\Fut_{1-\mu}( \cX_K, \Delta_{\cX_K};\cL_K)
&=
\Fut( \cX_K, \Delta_{\cX_K} +(1-\mu) \cH_K; \cL_K)\\
&=\Fut( \cX_\kappa, \Delta_{\cX_\kappa} +(1-\mu) \cH_\kappa ; \cL_\kappa )\\
& = \Fut_{1-\mu}( \cX_\kappa, \Delta_{\cX_\kappa} ; \cL_\kappa)
\end{align*}
where the third equality uses that intersection numbers of line bundles are locally constant in flat proper families.
 This implies  $(\cX_\kappa, \Delta_{\cX_\kappa})$  is special   (Corollary \ref{cor:twisted-special}) and, hence, all  geometric fibers of $(\cX , \Delta_{\cX}) \to \A^1_R$  are log Fano pairs.
 In addition, the fiber over $0_\kappa$ is $\mu$-twisted K-semistable  (Lemma \ref{lem:lwx3.1}), which completes the proof.
\end{proof}

\subsubsection{Proof of Theorem  \ref{t-thetared}}
We will now deduce Theorem \ref{t-thetared} from 
Proposition \ref{p-extendval}.

\begin{proof}[Proof of Theorem \ref{t-thetared}]
Fix  a positive integer $b$ and 
a divisor $F_K$ over $X_K$ so that $r(\ord_{\cX_{K,0_K}})$
equals $b \ord_{F_K}$.
If we write $\cF_{K}^\bullet V_K$ for  the filtration of $V_K$
induced by $E_K$, then
$\cX_K$ is equivariantly isomorphic to  $\Proj_{\mathbb{A}^1_K}( \bigoplus_{m,p} \cF_K^{ \lceil p/b\rceil} V_{K,m} x^{-p} )$ by \cite[Section 3.2]{Fuj17b}. 

Now, note that  Proposition \ref{prop:twFut=A-muS} implies
\[\Fut_{1-\mu}(\cX_K,\Delta_{\cX_K})=b (A_{X_K,\Delta_K}(F_K) - \mu S(F_K)),\] which is zero by assumption.
Therefore, Proposition \ref{p-extendval} (1) implies $ \bigoplus_{m,p} \cF^p V_m x^{-p}$, where 
$\cF^p V_m : = \cF^p_K V_{K,m}\cap V_m$, is a finite type $R[x]$-algebra. Hence, $ \bigoplus_{m,p} \cF^{\lceil p/b \rceil} V_mx^{-p} $ is a finite type $R[x]$-algebra as well. 

If we set $\cX := \Proj_{\mathbb{A}^1_R} ( \bigoplus_{m,p}  \cF^{\lceil p/b \rceil} V_m x^{-p})$ and write $\Delta_{{\cX}}$ for the closure of $\Delta\times (\mathbb{A}^1 \setminus 0)$ in $ {\cX}$, then  
 $f:(\cX,\Delta_{{\cX}})\to \mathbb{A}^1_R$ is an extension of the family $f^\circ$; see  \cite[Section 5.2]{ABHLX19}.
 Since  $({\cX},\Delta_{{\cX}})\to \mathbb{A}^1_R$ is a base change   of the family  in Proposition \ref{p-extendval} by the map $\mathbb{A}^1_R\to \mathbb{A}^1_R$ sending $x\mapsto x^b$, $f$ is a $\Q$-Gorenstein family of log Fano pairs and the fiber over $0_\kappa$ is $\mu$-twisted K-semistable.
\end{proof}

\subsubsection{Applications of  Theorem \ref{t-thetared}}\label{ss:consequencethetared}
Using Theorem \ref{t-thetared}, we show that certain certain classes of test configurations can be connected via a family over $\A^2$.

Let $(X,\Delta)$ be a log Fano 
pair and $(\cX^{i},\Delta_{\tc}^{i})$ for $i=1,2$ be special test configurations of $(X,\Delta)$. 
From this data, we can construct\footnote{
This is constructed from gluing the products
$(\cX^{i}, \Delta_{\tc}^i) \times ({\mathbb{A}^1\setminus 0}) \to \mathbb{A}^1 \times (\mathbb{A}^1 \setminus 0) $ for $i=1,2$
along their restrictions to $(\mathbb{A}^1 \setminus 0)^2$.
} a $\mathbb{G}_m^2$-equivariant $\Q$-Gorenstein family of log Fano pairs 
\begin{equation}\label{e-familyA^2minus00}
g^\circ :(\mathfrak{X}^\circ, \Delta_{\mathfrak{X}^\circ})\to \mathbb{A}^2 \setminus (0,0),\end{equation}
so that
\begin{enumerate}
\item $\G_m^2$ acts on $\A^2\setminus (0,0)$ via $(s,t)\cdot (x,y) = (sx,ty)$,
    \item there is a fixed isomorphism $(\fX^\circ ,\Delta_{\fX^\circ})_{(1,1)} \simeq (X,D)$,
    \item the restrictions 
    \[
    (\mathfrak{X}^\circ, \Delta_{\mathfrak{X}^\circ})\vert_{\mathbb{A}^1 \times \{1 \} }
    \quad \text{ and } \quad 
    (\mathfrak{X}^\circ, \Delta_{\mathfrak{X}^\circ})\vert_{\{1\} \times \mathbb{A}^1  },\]
      with  induced $\G_m$ actions given by  $\G_m\times \{1\}$  and $\{1 \} \times \G_m$, are isomorphic as a test configurations to $(\cX^{1}, \Delta_{\tc}^{1})$ and $(\cX^{2}, \Delta_{\tc}^{2})$.
\end{enumerate}

We seek to show $g^\circ$ extends to a
$\mathbb{G}_m^2$-equivariant $\Q$-Gorenstein family  of log Fano pairs
${g:({\mathfrak{X}},\Delta_{ {\mathfrak{X}}})\to \mathbb{A}^2}$ under certain hypotheses. 
When such a ${g}$ exists, the restrictions  $({\mathfrak{X}},\Delta_{ {\mathfrak{X}}})_{\mathbb{A}^1\times \{0 \}}$
and $({\mathfrak{X}},\Delta_{ {\mathfrak{X}}})_{\{0 \}\times \A^1}$, 
 are naturally special test configurations of $(\cX^2_0,\Delta_{\tc,0}^2)$
 and $(\cX^1_0,\Delta_{\tc,0}^1)$
 and have the same central fiber.
 
The following statement was proven in  \cite[Section 3]{LWX18}
when $\mu=1$.

\begin{theorem}\label{c:familyoverA^2}
If $(\cX^1_0, \Delta_{\tc,0}^1)$
is $\mu$-twisted K-semistable and $\Fut_{1-\mu}(\cX^{2}, \Delta_{\tc}^2)=0$, 
then $g^\circ $ extends uniquely to 
a $\mathbb{G}_m^2$-equivariant $\Q$-Gorenstein family  of log Fano pairs 
$${g}:({\mathfrak{X}},\Delta_{ {\mathfrak{X}}})\to \mathbb{A}^2.$$
Furthermore, 
the special test configuration  $
({\mathfrak{X}},\Delta_{ {\mathfrak{X}}})_{\{0 \}\times \mathbb{A}^1}$
has vanishing $\Fut_{1-\mu}$ and the fiber over $(0,0)$ is $\mu$-twisted K-semistable. 
\end{theorem}

 We will deduce the result from Theorem \ref{t-thetared}.
Throughout the proof, we rely on the elementary fact that if $o\in S$ is a closed point on a regular surface, then there is a bijection between vector bundles on $S$ and $S\setminus o$ 
given by pushing forward and pulling back.

\begin{proof} 
Fix a positive integer $r$ so that $\mathfrak{L}^\circ : =-r(K_{\mathfrak{X}^\circ / ( \mathbb{A}^2\setminus (0,0)} + \Delta _{\mathfrak{X}^\circ})$
is a Cartier divisor. 
By \cite[Lemma 2.15]{ABHLX19}, 
$(\mathfrak{X}^\circ , \mathfrak{L}^\circ)$ extends to a  polarized family of schemes over $\mathbb{A}^2$ with a $\mathbb{G}_m^2$-action if and only if 
$$ \mathcal{A}= 
\bigoplus_{m\in \N}
  H^0( \mathbb{A}^2, i_* g^\circ _* \cO_{\mathfrak{X}^\circ }(m\mathfrak{L}^\circ))
 $$ is a finite type $\cO_{\mathbb{A}^2}$-algebra, where $i: \mathbb{A}^2 \setminus (0,0) \to \mathbb{A}^2$ is the inclusion map, 
 and the latter holds if an only if $\cA\vert_{(0,0)}$
is a finite type $k$-algebra. 
Additionally, if the extension exists, then 
the extension is given by $(\Proj_{\mathbb{A}^2}(\cA) , \cO(1))$. 

Set $R: = \cO_{\mathbb{A}^1,0}$ and consider the natural map
$\mathbb{A}^1_R : = \Spec(R) \times  \mathbb{A}^1\to \mathbb{A}^2$, which induces a map
$\mathbb{A}^1_R \setminus 0_\kappa \to \mathbb{A}^2 \setminus (0,0)$. Restricting $(\mathfrak{X}^\circ, \Delta_{\mathfrak{X}^\circ}) $ via the last  map, gives 
a $\Q$-Gorenstein family of log Fano pairs $f:(\cX^\circ , \Delta_{\cX^\circ }) \to   \mathbb{A}^1_R \setminus 0_\kappa $.
Set $\cL^\circ :=\mathfrak{L}^\circ \vert_{\cX^\circ}$
and
$$\cB: = \bigoplus_{m \in \N} H^0(\mathbb{A}^1_R, j_* f_*^\circ \cO_{\cX^\circ}(m \cL^\circ))
,$$ where $j : \mathbb{A}^1_R \setminus 0_\kappa \to \mathbb{A}^1_R$ is the natural inclusion. 
Observe that $\cB^\circ \simeq  \cA\vert_{\mathbb{A}^1_R}$  and, hence, $\cB\vert_{0_\kappa} \simeq \cA\vert_{0,0}$.

Note that the family  $f^\circ$ can be constructed as in Section \ref{s:thetareductivity} 
using the family $(\cX^{1}, \Delta_{\tc}^1)\times _{\mathbb{A}^1} \Spec(R) \to \Spec(R)$ and 
the test configuration $(\cX^{2}, \Delta_{\tc}^2)\times \Spec(K) \to \mathbb{A}^1_K$ of $(X,\Delta)\times K$.
Therefore,  Theorem \ref{t-thetared} implies $f^\circ$ extends to a $\Q$-Gorenstein family of log Fano pairs $({\cX}, \Delta_{{\cX}})\to \mathbb{A}^1_R$, 
$\Fut_{1-\mu}({\cX}_\kappa , \Delta_{{\cX}_\kappa})=0 $,
and the fiber over $0_\kappa$ is $\mu$-twisted K-semistable.
Since the family extends,
$\cB$ is finite type  $\cO_{\mathbb{A}^1_R}$-algebra and ${\cX} \simeq \Proj_{\mathbb{A}^1_R}(\cB)$ by  \cite[Lemma 2.15]{ABHLX19}.  Therefore, $\cB\vert_{0_\kappa}\simeq \cA\vert_{0,0}$ are finite type as well and so is $\cA$. 
Hence, $g$ extends to a projective family ${g}:{\mathfrak{X}}\to \mathbb{A}^2$.

Finally, write $\Delta_{{\mathfrak{X}}}$ 
for the $\Q$-divisor on  ${\mathfrak{X}}$ that is the componentwise closure of $\Delta_{\mathfrak{X}}$
and observe that  
$({\cX}, \Delta_{{\cX}}) \simeq
({\mathfrak{X}}, \Delta_{\mathfrak{X}}) \vert_{\mathbb{A}^1_R}$. 
Therefore, $({\mathfrak{X}}, \Delta_{\mathfrak{X}}) \to \mathbb{A}^2$ 
is a $\Q$-Gorenstein family of log Fano pairs
and 
$({\mathfrak{X}}, \Delta_{{\mathfrak{X}}})\vert_{ \{0 \} \times \mathbb{A}^1}$ 
is naturally a test configuration with vanishing $\Fut_{1-\mu}$. 
\end{proof}

\begin{remark}\label{rem:familyoverA^2}
We describe additional properties of the family ${g}:({ \mathfrak{X}}, \Delta_{{ \mathfrak{X}}})\to \A^2$ constructed in Theorem \ref{c:familyoverA^2}.
\begin{enumerate}
    \item If we add  to the assumptions of Theorem \ref{c:familyoverA^2}  that $\Fut_{1-\mu}(\cX^1,\Delta_{\tc}^1) =0$, 
then the restriction  $({ \mathfrak{X}}, \Delta_{
{ \mathfrak{X}}}) \vert_{\mathbb{A}^1 \times \{0\}}$, with the $\G_m$-action given by $\G_m\times \{1\}$,
 also has vanishing $\Fut_{1-\mu}$. 
This is a consequence of the uniqueness of the extension ${g}$ and repeating the argument with the test configurations switched.

\item The test configurations $({ \mathfrak{X}}, \Delta_{{ \mathfrak{X}}})_{\A^1\times \{0 \}}$
and $({ \mathfrak{X}}, \Delta_{{ \mathfrak{X}}})_{\{0\} \times \A^1}$ are equivariant with respect to the $\mathbb{G}_m$-actions on $(\cX_{0}^i,\Delta_{\tc,0}^i)$. 
This follows from the fact that the actions by  $\mathbb{G}_m \times \{1\}$ and $\{1\}\times \G_m$ on $({ \mathfrak{X}}, \Delta_{{ \mathfrak{X}}})$ commute.

\item If $(X,\Delta)$ admits a $\T:= \G_m^r$-action and $(\cX^i,\Delta_{\tc}^i)$ for $i=1,2$ are $\T$-equivariant test configurations, then the $\T$-action on $$(\mathfrak{X}^\circ,\Delta_{\mathfrak{X}^\circ})\vert_{(\A^1 \setminus 0)^2} \simeq (X,\Delta) \times (\A^1\setminus0)^2,$$
given by the action of $\T$ on $(X,\Delta)$ and the trivial action on 
$(\mathbb{A}^1\setminus0)^2,$ extends to a $\T$-action on $(\mathfrak{X},\Delta_{\mathfrak{X}})$ commuting with the $\G_m^2$-action on $({X},\Delta_{\mathfrak{X}})$. Furthermore, this $\T\times \G_m^2$-action extends to $({\mathfrak{X}},\Delta_{{\mathfrak{X}}})$ by  the proof of \cite[Lemma 2.15]{ABHLX19}. 

Now, note that 
$(\cX^i_0,\Delta_{\tc,0}^i)$ for $i=1,2$
admit $\G_m \times \T $-actions, 
since they 
are the special fibers of $\T$-equivariant test configurations. 
The previous paragraph implies that the test configurations
$({ \mathfrak{X}}, \Delta_{{ \mathfrak{X}}})_{\{0\} \times \A^1}$  and 
$({ \mathfrak{X}}, \Delta_{{ \mathfrak{X}}})_{  \A^1\times \{0\}}$ are $\G_m\times \T$ equivariant.
\end{enumerate}
\end{remark}

\subsection{Proof of Theorem \ref{thm:kpolydegen2}}
We will proceed to use Theorem \ref{c:familyoverA^2} to deduce Theorem \ref{thm:kpolydegen2}.
First, we prove the following propositions that will be used to show the uniqueness and existence of the  twisted K-polystable optimal degeneration.

\begin{prop}\label{prop:uniquenessKps}
Let $(X,\Delta)$ be a $\mu$-twisted K-semistable log Fano pair. 
Assume 
  $(\cX^{1},\Delta_{\tc}^1)$
 and $(\cX^{2},\Delta_{\tc}^2)$ are
 special test configurations of $(X,\Delta)$
with vanishing $\Fut_{1-\mu}$. If
$(\cX^{1}_0,\Delta_{\tc,0}^{1})$ and
$(\cX^{2}_0,\Delta_{\tc,0}^{2})$ are  $\mu$-twisted K-polystable, then they are isomorphic.
\end{prop}

\begin{proof}
Following Section \ref{ss:consequencethetared},
we consider the $\G_m^2$-equivariant $\Q$-Gorenstein family of log Fano pairs $g^\circ:({\mathfrak{X}}^\circ, \Delta_{{\mathfrak{X}^\circ}})\to \mathbb{A}^2 \setminus (0,0)$
whose restrictions to $\mathbb{A}^1\times \{1 \}$ and $\{1\}\times \mathbb{A}^1$
are isomorphic to the special test configurations $(\cX^{1},\Delta_{\tc}^1)$
and $(\cX^{2},\Delta_{\tc}^2)$. 
 Theorem  \ref{c:familyoverA^2}  (with Remark \ref{rem:familyoverA^2} (1)) implies $g^\circ$ extends to a family ${g}: ({\mathfrak{X}}, \Delta_{{\mathfrak{X}}})\to \mathbb{A}^2$
and the restrictions
$$
({\mathfrak{X}}, \Delta_{{\mathfrak{X}}})\vert_{ 
\mathbb{A}^1 \times \{0\}}
\quad 
\text{ 
and }
\quad 
({\mathfrak{X}}, \Delta_{{\mathfrak{X}}})\vert_{ 
\{0 \} \times 
\mathbb{A}^1}$$ are special test configurations of 
$(\cX^{1}_0,\Delta_{\tc,0}^{1})$ and
$(\cX^{2}_0,\Delta_{\tc,0}^{2})$ with
vanishing $\Fut_{1-\mu}$. 

Since  $(\cX^{1}_0,\Delta_{\tc,0}^{1})$ and
$(\cX^{2}_0,\Delta_{\tc,0}^{2})$  are    $\mu$-twisted K-polystable by assumption, the test configurations constructed above must be product type. Therefore, $(\cX^{1}_0,\Delta_{\cX^{1},0})\simeq
(\mathfrak{X}, \Delta_{\mathfrak{X}})\vert_{ 
(0,0)}\simeq (\cX^{2}_0,\Delta_{\cX^{2},0})$ and the proof is complete. 
\end{proof}

\begin{prop}\label{lem:Fut0=>equivariant}
Let $(X,\Delta)$ be a $\mu$-twisted K-semistable log Fano pair with a $\mathbb{T}=\mathbb{G}_m^r$ action. 
If $(X,\Delta)$ is not $\mu$-twisted K-polystable, then there exists a $\mathbb{T}$-equivariant  special test configuration $(\cX, \Delta_{\tc})$
such that $\Fut_{1-\mu}(\cX, \Delta_{\tc})=0$ and $(\cX_0, \Delta_{\tc,0}) \not\simeq (X,\Delta)$.
\end{prop}

\begin{proof}
Since $(X,\Delta)$ is $\mu$-twisted K-semistable, but not polystable, 
Corollary \ref{cor:twisted-special} implies there exists a non-product type
special test configuration $(\cX,\Delta_{\tc})$ of 
$(X,\Delta)$ with vanishing $\Fut_{1-\mu}(\cX,\Delta_{\tc})$. 
The non-product type condition implies $(\cX_0,\Delta_{\tc,0})\not\simeq (X,\Delta)$.

We seek to replace $(\cX,\Delta_{\tc})$ with a $\T$-equivariant test configuration.
Arguing by induction, fix $0 \leq i < r$ and assume $(\cX,\Delta_{\tc})$ 
is equivariant with respect to  the $ \mathbb{G}^i$-action on $X$ determined by
the first $i$-copies of $\mathbb{G}_m$ in $\mathbb{T}$. 
It suffices to show that we can replace $(\cX,\Delta_{\tc})$ with a $\mathbb{G}_m^{i+1}$-equivariant test configuration with vanishing $\Fut_{1-\mu}$ and special fiber not isomorphic to $(X,\Delta)$. 

To proceed, let $\lambda: \mathbb{G}_m\to \Aut(X,\Delta)$ denote the 1-parameter subgroup determined by the $1+i$-th copy of $\mathbb{G}_m$ in $\mathbb{G}_m^r$ and write 
$(\cX_\lambda , \Delta_{\tc,\lambda})$
for the corresponding  product type test configuration. 
Following Section \ref{ss:consequencethetared}, 
we consider the $\mathbb{G}_m^2$-equvariant $\Q$-Gorenstein family of log Fano pairs $g:({\mathfrak{X}}, \Delta_{{\mathfrak{X}}})\to \mathbb{A}^2 \setminus (0,0)$ such that 
$$({\mathfrak{X}}, \Delta_{{\mathfrak{X}}})_{\mathbb{A}^1\times \{1 \}}
 \simeq (\cX_\lambda,\Delta_{\tc,\lambda})
 \quad \text{ and }
 \quad
 ({\mathfrak{X}}, \Delta_{{\mathfrak{X}}})_{\{1\}\times \mathbb{A}^1}\simeq(\cX,\Delta_{\tc}).$$
By Theorem \ref{c:familyoverA^2}, $g$ extends to a $\mathbb{G}_m^2$-equivariant family of log Fano pairs  $\widetilde{g}: (\widetilde{\mathfrak{X}}, \Delta_{\widetilde{\mathfrak{X}}}) \to \mathbb{A}^2$ and $(\cX',\Delta'_{\tc}):= (\widetilde{\mathfrak{X}},\Delta_{\widetilde{\mathfrak{X}}})_{\{0\} \times \mathbb{A}^1}$ is naturally a special test configuration with  vanishing $\Fut_{1-\mu}$.
Furthermore, since $(\cX_{\lambda},\Delta_{\tc,\lambda})$
and $(\cX,\Delta_{\tc})$-are $\G_m^i$-equivariant, 
Remark \ref{rem:familyoverA^2} (3)
implies $(\cX',\Delta'_{\tc})$ is $\G_m^{i+1}$-equivariant.
Note that $(\cX'_0,\Delta'_{\tc,0})\not \simeq (X,\Delta)$, since  there exist degenerations $(X,\Delta) \leadsto (\cX_0,\Delta_{\tc,0}) \leadsto (\cX'_0,\Delta'_{\tc,0})$. 
Hence, the proof is complete.
\end{proof}

\begin{proof}[Proof of Theorem \ref{thm:kpolydegen2}]
Since $(X,\Delta)$ is not uniformly K-stable, $\delta: = \delta(X,\Delta)\leq 1$.
By Theorem \ref{thm:main-twisted=delta}, 
 $(X,\Delta)$ is $\delta$-twisted K-semistable.
If the pair admits an optimal destabilization $(X,\Delta) \leadsto (X_0,\Delta_0)$ such that $(X_0,\Delta_0)$ is $\delta$-twisted K-polystable, then $(X_0,\Delta_0)$ is unique by Proposition \ref{prop:uniquenessKps}. 

We are left to prove the existence of such a degeneration. 
Theorem \ref{thm:main1} combined with the assumption that Conjecture \ref{conj:odc} holds
 implies the existence of a
test configuration $(\cX,\Delta_{\tc})$
of $(X,\Delta)$ such that $\Fut_{1-\delta}(\cX, \Delta_{\tc})=0$. 
Note that $(X_0,\Delta_0):=(\cX_0,\Delta_{\tc,0})$ has a $\mathbb{G}_m$-action given by restricting the $\G_m$-action on $(\cX,\Delta_{\tc})$ and, by Lemma \ref{lem:lwx3.1}, 
is $\delta$-twisted K-semistable. 

If $(X_0,\Delta_0)$ is $\delta$-twisted K-polystable, we are done. 
If not,   apply  Proposition  \ref{lem:Fut0=>equivariant} to find a $\mathbb{G}_m$-equivariant special test configuration $(\cX',\Delta'_{\tc})$ of $(X_0,\Delta_0)$
with $\Fut_{1-\delta}(\cX', \Delta'_{\tc})=0$ and special fiber $(X_1,\Delta_1)$
not equal to $(X_0,\Delta_0)$. 
 Hence, $(X_1,\Delta_1)$ admits a $\mathbb{G}_m^2$-action. Again, the pair is 
is $\delta$-twisted K-semistable  by Lemma \ref{lem:lwx3.1}. If $(X_1,\Delta_1) $ is not $\delta$-twisted K-polystable, we can apply  Proposition  \ref{lem:Fut0=>equivariant}
to find a $\G_m^2$-equivariant non-product type special test configuration of $(X_1,\Delta_1)$ with vanishing $\Fut_{1-\delta}$. 
This process must  terminate after $r \leq \dim(X)$ steps with a $\delta$-twisted K-polystable pair, since the $i$-th degeneration
$(X_i,\Delta_i)$ will be admit an effective $\mathbb{G}_m^i$-action. 
Hence, we get a sequence of special degenerations $$(X,\Delta)\leadsto (X_0,\Delta_0) \leadsto \cdots \leadsto (X_r,\Delta_r)$$  with vanishing $\Fut_{1-\delta}$ such that $(X_r,\Delta_r)$ is $\delta$-twisted K-polystable. 

We seek to construct a  special test configuration with  $\Fut_{1-\delta}=0$ and degenerating  $(X,\Delta)\leadsto (X_r,\Delta_r)$ (for this we use an argument from the proof of Lemma \ref{lem:lwx3.1}).
Fix a positive integer $m$ such that $-m(K_X+\Delta)$ is a Cartier divisor and $|-m(K_X+\Delta)|$ is base point free. 
Since $(\cX',\Delta'_{\tc})$ is equivariant with repsect to the $\G_m$-action on $(X_0,\Delta_0)$, 
the proof of \cite[Lemma 3.1]{LWX18} implies
there exists a positive integer $k$ and a special test configuration  $(\cX'',\Delta''_{\tc})$ of $(X,\Delta)$ such that the following hold:
\begin{itemize}
    \item $(X_1,\Delta_1)\simeq(\cX''_0,\Delta''_{\tc,0})$;
    \item $\Fut(\cX'', \Delta_{\tc}''+(1-\delta)\cD'')
=
k\cdot \Fut(\cX', \Delta_{\tc}'+(1-\delta)\cD')
+
\Fut(\cX, \Delta_{\tc}+(1-\delta)\cD)
$ for a general $D\in \frac{1}{m} |-m(K_X+\Delta)|$.
\end{itemize}
Since $D\in \frac{1}{m} |-m(K_X+\Delta)|$ is general, we may assume it is compatible with $\cX$ and $\cX''$. 
Hence, 
\[
 \Fut_{1-\delta}(\cX'', \Delta_{\tc}'')
=
k\cdot \Fut(\cX', \Delta_{\tc}'+(1-\delta)\cD')
+
\Fut_{1-\delta}(\cX, \Delta_{\tc})
\leq 0.\]
Since $(X,\Delta)$ is $\delta$-twisted K-semistable, we conclude $ \Fut_{1-\delta}(\cX'', \Delta_{\tc}'';\cL'')=0$.  
Repeating this argument  $r-1$ additional times gives a special test configuration with $\Fut_{1-\delta}=0$ and degenerating $(X,\Delta) \leadsto (X_r,\Delta_r)$.
\end{proof}

\section{Applications}\label{sec:greatestRic}

In this section, we prove the following result on stability thresholds in families from which Theorem \ref{thm:greatestRic} follows easily. 

\begin{thm}\label{thm:odc-openness}
 Let $\pi:(\cX,\Delta_T)\to T$ be a $\bQ$-Gorenstein flat family of log Fano pairs over a normal base $T$.
 Assume Conjecture \ref{conj:odc} is true for all K-unstable fibers $(\cX_t,\Delta_t)$
 of $\pi$. Then the function $T\ni t\mapsto \min\{1,\delta(\cX_t,\Delta_t)\}$
 is a lower semi-continuous constructible rational-valued function.
 In particular, the locus $\{t\in T\mid (\cX_t,\Delta_t)\textrm{ is K-semistable}\}$
 is a Zariski open subset of $T$.
\end{thm}

\begin{remark}
Note that in \cite{BLX19} it is shown that $t\mapsto\min\{1,\delta(\cX_t,\Delta_t)\}$ is a constructible function without assuming Conjecture \ref{conj:odc}. Moreover, the openness of K-semistability is confirmed in \cite{Xu19, BLX19}. Since the approaches from \cite{Xu19, BLX19} are quite different from this paper, we provide a proof of Theorem \ref{thm:odc-openness} as it may have independent interest.
\end{remark}

We first start with a few intermediate results.

\begin{prop}\label{prop:bdd-test-twistedK}
Let $\pi:(\cX,\Delta_T)\to T$ be a $\bQ$-Gorenstein flat family of log Fano pairs over a connected normal base $T$. Assume Conjecture \ref{conj:odc} is true for all K-unstable fibers $(\cX_t,\Delta_t)$ of $\pi$. Then there exists a positive integer $m$ such that the following hold:
\begin{enumerate}
    \item The divisor $-m(K_{\cX/T}+\Delta_T)$ is very ample over $T$. Moreover, $N:=h^0(\cX_t,\cO_{\cX_t}(-m(K_{\cX_t}+\Delta_t)))-1$ is independent of the choice of $t\in T$. Thus we have an embedding $(\cX_t,\Delta_t)\hookrightarrow \bP^N$ for each $t\in T$.
    \item If $(\cX_t,\Delta_t)$ is $\mu$-twisted K-unstable for some $t\in T$ and $\mu\in (0,1]$, then there exists a $1$-PS in $\SL(N+1)$ whose induced $\mu$-twisted generalized Futaki invariant is negative.
\end{enumerate}
\end{prop}

\begin{proof}
The proof of (1) is straightforward since $-(K_{\cX/T}+\Delta_T)$ is a $\pi$-ample $\bQ$-Cartier $\bQ$-divisor and $\pi$ is flat. Let us take a positive integer $m_1$ so that all multiples $m$ of $m_1$ satisfies (1). We now prove (2).
Firstly, by \cite[Proposition 5.3]{BL18b}, there exists a positive constant $c_1$ such that $\alpha(\cX_t,\Delta_t)\geq c_1$ for any $t\in T$. Assume that $(\cX_t,\Delta_t)$ is $\mu$-twisted K-unstable for some $\mu\in (0,1]$. Hence we know that it is K-unstable as well. By assumption, Theorem \ref{thm:main1} and \ref{thm:twisted-iff}, there exists a non-trivial special test configuration $(\cY_t,\Gamma_t)/\bA^1$ of $(\cX_t,\Delta_t)$ such that $\delta(\cY_{t,0},\Gamma_{t,0})=\delta_t$ and $\Fut_{1-\delta_t}(\cY_t,\Gamma_t)=0$ where $\delta_t:=\delta(\cX_t,\Delta_t)$. Hence we have
\[
\alpha(\cY_{t,0},\Gamma_{t,0})\geq (n+1)^{-1}\delta(\cY_{t,0},\Gamma_{t,0})=(n+1)^{-1}\delta_t\geq (n+1)^{-1}\alpha(\cX_t,\Delta_t)\geq \frac{c_1}{n+1}.
\]
Here we use the comparison between $\alpha$-invariants and stability thresholds from \cite[Theorem A]{BJ17}. By boundedness results due to Jiang \cite{Jia17}, Chen \cite{Che18} and Li-Liu-Xu \cite[Corollary 6.14]{LLX18}, we can take a multiple $m$ of $m_1$ such that $-m(K_{\cY_{t,0}}+\Gamma_{t,0})$ is very ample for any $t\in T$. Thus $-m(K_{\cY_t/\bA^1}+\Gamma_t)$ is also very ample over $\bA^1$ for any $t\in T$. By Kawamata-Viehweg vanishing theorem, we know that $N+1=h^0(\cY_{t,0},\cO_{\cY_{t,0}}(-m(K_{\cY_{t,0}}+\Gamma_{t,0})))$. Hence there exists a $\bG_m$-equivariant embedding $(\cY_t,\Gamma_t)\hookrightarrow\bP^N\times\bA^1$. As a result, there exists a $1$-PS $\lambda$ in $\SL(N+1)$ whose induced test configuration of $(\cX_t,\Delta_t)$ is exactly $(\cY_t,\Gamma_t)$. Since $(\cX_t,\Delta_t)$ is $\mu$-twisted K-unstable, we have $\mu>\delta_t$ by Theorem \ref{thm:twisted=delta}. Since $(\cY_t,\Gamma_t)$ is non-trivial, we have
\[
\Fut_{1-\mu}(\cY_t,\Gamma_t)<\Fut_{1-\delta_t}(\cY_t,\Gamma_t)=0.
\]
This finishes the proof of (2).
\end{proof}

Next we use Hilbert schemes and Chow schemes to study twisted K-stability. We first introduce our parameter space. We follow notation from \cite[Chapter I]{Kol96}.

\begin{defn}
Let $\chi(k):=\chi(\cX_t,\cO_{\cX_t}(k))$ be the Hilbert polynomial of $\cX_t\hookrightarrow\bP^N$ for $t\in T$. Denote by $\bfH:=\Hilb_{\chi}(\bP^N)$ the Hilbert scheme of $\bP^N$ with Hilbert polynomial $\chi$. Let $\cX_{\bfH}\hookrightarrow \bP^N\times\bfH$ be the universal family.
Let $\bfH^*$ be the seminormalization of the reduced scheme $\bfH_{\red}$ supported on $\bfH$. Denote by $\cX_{\bfH^*}:=\cX_{\bfH}\times_{\bfH}\bfH^*$.

Suppose $\Delta_T=\sum_{i=1}^{l} c_i\Delta_{T,i}$ where $\Delta_{T,i}$ are distinct prime divisors and $c_i\in\bQ_{>0}$. Then each $\Delta_{t,i}$ is a Weil divisor on $\cX_t$. Denote by $d_i$ (resp. $d_0$) the degree of $\Delta_{t,i}$ (resp. $\cX_t$) as a cycle in $\bP^N$. 
Consider the Chow schemes 
\[
W:=\Chow_{n-1,d_1}(\cX_{\bfH^*}/\bfH^*)\times_{\bfH^*}\cdots\times_{\bfH^*}\Chow_{n-1,d_l}(\cX_{\bfH^*}/\bfH^*)
\]
and
\[
Z:=\Chow_{n-1,d_0}(\cX_{\bfH^*}/\bfH^*)\times_{\bfH^*} W.
\]
Hence we have universal families $\cX_W\to W$ and $\cX_Z\to Z$ over projective seminormal schemes $W$ and $Z$ by pulling back $\cX_{\bfH^*}\to \bfH^*$, together with   $\Delta_{W,i}/W$ and $\Delta_{Z,i}/Z$ as universal families of $(n-1)$-cycles on $\cX_{\bfH^*}/\bfH^*$ for $1\leq i\leq l$, respectively. 
Moreover, there is one more universal family $\cH_Z/Z$ of $(n-1)$-cycles on $\cX_{\bfH^*}/\bfH^*$. Denote by $\Delta_Z:=\sum_{i=1}^l c_i\Delta_{Z,i}$.
\end{defn}

\begin{prop}
Let $P$ be the $\PGL(N+1)$-torsor over $T$ associated to the vector bundle $\pi_*\cO_{\cX}(-m(K_{\cX/T}+\Delta_T))$.
Then we have the following commutative diagram 
\[
 \begin{tikzcd}
   (\bP^N)^{\vee}\times P\arrow{d}{\tPhi}\arrow{r}{} & P \arrow{d}{\Phi}\arrow{r}  &  T\arrow{d}{\varphi}\\
   Z \arrow{r}{} & W\arrow{r} & {[W/\PGL(N+1)]}
        \end{tikzcd}
 \]
\end{prop}

\begin{proof}
Let $(\cX_P,\Delta_P):=(\cX,\Delta_T)\times_T P$.
From the definition of $P$, it is clear that there is a $\PGL(N+1)$-equivariant embedding $\iota: (\cX_P,\Delta_P)\hookrightarrow \bP^N\times P $. Hence the morphism $\Phi$ is defined by the universality of $W$ as a Hilbert-Chow scheme. To define $\tPhi$, note that each point $h\in(\bP^N)^\vee$ corresponds to a hyperplane $H_h$ of $\bP^N$. Since $\iota(\cX_p)$ is a  non-degenerate subvariety of $\bP^N$ for any $p\in P$ by Proposition \ref{prop:bdd-test-twistedK}, we know that $H_h|_{\iota(\cX_p)}$ is an effective Cartier divisor of degree $d_0$. Hence $\tPhi$ is defined by the universality of $Z$ as a Hilbert-Chow scheme. It is quite obvious that both $\Phi$ and $\tPhi$ are $\PGL(N+1)$-equivariant morphisms. Thus the proof is finished.
\end{proof}

\begin{proof}[Proof of Theorem \ref{thm:odc-openness}]
Our strategy is similar to \cite[Proof of Proposition A.4]{LWX14}.
Let $\pi_{Z}:(\cX_Z,\Delta_Z+\cH_Z)\to Z$ be the universal family over $Z$. Denote by $\cL$ the line bundle $\cO_{\cX_Z}(1)$. Consider the following $\bQ$-line bundles coming from the CM line bundle and the Chow line bundle (see \cite[Section 2.4]{ADL19}):
\[
M_1:=\lambda_{\CM,\pi_Z,\Delta_Z,\cL}\quad \textrm{ and }\quad
M_2:=\frac{n(\cL_z^{n-1}\cdot\cH_z)}{(\cL_z^n)}\lambda_{\Chow,\pi_Z,\cL}-(n+1)\lambda_{\Chow, \pi_Z|_{\cH_Z},\cL|_{\cH_Z}}.
\]
Then we see that $\lambda_{\CM,\pi_Z,\Delta_Z+(1-\mu)\cH_Z,\cL}=M_1-(1-\mu)M_2$. Let $Z^\circ:=\{z\in Z\mid \cX_z\textrm{ is normal}\}$. Then it is clear that $Z^\circ$ is a Zariski open subset of $Z$ and 
$\tPhi((\bP^N)^\vee\times P)\subset Z^\circ$. Denote by $G:=\SL(N+1)$ for simplicity. Then for any $z\in Z^\circ$ and any $1$-PS $\lambda$ of $G$, we have
\begin{equation}\label{eq:flagcx1}
\Fut(\cX_{z,\lambda},\Delta_{z,\lambda}+(1-\mu)\cH_{z,\lambda};\cL_{z,\lambda})=\frac{1}{(n+1)(\cL_z^n)}\left(\upmu^{M_1}(z,\lambda)-(1-\mu)\upmu^{M_2}(z,\lambda)\right).
\end{equation}
Here $\upmu^*(*,*)$ represents the GIT weight (see \cite[Definition 2.2]{GIT}), and $(\cX_{z,\lambda},\Delta_{z,\lambda}+(1-\mu)\cH_{z,\lambda};\cL_{z,\lambda})$ denotes the test configuration of $(\cX_z,\Delta_z+(1-\mu)\cH_z, \cL_z)$ induced by $\lambda$ (see e.g. \cite[Section 2.3]{BHJ17}).
Let $M$ be a sufficiently ample $\SL(N+1)$-line bundle on $Z$ such that $L_1:=M\otimes M_1$ and $L_2:=M\otimes M_2$ are both ample $\bQ$-line bundles. Then we have
\[
\upmu^{M_1}(z,\lambda)-(1-\mu)\upmu^{M_2}(z,\lambda)=\left(\upmu^{L_1}(z,\lambda)-\upmu^{M}(z,\lambda)\right)-(1-\mu)\left(\upmu^{L_2}(z,\lambda)-\upmu^{M}(z,\lambda)\right).
\]
We fix a maximal torus $\bT\subset G$.
Hence \cite[Lemma 9.3]{LWX14} implies that there exists a decomposition $Z^{\circ}=\sqcup_I Z_I$ to constructible subsets $Z_I$ where $I$ belongs to some finite index set $\fI$, such that  $\upmu^{M_1}(z,\cdot)$ and $\upmu^{M_2}(z,\cdot)$ are continuous rational piecewise linear functions on $\Hom_{\bQ}(\bG_m, \bT)$ for any $z\in Z_I$, and they are independent of the choice of $z\in Z_I$. We denote these two functions by $\upmu_{1,I}(\cdot)$ and $\upmu_{2,I}(\cdot)$, respectively. 

Given a point $t\in T$, let us take $p\in P$ lying over $t$. Then $\Phi(p)\in W$ is the Hilbert-Chow point of an embedding of $\iota_p:(\cX_t, (\Delta_{t})_{\red}=\sum_{i=1}^l \Delta_{t_i})\hookrightarrow\bP^N$. Each point $h\in (\bP^N)^\vee$ corresponds to a hyperplane $H_h$ in $\bP^N$. Let $z:=\tPhi(h,p)\in Z^{\circ}$. Then it is clear that $(\cX_z,\Delta_z,\cH_z)= (\iota_{p}(\cX_t,\Delta_t), H_h|_{\iota_{p}(\cX_t)})$. Thus by Proposition \ref{prop:bdd-test-twistedK}, we know that $(\cX_t,\Delta_t)$ is $\mu$-twisted K-semistable if and only if 
\begin{equation}\label{eq:flagcx2}
\Fut_{1-\mu}(\cX_{z,\lambda},\Delta_{z,\lambda})=\max_{h\in (\bP^N)^\vee}\Fut(\cX_{z,\lambda},\Delta_{z,\lambda}+(1-\mu)\cH_z)\geq 0
\end{equation} for any $1$-PS $\lambda$ of $G$. Combining \eqref{eq:flagcx1} and \eqref{eq:flagcx2} yields that $(\cX_t,\Delta_t)$ is $\mu$-twisted K-semistable if and only if 
\[
\max_{h\in(\bP^{N})^\vee}\big(\upmu^{M_1}(z,\lambda)- (1-\mu)\upmu^{M_2}(z,\lambda)\big)\geq 0 \quad\textrm{ for any $1$-PS $\lambda$ of $G$.}
\]
Since any $1$-PS $\lambda$ of $G$ is conjugated to a $1$-PS $g\lambda g^{-1}$ of $\bT$ for some $g\in G$, and $\upmu^{M_i}(z, \lambda)=\upmu^{M_i}(g\cdot z, g\lambda g^{-1})$ for $i=1,2$, we know that $(\cX_t,\Delta_t)$ is $\mu$-twisted K-semistable if and only if for any $1$-PS $\lambda$ in $\bT$ and any $g\in G$, there exists $h\in (\bP^N)^\vee$ such that
\begin{equation}\label{eq:flagcx3}
\upmu^{M_1}(g\cdot z,\lambda)\geq (1-\mu)\upmu^{M_2}(g\cdot z,\lambda) \quad \textrm{ where $z=\tPhi(h,p)$}.
\end{equation}
For each $I\in\fI$, define $C_I:=\{\lambda\in\Hom_{\bQ}(\bG_m,\bT)\mid \upmu_{2,I}(\lambda)<0\}$. Since $\upmu_{2,I}$ is continuous rational piesewise linear, we know $C_I$ is a finite union of open rational polyhedral cones in $\Hom_{\bQ}(\bG_m,\bT)$. 
We also define the function $\mu_{\max,I}:\Hom_{\bQ}(\bG_m,\bT)\to (-\infty,+\infty]$ as
\[
\mu_{\max,I}(\lambda):=\begin{cases}
1-\frac{\upmu_{1,I}(\lambda)}{\upmu_{2,I}(\lambda)} & \textrm{if }\lambda\in C_I\\
+\infty & \textrm{otherwise}
\end{cases}
\]
For each $p\in P$, we assign a subset $\fI_p$ of $\fI$ as
\[
\fI_p:=\{I\in\fI\mid \textrm{For any $g\in G$, there exists $h\in (\bP^N)^{\vee}$ such that $g\cdot\tPhi(h,p)\in Z_I$}\}.
\]
We claim that $(\cX_t,\Delta_t)$ is $\mu$-K-semistable if and only if
\begin{equation}\label{eq:flagcx4}
    \mu\leq \inf_{\lambda\in \Hom_{\bQ}(\bG_m, \bT)}\max_{I\in\fI_p} \mu_{\max,I}(\lambda).
\end{equation}
We choose $c_1>0$ as in the proof of Proposition \ref{prop:bdd-test-twistedK} such that $\delta(\cX_t,\Delta_t)>c_1$ for any $t\in T$. Then Theorem \ref{thm:twisted=delta} implies that $(\cX_t,\Delta_t)$ is $\mu$-K-semistable for any $t\in T$ and any $\mu\in (0,c_1]$. Hence \eqref{eq:flagcx3} is satisfied for $\mu\in (0,c_1]$, which means that the restrictions of $\mu$ from \eqref{eq:flagcx3} when $\mu>c_1$ only come from the case when $\mu^{M_2}(g\cdot z,\lambda)<0$. This is precisely addressed in the definitions of $C_I$, $\fI_p$, and $\mu_{\max,I}(\cdot)$.

Next we show that $p\mapsto \fI_p$ is a constructible function from $P$ to $2^\fI$. Define $\tPhi_G:G\times(\bP^N)^\vee\times P\to Z^\circ$
as $\tPhi_G(g, h, p):=g\cdot \tPhi(h,p)$. Let $\pr_{G\times P}:G\times (\bP^N)^\vee\times P\to G\times P$ and $\pr_P:G\times P\to P$ be projection morphisms.
Given $I\in\fI$, from \eqref{eq:flagcx4} we know that 
\begin{align*}
 I\in \fI_p &\Longleftrightarrow \textrm{ for any $g\in G$ there exists $h\in(\bP^N)^\vee$ such that }(g,h,p)\in\tPhi_G^{-1}(Z_I) \\
 & \Longleftrightarrow G\times\{p\}\subset\pr_{G\times P}(\tPhi_G^{-1}(Z_I))\\
 & \Longleftrightarrow p\not\in\pr_P\big((G\times P)\setminus\pr_{G\times P}(\tPhi_G^{-1}(Z_I))\big)
\end{align*}
Since $Z_I$ is a constructible subset of $Z^\circ$, by Chevalley's lemma \cite[Exercise II.3.19]{Har77} we know that  $\pr_P\big((G\times P)\setminus\pr_{G\times P}(\tPhi_G^{-1}(Z_I))\big)$ is also a constructible subset of $P$. Hence $p\mapsto \fI_p$ is a constructible function from $P$ to $2^\fI$. 

If $(\cX_t,\Delta_t)$ is a K-unstable fiber satisfying Conjecture \ref{conj:odc}, then from Theorem \ref{thm:delta-min-preserve} we know that $\delta(\cX_t,\Delta_t)$ is a rational number.
Since stability thresholds are lower semicontinuous in families \cite{BL18b}, to prove the theorem it suffices to show that $t\mapsto \min\{1,\delta(\cX_t,\Delta_t)\}$ is constructible. From the discussion above, we have that 
\[
\min\{1,\delta(\cX_t,\Delta_t)\}=\min\left\{1, \inf_{\lambda\in \Hom_{\bQ}(\bG_m, \bT)}\max_{I\in\fI_p} \mu_{\max,I}(\lambda)\right\},
\]
where the function $p\mapsto \fI_p$ is constructible. Hence $p\mapsto \min\{1,\delta(\cX_t,\Delta_t)\}$ is a $\PGL(N+1)$-invariant constructible function on $P$. Since $P$ is a $\PGL(N+1)$-torsor over $T$, we know that $t\mapsto \min\{1,\delta(\cX_t,\Delta_t)\}$ is also constructible. This finishes the proof.
\end{proof}

Next we recall the definition of the greatest Ricci lower bound of a Fano manifold.

\begin{defn}
The \emph{greatest Ricci lower bound} of a Fano manifold $X$ is defined as
\[
\beta(X):=\sup\{\beta\in [0,1]\mid \textrm{there exists a K\"ahler metric $\omega\in c_1(X)$ such that $\mathrm{Ric}(\omega)>\beta\omega$}\}.
\]
\end{defn}

This invariant was studied by Tian in \cite{Tia92}, although it was not explicitly defined there. It was first
explicitly defined in \cite{Rub08, Rub09} and further studied by Sz\'ekelyhidi \cite{Sze11}.
Recently, it was proven by Berman, Boucksom, and Jonsson \cite{BBJ18} and independently by Cheltsov, Rubinstein, and Zhang \cite{CRZ19} that $\beta(X)=\min\{1,\delta(X)\}$. Note that $\beta(X)$ is also denoted by $R(X)$ in some papers.

The following result  confirms Conjecture \ref{conj:odc} for Fano manifolds. It essentially follows from the work of Datar and Sz\'ekelyhidi \cite{DS16}, Ross and Sz\'ekelyhidi \cite{RS19}, and the relation between $\beta(X)$ and $\delta(X)$ from \cite{BBJ18, CRZ19}.

\begin{thm}\label{thm:odc-Fano-manifold}
 Let $X$ be a Fano manifold that is not uniformly K-stable. Then there exists a prime divisor $F$ over $X$ computing the stability threshold $\delta(X)$.
\end{thm}

\begin{proof}
Consider the anti-pluricanonical embedding $X\hookrightarrow\bP^N$ given by the linear system $|-mK_X|$ for $m$ sufficiently large.
Choose a K\"ahler metric $\sigma:=\frac{1}{m}\omega_{\rm FS}|_{X}\in c_1(X)$. By \cite{Sze11}, for any $\beta\in (0, \beta(X))$ there exists a K\"ahler metric $\omega_\beta\in c_1(X)$ on $X$ such that $\mathrm{Ric}(\omega_\beta)=\beta\omega_\beta+(1-\beta)\sigma$. 
In \cite[Section 3]{DS16}, through taking the Gromov-Hausdorff limit of $(X,\omega_\beta)$ as $\beta\nearrow\beta(X)$, the authors constructed a sequence of non-trivial special test configurations $(\cX_i;\cL_i)\to \bA^1$ of $X$ such that $\lim_{i\to\infty}\Fut_{1-\beta(X)}(\cX_i;\cL_i)/||(\cX_i;\cL_i)||_{L^2}=0$ and each $(\cX_i;\cL_i)$ is induced by a $1$-PS $\lambda_i$ of $\mathrm{SL}(N+1)$.
Then \cite[Corollary 9]{RS19} implies that $X$ is not $\beta(X)$-twisted K-stable.
(Note that $\Fut_{1-\beta}(\cX_i, \cL_i)$ equals  $\Fut_{(1-\beta) \sigma } (\cX_i, \cL_i)$ in the notation of Ross-Sz\'ekelyhidi's; see  \cite[page 3]{RS19}.)
Hence,   there exists a non-trivial normal test configuration $(\cX;\cL)$ induced by a $1$-PS $\lambda$ of $\SL(N+1)$ such that  $\Fut_{1-\beta(X)}(\cX;\cL)=0$. 
Since $\beta(X)=\delta(X)$ by \cite{BBJ18, CRZ19},
Theorem \ref{thm:twisted-iff} implies $v_{\cX_0}$ is a  divisorial valuation computing $\delta(X)$. Hence, the proof is finished.
\end{proof}

\begin{remark}\label{rem:twisted-Kps}
In the above proof, it is natural to expect that the optimal degeneration $\cX_0$ of $X$ constructed in \cite{DS16} is independent of the choice of $\sigma$. Indeed, it suffices to show that  $\cX_0$ is $\beta(X)$-twisted K-polystable since then Theorem \ref{thm:kpolydegen} implies that $\cX_0$ is uniquely determined. 
\end{remark}

\begin{proof}[Proof of Theorem \ref{thm:greatestRic}]
Since  Fano manifolds of a fixed dimension $n$ form a bounded family by \cite{KMM92, Cam92}, we can find a smooth morphism $\pi:\cX\to T$ parametrizing all of them. By Theorem \ref{thm:odc-Fano-manifold}, we know that Conjecture \ref{conj:odc} holds for any K-unstable fiber of $\pi$. Thus the theorem directly follows by applying Theorem \ref{thm:odc-openness} to the family $\pi:\cX\to T$.
\end{proof}

\bibliography{reference}

\end{document}